\documentclass[preprint,3p]{elsarticle}

\usepackage{graphicx}
\usepackage{subfigure}
\usepackage{dsfont}
\usepackage[inkscapeformat=png]{svg}
\usepackage{hyphenat}

\usepackage{amsmath}
\usepackage{tikz}

\usepackage{siunitx}

\usepackage{xcolor}
\definecolor{tba}{RGB}{1,164,154}
\definecolor{bw}{RGB}{237,102,0}
\definecolor{rdi}{RGB}{145,0,237}
\definecolor{red}{RGB}{255,0,0}
\definecolor{tum}{RGB}{48,112,179}

\renewcommand{\vec}[1]{\boldsymbol{#1}}
\newcommand{\tens}[1]{\boldsymbol{#1}}

\newcommand{\R}{\mathds{R}}

\newcommand{\dt}[1]{\dot{#1}}							
\newcommand{\transpose}[1]{{#1}^T}						
\newcommand{\transposeminus}[1]{{#1}^{-T}}				
\newcommand{\inv}[1]{#1^{-1}}							

\newcommand{\bilinear}[3]{\klr{ {#1} \text{,} {#2} }_ {{#3}}}             

\newcommand{\kl}[1]{\left( #1 \right)}
\newcommand{\klr}[1]{\left( #1 \right)}
\newcommand{\kle}[1]{\left[\ #1\right]}

\newcommand{\mytime}{t}
\newcommand{\infinetismald}{\mathrm{d}}

\newcommand{\partialfrac}[2]{\frac{\partial {#1}}{\partial {#2}}}

\newcommand{\grad}{\boldsymbol{\nabla}}
\newcommand{\matgrad}{\grad_0 }

\newcommand{\nsd}{n_{sd}}

\newcommand{\fluid}{f}									
\newcommand{\solid}{s}									

\newcommand{\dis}{\boldsymbol{d}}					    
\newcommand{\us}{\dis^\solid}							

\newcommand{\vel}{\boldsymbol{v}}						
\newcommand{\vf}{\vel^\fluid}							
\newcommand{\vs}{\vel^\solid}							

\newcommand{\hatdis}{\hat{\dis}}						
\newcommand{\hatvf}{\hat{\vel}^\fluid}							
\newcommand{\hatvs}{\hat{\vel}^\solid}							

\newcommand{\initial}{0}
\newcommand{\vfini}{\hatvf_{\initial}}				    
\newcommand{\dsini}{\hatdis^\solid_{\initial}}					

\newcommand{\tstart}{\mytime_\initial}                         
\newcommand{\tfinal}{\mytime_e}                         
\newcommand{\timeinterval}{\tstart,\tfinal}       

\newcommand{\p}{p}										
\newcommand{\pf}{\p^\fluid}								

\newcommand{\den}{\rho}									
\newcommand{\rhof}{\den^\fluid}							
\newcommand{\rhos}{\den^\solid}							

\newcommand{\xref}{\boldsymbol{X}}
\newcommand{\xcurr}{\boldsymbol{x}}

\newcommand{\defgrad}{\boldsymbol{F}}					          
\newcommand{\detdefgrad}{J}								

\newcommand{\GreenLagrangestrain}{\boldsymbol{E}}		             

\newcommand{\strainenergysolid}{\Psi^\solid}
\newcommand{\strainenergysolidEJporosity}{\strainenergysolid (\GreenLagrangestrain, \detdefgrad,\porosity)}

\newcommand{\macrostrainskel}{\Psi^{skel}\klr{\GreenLagrangestrain}}
\newcommand{\macrostrainvol}{\Psi^{vol}\klr{\detdefgrad,\porosity}}
\newcommand{\macrostrainpen}{\Psi^{pen}\klr{\GreenLagrangestrain, \detdefgrad,\porosity}}

\newcommand{\bulkmod}{\kappa}
\newcommand{\poromatpenalty}{\eta}

\newcommand{\Pktwostress}{\boldsymbol{S}}

\newcommand{\fluidstressvisc}{\boldsymbol{\sigma}^\fluid_{visc}}
\newcommand{\Pktwofluidvisc}{\boldsymbol{S}^\fluid_{visc}}

\newcommand{\rightcauchgreendeftensor}{\boldsymbol{C}}

\newcommand{\bodyforce}{\boldsymbol{b}} 				
\newcommand{\bodyforcefluid}{\bodyforce^\fluid}		    
\newcommand{\bodyforcesolid}{\bodyforce^\solid}		    

\newcommand{\tractionfluid}{\hat{\boldsymbol{t}}^\fluid}
\newcommand{\tractionfluidneuman}{\hat{\boldsymbol{t}}^\fluid}

\newcommand{\tractionsolid}{\hat{\boldsymbol{t}}^\solid}

\newcommand{\porosity}{\phi}							
\newcommand{\porosityinitial}{\phi_\initial}					
\newcommand{\hatporosityinitial}{\hat{\porosity}_\initial}					

\newcommand{\vis}{\mu}									
\newcommand{\visfluid}{\mu^\fluid}						

\newcommand{\permea}{\boldsymbol{k}}						
\newcommand{\permeainitial}{\boldsymbol{K}}	

\newcommand{\virtual}{\delta}							

\newcommand{\virtualvf}{\delta \boldsymbol{w}^{\fluid}}                   
\newcommand{\virtualpf}{\delta \pf}
\newcommand{\qf}{q^{\fluid}}
\newcommand{\virtualus}{\delta \us}
\newcommand{\virtualE}{\delta \GreenLagrangestrain}
\newcommand{\virtualporosity}{\delta \porosity}

\newcommand{\discrete}{h}								

\newcommand{\discreteresidualPFM}{\mathcal{R}^\discrete_\mathcal{PFM}}

\newcommand{\nele}{n_{ele}}                             

\newcommand{\domain}{\Omega}							
\newcommand{\porodomain}{\Omega^p}						
\newcommand{\domaininitial}{\Omega_\initial}					
\newcommand{\domaindeformed}{\Omega_t}					
\newcommand{\domaindeformedfluid}{\Omega_t^{\fluid}}			

\newcommand{\fluiddomain}{\Omega^{\fluid}}				
\newcommand{\soliddomain}{\Omega^{\solid}}              

\newcommand{\bcsurfdiri}{\Gamma_{D}}
\newcommand{\bcsurfneu}{\Gamma_{N}}

\newcommand{\boundarydirichtime}{\Gamma^{D}_\mytime}

\newcommand{\boundaryconstraint}{\Gamma^{C}}
\newcommand{\boundaryconstrainttime}{\Gamma^{C}_\mytime}

\newcommand{\boundaryfluiddirichtime}{\Gamma^{D\fluid}_\mytime}
\newcommand{\boundarysoliddirichtime}{\Gamma^{D\solid}_\mytime}

\newcommand{\boundaryfluidneumantime}{\Gamma^{N\fluid}_\mytime}
\newcommand{\boundarysolidneumantime}{\Gamma^{N\solid}_\mytime}

\newcommand{\boundarysoliddirichinitial}{\Gamma^{D\solid}_\initial}

\newcommand{\boundaryfluidneumaninitial}{\Gamma^{N\fluid}_\initial}
\newcommand{\boundarysolidneumaninitial}{\Gamma^{N\solid}_\initial}

\newcommand{\boundaryporositydirichinitial}{\Gamma^{D\porosity}_\initial}

\newcommand{\normalvectorcurr}{\boldsymbol{n}}
\newcommand{\normalvectorinitial}{\boldsymbol{N}}

\newcommand{\Identity}{\boldsymbol{I}}
\newcommand{\vecnull}{\boldsymbol{0}}

\newcommand{\pcttrans}{\dt{\porosity}}

\newcommand{\pctsolid}{\porosity \boldsymbol{\nabla} \cdot \vs }

\newcommand{\pctfluid}{\boldsymbol{\nabla} \cdot \klr{\porosity \klr{\vf-\vs} }}

\newcommand{\fpmttrans}{\rhof \dt{\vf}}

\newcommand{\fpmtconv}{\rhof \kl{\kl{\vf -\vs}\cdot \boldsymbol{\nabla} }\vf}

\newcommand{\fpmtpress}{\boldsymbol{\nabla} \pf}

\newcommand{\fpmtbody}{\rhof \bodyforcefluid}

\newcommand{\fpmtdarcy}{\visfluid \porosity \inv{\permea} \cdot \kl{\vf -\vs}}

\newcommand{\fpmtvisc}{\frac{1}{\porosity} \boldsymbol{\nabla} \cdot \left(\phi \boldsymbol{\sigma}_{{visc}}^{\fluid}\right)}

\newcommand{\sigmafvis}{\sigma_{visc}^\fluid}
\newcommand{\Pktwostressterm}{
\porosity\Pktwofluidvisc -\pf\detdefgrad\inv{\rightcauchgreendeftensor}+\frac{\partial \strainenergysolidEJporosity}{\partial \GreenLagrangestrain}}

\newcommand{\defgradterm}{\frac{\partial \xcurr}{\partial \xref}}

\newcommand{\greenlagrangeterm}{\frac{1}{2}\klr{\transpose{\defgrad}\defgrad-\boldsymbol{1}}}

\newcommand{\spmttrans}{\rhos_0(1-\porosityinitial) \dt{\vs}}

\newcommand{\spmtstress}{\matgrad \cdot \klr{\defgrad \Pktwostress } }

\newcommand{\spmtbody}{\rhos_0(1-\porosityinitial) \bodyforce}

\newcommand{\spmtpress}{\detdefgrad\porosity \transposeminus{\defgrad} \cdot \matgrad \pf}

\newcommand{\spmtdarcy}{\visfluid \detdefgrad \porosity^2 \inv{\permea} \cdot \kl{\vf -\vs}}

\newcommand{\macrostrainvolterm}{\bulkmod \klr{\frac{ \detdefgrad \klr{1-\porosity}}{1-\porosityinitial} - 1 - ln \klr{ \frac{ \detdefgrad \klr{1-\porosity} }{1-\porosityinitial}}}}

\newcommand{\macrostrainpenterm}{\poromatpenalty\klr{-ln\klr{\frac{\detdefgrad\porosity}{\porosityinitial}}+\frac{\detdefgrad \porosity}{\porosityinitial}-\frac{1}{\porosityinitial}}}

\newcommand{\taumstab}{\tau_m}
\newcommand{\PSPGterm}{\sum_{e=1}^{\nele}  \bilinear{\grad {\virtualpf} }{\taumstab \discreteresidualPFM }{\Omega_e} } 

\newcommand{\PSPGDarcyterm}{\sum_{e=1}^{\nele}  \bilinear{ \porosity \visfluid\inv{\permea}  {\virtualvf} }{-\taumstab \discreteresidualPFM}{\Omega_e} }

\newcommand{\pushforward}[1]{\frac{1}{\detdefgrad}\defgrad\cdot {#1}\cdot \transpose{\defgrad}}

\usepackage[acronyms,nopostdot,nogroupskip,section=section,shortcuts]{glossaries}

\newacronym{ac:web}{WEB}{Woven Endo Bridge}
\newcommand{\web}{\gls{ac:web}}
\newcommand{\webs}{\glspl{ac:web}}

\newacronym{ac:fd}{FD}{Flow Diverter}
\newcommand{\FD}{\gls{ac:fd}}
\newcommand{\FDs}{\glspl{ac:fd}}

\newacronym{ac:ca}{CA}{Cerebral Aneurysm}
\newcommand{\CA}{\gls{ac:ca}}
\newcommand{\CAs}{\glspl{ac:ca}}

\newacronym{ac:lbm}{LBM}{Lattice Boltzmann Method}
\newcommand{\LBM}{\gls{ac:lbm}}
\newcommand{\LBMs}{\glspl{ac:lbm}}

\newacronym{ac:fem}{FEM}{Finite Element Method}
\newcommand{\FEM}{\gls{ac:fem}}

\newacronym[firstplural=porous media (PMs)]{ac:pm}{PM}{Porous Medium}
\newcommand{\PM}{\gls{ac:pm}}
\newcommand{\PMs}{\glspl{ac:pm}}

\newacronym{ac:wss}{WSS}{wall shear stress}
\newcommand{\WSS}{\gls{ac:wss}}

\newacronym{ac:ale}{ALE}{Arbitrary-Lagrangian-Eulerian}

\newacronym{ac:pspg}{PSPG}{Pressure Stabilizing Petrov-Galerkin}
\newcommand{\PSPG}{\gls{ac:pspg}}

\newacronym{ac:supg}{SUPG}{Streamline Upwind Petrov-Galerkin}
\newcommand{\SUPG}{\gls{ac:supg}}
\usepackage{algorithm}
\usepackage{algpseudocode}
 \usepackage[nodots,nocompress]{numcompress}
\usepackage{cleveref}

\newtheorem{Remark}{Theorem}
\journal{arxiv }
\begin{document}
\begin{frontmatter}
\title{Numerical simulation of endovascular treatment options for cerebral aneurysms }

\author[0]{Martin Frank}

\author[1]{Fabian Holzberger}

\author[1]{Medeea Horvat}

\author[2]{Jan Kirschke}

\author[0,3]{Matthias Mayr\corref{cor1}}
\cortext[cor1]{Corresponding author}
\author[1]{Markus Muhr}

\author[1]{Natalia Nebulishvili}

\author[0]{Alexander Popp}

\author[2]{Julian Schwarting}

\author[1]{Barbara Wohlmuth}



\affiliation[0]{organization={University of the Bundeswehr Munich, Institute for Mathematics and Computer-Based Simulation},
	addressline={Werner-Heisenberg-Weg~39}, 
	city={Neubiberg},
	postcode={85577}, 
	country={Germany}}


\affiliation[1]{organization={Technical University of Munich, School of Computation, Information and Technology, Department of Mathematics,Chair for Numerical Mathematics},
	addressline={Boltzmannstr.~3}, 
	city={Garching},
	postcode={85748}, 
	country={Germany}}


\affiliation[2]{organization={Technical University of Munich, Department of Neuroradiology},
	addressline={Ismaninger~Str.~22}, 
	city={Munich},
	postcode={81675}, 
	country={Germany}}

\affiliation[3]{organization={University of the Bundeswehr Munich, Data Science \& Computing Lab},
	addressline={Werner-Heisenberg-Weg~39}, 
	city={Neubiberg},
	postcode={85577}, 
	country={Germany}}
\begin{keyword}Cerebral Aneurysm \sep Endovascular Intervention \sep  Lattice Boltzmann Method \sep  Finite Elements \sep  Porous Medium \sep  Patient-Specific Simulation
\end{keyword}




\begin{abstract}
Predicting the long-term success of endovascular interventions in the clinical management of cerebral aneurysms
requires detailed insight into the patient-specific physiological conditions.
In this work,
we not only propose numerical representations of endovascular medical devices such as coils, flow diverters or Woven EndoBridge
but also outline numerical models for the prediction of blood flow patterns in the aneurysm cavity right after a surgical intervention.
Detailed knowledge about the post-surgical state then lays the basis to assess the chances of a stable occlusion of the aneurysm
required for a long-term treatment success.
To this end,
we propose mathematical and mechanical models of endovascular medical devices made out of thin metal wires.
These can then be used for fully resolved flow simulations of the post-surgical blood flow,
which in this work will be performed by means of a Lattice Boltzmann method
applied to the incompressible Navier-Stokes equations and patient-specific geometries.
To probe the suitability of homogenized models,
we also investigate poro-elastic models to represent such medical devices.
In particular, we examine the validity of this modeling approach for flow diverter placement across the opening of the aneurysm cavity.
For both approaches,
physiologically meaningful boundary conditions are provided from reduced-order models of the vascular system.
The present study demonstrates our capabilities to predict the post-surgical state
and lays a solid foundation to tackle the prediction of thrombus formation and, thus, the aneurysm occlusion in a next step.
\end{abstract}
\end{frontmatter}







\section{Introduction}

{\CAs} are abnormal focal outpouchings of large intracranial arteries that result from a structural weakening of the arterial wall~\cite{Claassen2022a}.
Their rupture leads to a severe type of intracranial hemorrhage, subarachnoid hemorrhage (SAH),
which is a severe subtype of stroke with a high mortality and morbidity resulting in a significant socioeconomic burden~\cite{Hoogmoed2019a,Macdonald2017a,Neifert2021a}.
In the European Union alone, there are an estimated 12.5 million individuals, or about 3\% of the population, who have one or more aneurysms.
Approximately 0.5\% of the world population die from aneurysm ruptures.
Patients, who survive subarachnoid hemorrhages,
experience in 30\% of all cases significant neurological impairments and disabilities~\cite{Pierot2017a}.
The development of aneurysms is influenced by various factors,
including hemodynamic stress, the degradation of arterial walls,
and conditions that result in high blood flow, such as hypertension, smoking, or connective tissue disorders~\cite{Etminan2016a}.

Along with microsurgical clipping of aneurysms through an open craniotomy,
endovascular interventions have emerged as an essential technique for the treatment of aneurysms.
In these minimally invasive procedures,
medical devices are inserted into the cerebral vessels and/or aneurysms through microcatheters
to trigger the occlusion of the aneurysm sac by thrombus formation~\cite{Pierot2013a}.
After the surgical intervention,
either the endothelium re-grows to complete a successful {\CA} occlusion within months after surgery
or inflammatory processes and/or high stresses in the {\CA} walls or compaction of the thrombus can trigger a recanalization,
manifesting a relapse for the patient.

Due to its impact on the homogeneity and stability of the thrombus,
and thereby the risk of aneurysm regrowth,
the choice of medical devices highly influences the overall treatment success.
Despite specialized endovascular procedures and advanced imaging techniques,
device selection is mostly governed by the personal experience of the neuroradiologist.
Furthermore, comparability of devices in specific geometries is limited due to intraindividual variability of the anatomy.
To support the neuroradiologists in their treatment planning and decision process,
computer simulations before device implantation would enable comparisons of different devices
and insertion techniques on a patient-specific basis,
which could substantially support interventionalists in the preparation of the procedure
and in the choice of the device for each patient.
Overall, assisting the neuroradiologist's planning process through numerical simulations
could have a substantial influence on long-term treatment outcomes,
especially in critical scenarios~\cite{Pierot2017a,Sindeev2019a}.


So far, numerical models to analyze the blood flow in {\CAs} are often restricted to computational fluid dynamics (CFD),
which can predict flow patterns in good agreement with post-surgery MRI scans~\cite{Jain2016a,Sindeev2018a}.
Various flow models (e.g. incompressible Navier-Stokes) and discretization techniques (e.g. {\FEM}) are available.
Alternatively, one can reconstruct the velocity and pressure field from the {\LBM} discretization of the Boltzmann equation~\cite{kruger2017lattice}.
In contrast, fluid/structure interaction (FSI) models take the vessel's compliance into account,
even with advanced isotropic or anisotropic constitutive laws~\cite{Tricerri2015a}.
Various endovascular treatment strategies for CAs have been modeled numerically so far,
most often related to either stent placement or coiling.
The novel {\web} device has not been studied in computer simulations, yet.
For example, \cite{Leng2018a} perform a ``virtual'' coiling procedure on a patient-specific CA geometry
using Euler-Bernoulli beam finite elements for the wire model for computational efficiency,
yet assuming a rigid vessel wall and neglecting the blood flow.

Regarding constitutive modeling of blood,
many preliminary studies use the simplifying assumption of a Newtonian fluid.
Yet, capturing the non-Newtonian character is crucial to mimic the shear-dependent viscosity of blood,
for example via the Casson or Carreau-Yasuda model (also applicable to the LBM~\cite{Boyd2007a,horvat2024lattice}).


To provide physiologically meaningful boundary conditions at the inflow and outflow cross sections,
but simultaneously limit the computational effort to the necessary minimum,
one does not represent the entire cardiovascular system with fully resolved three-dimensional models,
but rather resorts to networks of reduced\hyp{} and mixed\hyp{}dimensional models.
Such mixed-dimensional models, i.e. the coupling of fully resolved 3D models with reduced-dimensional 1D and/or 0D models,
allow to represent the effect of the entire vascular system with reasonable computational effort.
While 3D models resolve details of the vessel of interest,
flow within the remaining network of large vessels is represented by a 1D flow model,
i.e. a coupled system of partial differential equations with just a single space variable~\cite{vcanic2003mathematical}.
The remainder of the blood circulation and the heart are represented by ordinary differential equations and algebraic equations,
hence their naming as 0D models.
The lower-dimensional models then supply the boundary conditions for the higher-dimensional models.
Coupled 0D--1D--3D models are well-established in modeling the human vascular system~\cite{ambrosi2012modeling,Perdikaris2016a}.
In closed-loop form, the 3D model is embedded into a reduced-dimensional model of the circulation in the entire vascular system,
as e.g. proposed in~\cite{Zhang2016a} for mixed-dimensional simulation of flow information in the brain's circle of Willis.
In open-loop models, flow quantities are prescribed at the inflow boundary,
while the outflow boundary is coupled to a 0D Windkessel model~\cite{Westerhof2009a}
that needs to be calibrated to patient-specific data~\cite{Ismail2013a}.
The class of 3D--1D models has been successfully used for large arteries~\cite{Formaggia2001a}
as well as microvascular networks~\cite{d2008coupling,Koeppl2020a}.
3D--0D models are widely established in cardiovascular and pulmonary networks~\cite{Ismail2014b},
along with 3D--1D models of flow and transport processes in vascularized tissue~\cite{vidotto2019hybrid,Koeppl2020a}.

In this paper,
we will work towards the numerical modeling of different endovascular treatment options for {\CAs}.
We tackle this topic at several fronts:
On the one hand,
we propose geometric and/or mechanical models of the different endovascular devices, i.e. coiling, flow diverter and {\web}.
For the coiling wire,
a mechanical model will be derived that accounts for the naturally non-straight configuration of the coiling wire
and, thus, represents the self-induced folding into a coil when the coiling wire is deployed from the catheter.
For flow diverters and {\webs}, we show geometric descriptions of their reference configurations and outline a process to map them into patient-specific geometries.
On the other hand,
we consider different continuum models to assess the state of the blood flow in {\CAs} immediately after device placement.
For coiling, we generate a patient-specific coil geometry using our coil model
and then study the reduction in blood flow inside the aneurysm cavity through flow simulations using the {\LBM}.
As an alternative modeling approach,
we represent the medical device and the blood volume by a poro-elastic medium with different porosity and permeability characteristics
and exemplarily study the impact of a flow diverter on the flow field inside the {\CA}.
The results lay the basis for future model improvements as well as for the quantitative prediction of the thrombus formation inside the {\CA}
to forecast the long-term quality and stability of the {\CA} occlusion and, thus, of the treatment success.

The remainder of this paper is organized as follows:
In Section~\ref{sec:MedicalRelevanceAndChallenges},
we discuss different types of {\CAs} as well as their individually appropriate treatment methods and devices from a medical point of view
and list respective challenges and opportunities where numerical simulation can assist in medical decision making.
In Section~\ref{sec:FromImagingToGeometries},
we then introduce our simulation preprocessing pipeline to obtain three-dimensional geometries and computational domains and meshes from the original medical imaging data.
Section~\ref{sec:ContinuumMechanicsAndNumericsOfFlowInAneurysms} proposes to sets of models and numerical methods for the simulation of blood flow through cerebral vessels,
namely the solution of the incompressible Navier-Stokes equations via a {\LBM}
and a homogenized approach using poro-elastic media discretized with the {\FEM}.
In Section~\ref{sec:DeviceModelling}, we discuss concrete mechanical and geometric device models for endovascular coils, {\webs} and stents/flow diverters,
which will then be used together with the porous flow surrogates from Section~\ref{sec:ContinuumMechanicsAndNumericsOfFlowInAneurysms} in numerical experiments in Section~\ref{sec:NumericalSimulations}.
Section~\ref{sec:Conclusion} concludes this manuscript with some final remarks.

\section{Medical Relevance \& Challenges} \label{sec:MedicalRelevanceAndChallenges}

Aneurysms are focal bulges in arteries which can occur anywhere in the body.
Clinically relevant aneurysms can form in the aorta, coronary arteries and vessels supplying the brain. 
Pathophysiologically, these are bulges of the blood vessel wall involving all wall layers
as a result of mostly acquired wall weakness of the basal cerebral arteries caused by high hemodynamic stress, atherosclerosis, or vasculopathies~\cite{Timperman1995a}.
Preferentially, these arise in areas of high hemodynamic stress.
Predilection sites are bifurcations (for example, the division of the anterior cerebral artery and the anterior communicating artery, the division of the internal carotid artery and the posterior communicating artery, the bifurcation of the middle cerebral artery, and the tip of the basilar artery)
and large vessels with a small radius of curvature (for example, in the region of the carotid siphon). 
In 90\% of cases, the anterior circulation of the \emph{circulus arteriosus Willisii (circle of Willis)}
in the stromal area of the internal carotid artery is affected~\cite{Steiner2013a}.

Aneurysm therapy started in the late 1930s with predominantly transcranial surgeries,
where silver clips were first used to eliminate aneurysms~\cite{Dandy1938a}.
Since then, clipping, i.e. the placing of titanium clips at the aneurysm neck,
has become a standard treatment for the occlusion of aneurysms~\cite{Zhao2018a}.
In parallel, the first endovascular approaches to aneurysm occlusion began as early as the 1960s,
and were substantially improved, for instance by Guglielmi \emph{et al.}~\cite{Guglielmi1991a,Guglielmi1991b},
who developed the first platinum coils~\cite{Guglielmi2007a,Pierot2017a}.
The treatment modality of aneurysms is nowadays dependent on several factors,
such as aneurysm location, size, and configuration,
and is recommended individually for every patient by interdisciplinary boards~\cite{Pierot2013a}.
Due to recent innovations in the development of neurointerventional devices, however,
endovascular treatment can nowadays be offered to more and more patients~\cite{Pierot2017a}.

\begin{figure}
\centering
\subfigure[Aneurysm in the anterior circulation before coiling procedure ]{\label{fig:CAAngioCoiling}\includegraphics[height=4.5cm,width=4.5cm]{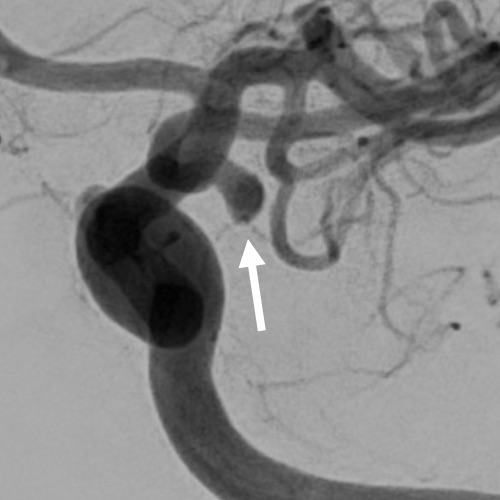}}\hspace{0.5cm}
\subfigure[Aneurysm in the internal carotid artery before \acrfull{ac:fd} implantation]{\label{fig:CAAngioFlowDiverter}\includegraphics[height=4.5cm,width=4.5cm]{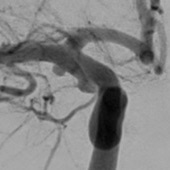}}\hspace{0.5cm}
\subfigure[Aneurysm before {\web} device implantation]{\label{fig:CAAngioWEB}
\includegraphics[height=4.5cm,width=4.5cm]{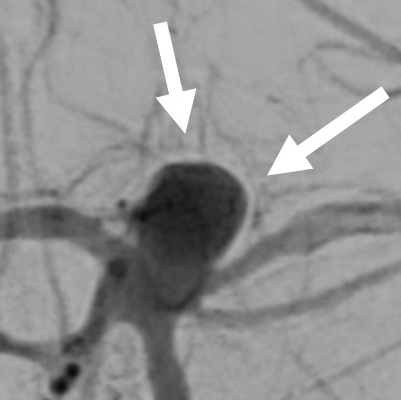}}\\
\subfigure[Placement of five coils inside a non-ruptured {\CA} ]{\label{fig:DeviceAngioCoiling}\includegraphics[height=4.5cm,width=4.5cm]{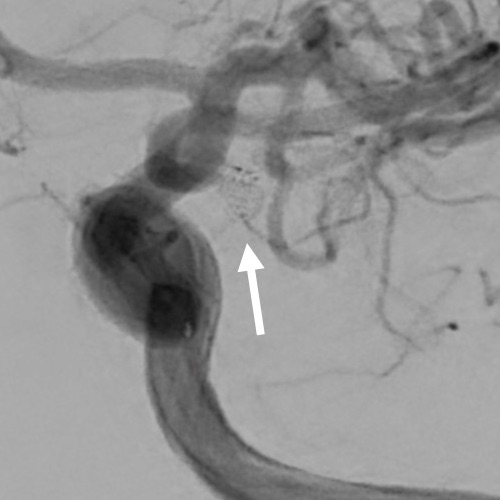}}\hspace{0.5cm}
\subfigure[Enlarged internal carotid artery after {\FD} implantation ({\FD} highligted for better visibility)]{\label{fig:DeviceAngioFlowDiverter}
\includegraphics[height=4.5cm,width=4.5cm]{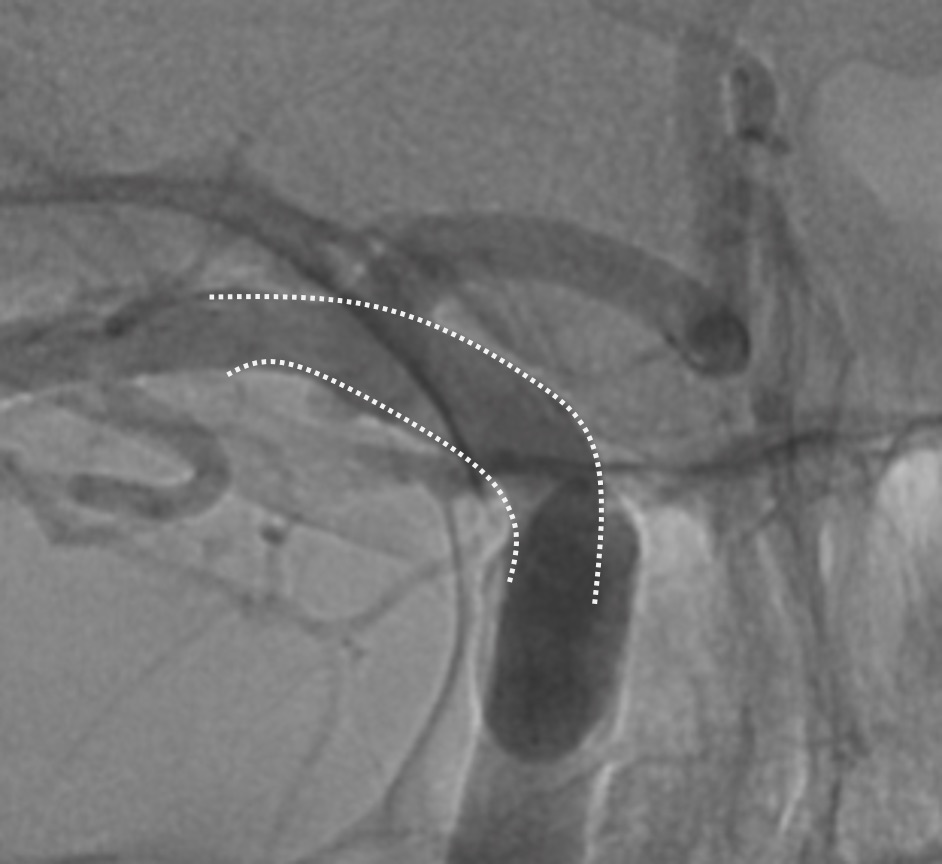}}\hspace{0.5cm}
\subfigure[{\web} implantation in a basilar artery aneurysm]{\label{fig:DeviceAngioWEB}
\includegraphics[height=4.5cm,width=4.5cm]{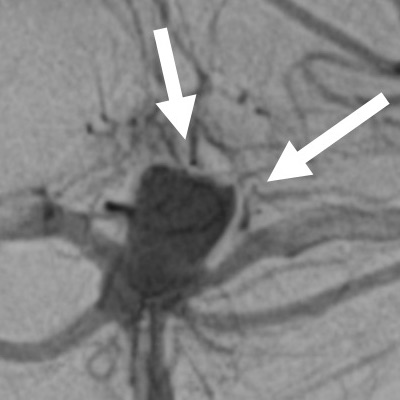}}
\caption{Pre- and post-operative angiography for different options of endovascular aneurysm occlusion with injected contrast medium (Images taken from our data collection described in Section~\ref{subsec:MedicalSelectionAndSegmentation})}
\label{fig:DeviceAngio}
\end{figure}

In principle, there are three different options for endovascular aneurysm occlusion:
coiling, flow diversion, and flow disruption (see Figure~\ref{fig:DeviceAngio}).

\subsection{Coiling}

Intraaneurysmal occlusion by detachable platinum spirals, i.e. coiling (see Figure~\ref{fig:DeviceAngioCoiling}),
was the first method of endovascular aneurysm occlusion and remains up to day the most commonly used technique,
because it enabled for the first time to insert, extract, or disconnect an embolic substance within the cerebral blood vessels,
significantly enhancing the safety of neuroradiological embolization procedures~\cite{Pierot2017a}.
Nonetheless, coiling has a major drawback in comparison to clipping,
which is the recurrence of aneurysms caused by either coil compression or aneurysm expansion~\cite{Pierot2013a}.
This has also been shown in a systematic meta-analysis by Lecler \emph{et al.}~\cite{lv2020review},
who describe in approximately 3000 cases an aneurysm recurrence in 11.4\%, see also~\cite{Pierot2017a}.

The optimal choice of coil material is therefore essential to prevent aneurysm recurrence.
Additionally, hydrogel-coated coils were developed
which can expand after coil implantation and fill more of the aneurysm lumen,
which reduced rates of aneurysm recurrence in recent studies.
A particular challenge for coiling is the treatment of wide-neck aneurysms
as coil loops can extend into the parent artery more commonly.
Also, the risk of aneurysm recurrence is higher in these aneurysm shapes.
In some cases, coils can therefore be combined with stents or balloons,
which can enable the use of coils in these cases.
Ultimately, devices for flow diversion and flow disruption were developed to address this challenge~\cite{Pierot2017a}.

\subsection{Flow diversion}
\label{sec:FlowDiversion}

\acrfullpl{ac:fd} represent a new generation of stent-like devices that are placed along the wall of the parent vessel (see Figure~\ref{fig:DeviceAngioFlowDiverter}).
{\FDs} alter the course of blood flow into the parent vessel,
leading to reduced blood movement within the aneurysm and the formation of a blood clot,
which eventually transforms into scar tissue.
Over time, the aneurysm's neck may also become covered with a new layer of tissue,
known as neointimal proliferation, ultimately sealing off the aneurysm~\cite{Briganti2015a}.
Contrary to microsurgical clipping and coil embolization,
the aneurysm closure occurs gradually over the course of 6 to 12 months after treatment~\cite{Pierot2017a,Szikora2015a}.

Flow diverters were initially used in aneurysms that are not suitable for coiling due to their size and shape,
i.e. large/giant wide-necked aneurysms.
Recently, their spectrum has extended also to fusiform aneurysms, dissecting aneurysms and recanalized aneurysms,
where coiling can only be performed with a high risk of failure.
Furthermore, flow diversion can also be used in very small aneurysms untreatable by standard coiling.
Flow diverters cannot be used for aneurysms that are located at vessel bifurcations
due to the risks of side branch occlusion~\cite{Pierot2011a}.
Currently, complications are higher than after coiling.
Possible risks include ischemic stroke, perforator infarctions,
vessel- or aneurysm ruptures or secondary flow diverter occlusion~\cite{Pierot2017a}.

\subsection{Flow disruption devices}

Flow disruption is the newest endovascular approach for the occlusion of aneurysms.
Here, an intrasaccular device is strategically positioned to modify blood flow dynamics at the aneurysm's neck.
This modification induces the formation of an aneurysmal thrombus,
functioning in a manner conceptually similar to intravascular flow diversion.
The most commonly used device for flow disruption is the {\web} device
which was initially developed in 2010~\cite{Pierot2016a}.
The device is predominantly used for medium to large wide-neck bifurcation aneurysms,
but recent technological advancements in the device have expanded its potential applications to more distal bifurcation aneurysms,
sidewall aneurysms and smaller aneurysms up to 2mm (see Figure~\ref{fig:DeviceAngioWEB}).
Because the device is relatively new, there are only limited long-term results available~\cite{Pierot2017a}.

\section{From imaging data to computational geometries} \label{sec:FromImagingToGeometries}


Since this work sets out to aid medical doctors in evaluating treatment options and predicting treatment outcome for individual patients,
the extraction of patient-specific vessel geometries with physiological relevant boundary conditions builds the foundation for computational modeling.
Hence, we outline our preprocessing pipeline to arrive at patient-specific geometries suitable to be used in numerical simulations.

\subsection{Case selection, medical imaging and segmentation} \label{subsec:MedicalSelectionAndSegmentation}


Patient-specific geometries were taken from patients with unruptured intradural aneurysms of the anterior or posterior circulation treated at the Klinikum rechts der Isar, Technical University of Munich, between 2013 and 2018 and continuous clinical and radiographic follow-up for $\geq$5 years.
Written informed consent was waived by the local institutional review board for this study because of its retrospective character and the analyses being based only on data acquired during the clinical routine.
Of all 314 patients, 267 aneurysms were coiled, 29 aneurysms were treated by implantation of a flow diverter,
and 18 with {\web} device implantations.

The procedure to obtain accurate vessel geometries is given by the following workflow:
First, one performs an image acquisition of the ipsilateral vascular tree,
beginning at the skull base with 3D rotational angiography.
For image acquisition, Klinikum rechts der Isar, Technical University of Munich,
uses a Philips Azurion 7 Neuro Suite (Philips Medical Systems, Netherlands).
The resulting image exhibits a size of 384 x 384 x 384 pixels with a resolution of resolution 3.7 pixels / mm.
In a second step,
one performs semi-automatic segmentation of the intracranial vasculature with ITK SNAP~\cite{py06nimg},
a free, open-source, multi-platform software application used to segment structures in 3D and 4D biomedical images.
Starting from a manual selection of a few seeds points to define the interior of the vessels,
a contour evolution algorithm expands the vessel contour to align it with the grey-scale contrast being present at the vessel-wall locations automatically
ITK SNAP applies a contour-evolution algorithm to expand the vessel contour to automatically align with the grey-scale contrast at the vessel wall locations.
Finally, the mask of the segmented vessels is exported in STL format,
yielding a comprehensive representation of the vessel geometries for further analysis.

\subsection{3D Geometry preprocessing \& meshing} \label{subsec:GeometryPreprocessing}

\subsubsection{Geometry extraction, adaption and alignment}
\label{sec:GeometryExtraction}
\begin{figure}
\centering
~\hfill~
\begin{minipage}{0.4\textwidth}
\subfigure[Original vascular tree from dataset, here \#C0074b from~\cite{aneurisk}]{\includegraphics[height=8cm]{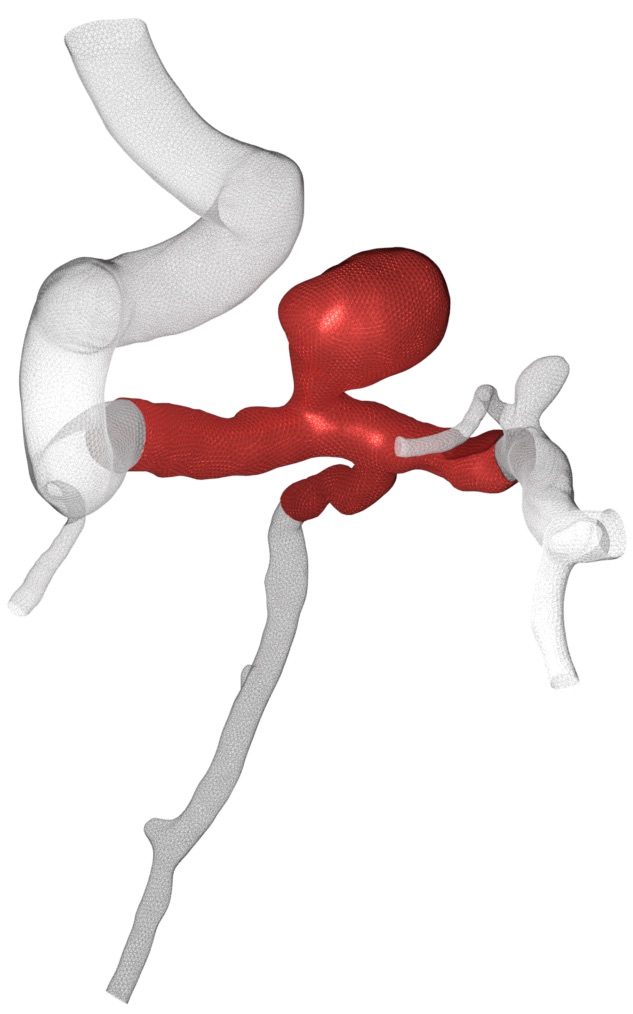}}
\end{minipage}\begin{minipage}{0.5\textwidth}
\subfigure[Domain of interest (red) with a flow extension (light blue)]{\includegraphics[height=3cm]{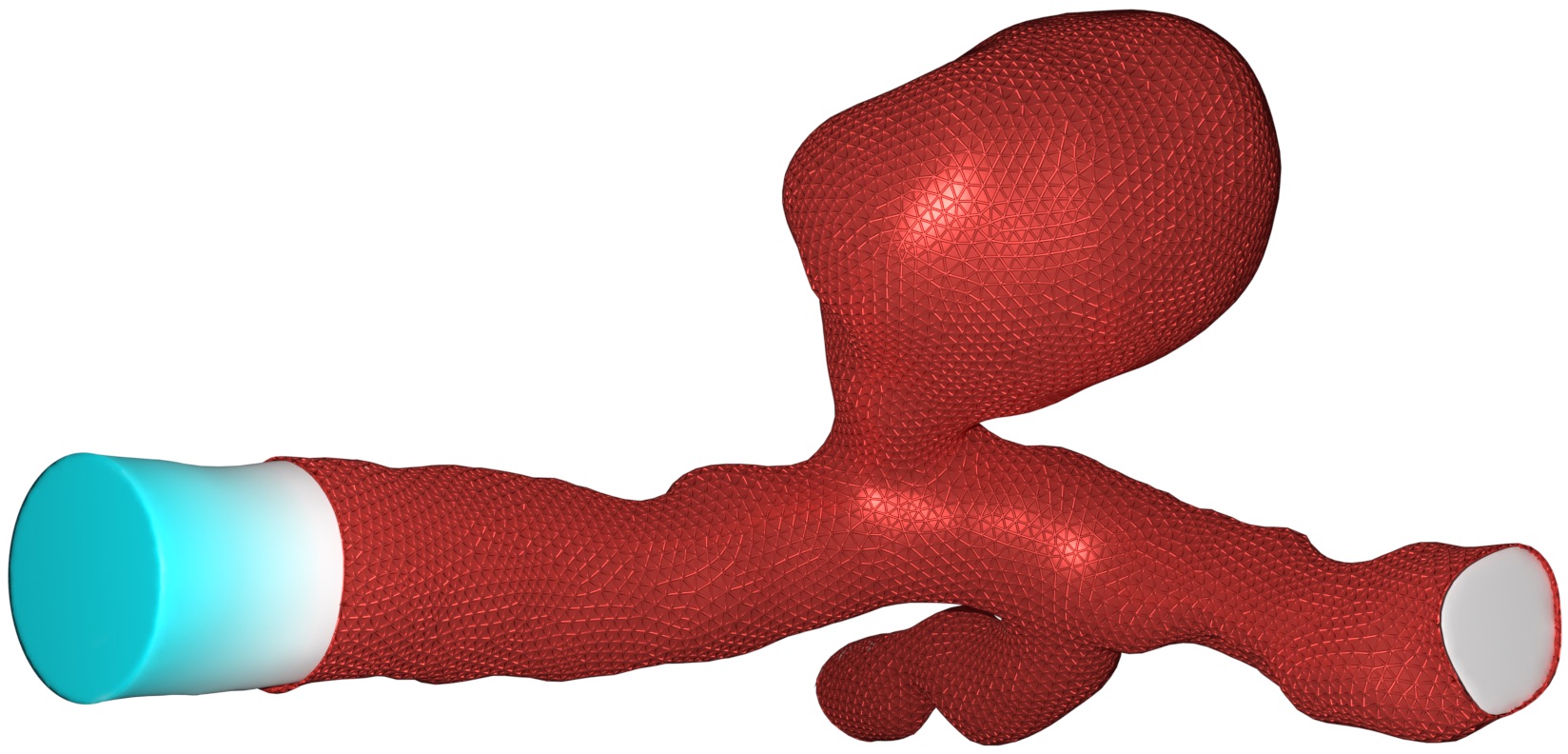}}
\subfigure[Inflow (orange) and outflow (green) surfaces (second outflow surface in the back not visisble in this perspective). A digital version of this data set is available in~\cite{Frank2024a}.]{\label{fig:GeometryPreProcessing}\includegraphics[height=3cm]{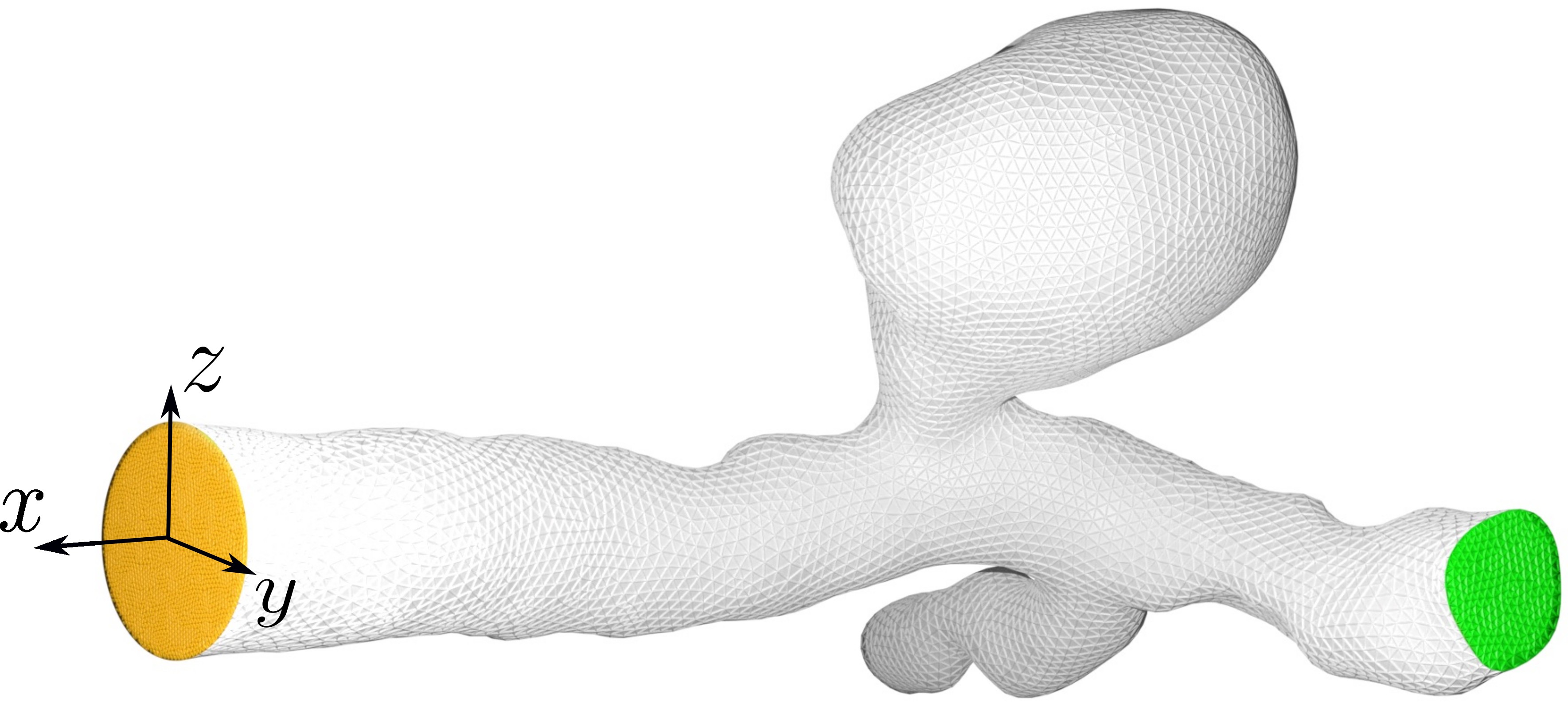}}
\end{minipage}
\caption{Pre-processing of vessel geometries: After extraction of the domain of interest, the inflow surface is extruded to create a flow extension with a circular inflow cross section. For the perception of the used colors we refer the reader to the online version of the article.}
\label{fig:CutNExtendGeometry}
\end{figure}

Besides the acquisition of patient-specific geometries as outlined in Section~\ref{subsec:MedicalSelectionAndSegmentation},
we also rely on publicly available online databases such as~\cite{aneurisk, PozoSoler2017, IntrA} for a rich set of different test-geometries, mostly available in \texttt{.vtk}, \texttt{.stl} or \texttt{.obj} file format. The vascular trees resolved in these datasets or by segmentation of angiographic imaging data are usually much larger than the vessel portion to be simulated in full three-dimensional (3D) resolution. Consequently, we first truncate the vessel geometry at roughly two aneurysm diameters up- and -- depending on the available dataset -- downstream the position of the aneurysm by a simple planar cut, see Figure~\ref{fig:CutNExtendGeometry}, with its exact position being up to human decision. For this, but also upcoming mesh manipulation tasks as well as final visualization results, we have found the software tools \textit{blender}\footnote{blender software webpage: \url{https://www.blender.org/}} and \textit{Vascular Modelling Toolkit (vmtk)}\footnote{vmtk software webpage: \url{http://www.vmtk.org/}}~\cite{vmtk} extremely useful. If the segmentation mesh data is slightly noisy, contains triangles with insufficient aspect ratio, or exhibits other forms of poor mesh quality, the Laplace-Taubin algorithm~\cite{taubin1995signal} for smoothing the discrete arterial surface has proven to be effective. For a simple and  optimal integration of the algorithm into the workflow, the implementations in \textit{MeshLab}\footnote{MeshLab software webpage: \url{https://www.meshlab.net/}} \cite{LocalChapterEvents:ItalChap:ItalianChapConf2008:129-136} or \textit{vmtk}~\cite{vmtk} are used where for example the choice of the scale parameters $\lambda=0.5$ and $\mu=-0.53$, see~\cite{taubin1995signal} and the standard settings in \textit{MeshLab}, for 10 refinement steps provides a good starting point for a smooth arterial surface. 

To simplify the imposition of physiologically meaningful inflow boundary conditions, we sometimes add flow extensions to the inflow surface, such that the flow field is fully developed when it reaches the beginning of the truncated vessel domain~\cite{moyle2006inlet, tang2020fluid}. These straight geometric extension are extruded in orthogonal direction of the original inflow surface and gradually interpolate between the original, not necessarily circular cut-off opening of the vessel and the artificially circular ending of the extension (light blue part in Figure~\ref{fig:CutNExtendGeometry}), which will then constitute the inlet boundary of the simulation model. Alternatively, we can prescribe arbitrary inflow profiles directly on the original cut-off boundary of the vessel. An extension of outflow areas is sometimes beneficial for flow simulations.


We then rotate the entire geometry, such that the outward normal $\Vec{n}'$ of the inlet cross section is aligned with the negative $x$-axis of the coordinate system and the geometric center~$\Vec{c}'$ of the inflow cross section resides at~$\Vec{0}$, see Figure~\ref{fig:GeometryPreProcessing}. Consequently, a prescribed inflow profile in positive $x$-direction implies an orthogonal inflow. The geometry's translation to map the inlet's geometric center \textit{before} translation~$\Vec{c}$ onto $\Vec{c}'=\Vec{0}$ is simply the addition of~$-\Vec{c}$ to each vertex coordinate vector~$\Vec{v}$ of the vessels-surface triangulation. To define an unique and objective rotation, we use  Euler's finite rotation formula (Rodrigues' formula). For convenience of the reader we state it as Remark~\ref{rem:RodriguesFormula}.
\begin{Remark}
\label{rem:RodriguesFormula}
    \textit{Let $\vec{a},\vec{b}\in\mathds{R}^3$ be two unit-vectors enclosing the angle $\vartheta=\angle(\Vec{a},\Vec{b})\in(0,\pi)$ between them. Defining 
    \begin{itemize}
        \item $\vec{w}:=\Vec{a}\times\Vec{b}$,
        \item $c:=\cos(\vartheta)=\langle \vec{a},\vec{b}\rangle$,
        \item  $s:=\sin(\vartheta)=\|\Vec{w}\|_2$ and
        \item $\hat{\Vec{w}}:=\Vec{w}/s$
    \end{itemize}
    the matrix representation $\tens{R}$ of the rotation by the angle $\vartheta$ around the rotational axis $\hat{\Vec{w}}$ that maps $\vec{a}$ onto $\Vec{b}$ is given by:
    \begin{equation}\label{eq:Rotation}
    \tens{R} = \tens{I}+s\,[\,\hat{\vec{w}}\,]_{\times}+(1-c)\,[\,\hat{\vec{w}}\,]_{\times}^2, \quad \textup{with} \quad [\,\hat{\vec{w}}\,]_{\times}:=\left(\begin{array}{ccc}
        0 & -\hat{w}_3 & \hat{w}_2 \\
        \hat{w}_3 & 0 & -\hat{w}_1 \\
        -\hat{w}_2 & \hat{w}_1 & 0
    \end{array}\right)
    \end{equation}
    For the special cases~$\vartheta=0$ (iff $c=1$) and~$\vartheta=\pi$ (iff $c=-1$),
    one instead uses~$\tens{R}=\tens{I}$ and~$\tens{R}=-\tens{I}$, respectively.}
\end{Remark}

This formula (in different notation) was already known to Euler~\cite{cheng1989historical} and is nowadays also often used in the differential-geometrical context of frames for space-curves~\cite[eq. 2.4]{jawed2018primer}, see Section~\ref{subsec:Stents} for such an application even in this work. In the present case however, we can simply use it for the desired geometry rotation by choosing $\Vec{a}=\vec{n}$ to be mapped onto $\Vec{b}=(-1,0,0)^{\top}$ and again apply the rotation matrix~$\tens{R}$ to each vertex coordinate vector~$\Vec{v}$ of the vessel's surface triangulation.

\subsubsection{Geometric measures of vessel and aneurysm}\label{subsec:Measures}

After obtaining a geometric representation of the vessel and its pre-processing, another important step in the geometry acquisition is the \textit{measurement} of the resulting geometry. This step can encompass for example the computation of the vessel's centerline together with its corresponding arc length, the extraction of the local vessel radius along the centerline, or assessing the diameter or volume of the aneurysm itself. Such quantities are required for the correct choice of design parameters and placement of devices in Section~\ref{sec:DeviceModelling}.
In this regard, again the \textit{vmtk}~\cite{vmtk} package has proven helpful not only for measuring cerebral aneurysms, but also for aortic aneurysms~\cite{vilalta2016correlation}.

Measurement based on the centerline will become useful in geometric stent modeling (see Section~\ref{subsec:Stents}). Geometric and volumetric information about the aneurysm itself will be required for treatment options such as coiling or {\web}-device insertion. For coiling, it is required to compute the volume packing density, i.e. the ratio of the volume occupied by the coiling wire to the entire volume of the aneurysm sac (see Section~\ref{subsec:Coils}), and hence occlusion ratio which can be used as an indicator for the success of the coiling procedure.

To extract volume information from the geometrical model, we virtually separate the sac of the aneurysm from the adjacent vessel and close it at the ostium with a triangulated surface patch, hence arriving at the capsulated aneurysm shapes depicted in Figure~\ref{fig:CoilModel}. These then serve as ``cages'' for the coil-placement simulation. If we denote by $\mathcal{T}$ the set of all triangles constituting the surface mesh of the (encapsulated) aneurysm and by $\mathcal{C}$ only the triangles used to close the aneurysm at the ostium, then the surface and volume of the aneurysm-sack can then be computed as 
\begin{subequations}
\begin{align*}
    S&=\sum_{T\in\mathcal{T}\backslash \mathcal{C}}A_T,\qquad\textup{with}\qquad A_T=\frac{1}{2}\,\|(\Vec{v}_{T,2}-\Vec{v}_{T,1})\times(\Vec{v}_{T,3}-\Vec{v}_{T,1})\|_2\\
    V&=\frac{1}{3}\sum_{T\in\mathcal{T}}\langle \Vec{v}_{T,1},\vec{n}_T\rangle~A_T 
\end{align*}
\end{subequations}
where $\Vec{n}_T$ is the outward unit normal vector of triangle $T\in\mathcal{T}$ and $\vec{v}_{T,i},i=1,2,3$ the spatial locations of the three vertices of triangle $T\in\mathcal{T}$,
also see~\cite[Section IV.1]{arvo2013graphics}.

\subsubsection{Meshing}
\label{sec:Meshing}

The acquisition of vessel geometries as outlined in Section~\ref{sec:GeometryExtraction} delivers a triangular mesh of the vessel's surface.
Such a pure surface representation is sufficient for some types of analysis, for example for the simulations of the coil deployment process (see Section~\ref{subsec:Coils}). There, we will neglect the presence of the blood for simplicity, but will account for the contact interaction between the coiling wire and the aneurysm wall.

When accounting for blood and its flow through the vessel and aneurysm, a pure surface representation is insufficient and a volumetric representation of the lumen (volume covered by blood flow) along with a suitable computational mesh is required. The actual choice of the computational mesh highly depends on the applied numerical methods. In this paper, we will employ the {\LBM} and the {\FEM}.

On the one hand, {\FEM} can deal with a variety of mesh topologies and cell types. When aiming at tetrahedral meshes for example, \textit{vmtk} offers a capable mesh generator to setup the overall mesh structure of the vessel (see Figure~\ref{fig:Grids}), which is then fine-tuned and completed e.g. by assigning flags and extracting node-sets at the different boundary parts, in \textit{cubit}\footnote{cubit software webpage: \url{https://coreform.com/}}\cite{cubit}. Depending on the previous used segmented data, the boundary parts can be identified and assigned within \textit{cubit}. Additionally, cubit offers the capabilities to create and refine meshes from parametrized virtual geometry aneurysms, as they will be used in Section~\ref{subsec:poroexample}.

On the other hand, {\LBMs} require a structured, \textit{Cartesian} lattice of cubes, that then represent any curved boundaries via a staircase-approximation. Although \textit{cubit} also offers such opens e.g. via its \textit{sculpt} functionality~\cite{owen2017hexahedral}, we rely on the built in meshing-tools of the {\LBM} software of our choice, \textit{waLBerla}\footnote{waLBerla software webpage: \url{https://walberla.net/}} \cite{bauer2021lbmpy, bauer2019code, bauer2021walberla}. This high-performance {\LBM} solver generates and also partitions such a lattice into subdomains ready for parallel computing. A signed distance octree algorithm delivers macro blocks assigned to each process, while their union covers the entire geometry, and distinguishes into two separate structures of {\LBM}-cells, namely fluid cells inside the geometry's boundary and solid cells outside, respectively. In this work, we focus on the fluid cells for {\LBM} and just use the solid cells to define the boundary.
See Figure~\ref{fig:Grids} for a depiction of an {\LBM}-suited mesh also roughly depicting its construction process by ``filling out'' the surface triangular mesh ``from the inside''.

\subsubsection{Calculation of surface normals}

In {\FEM}, the boundary fitted meshes facilitate a straightforward computation of a field of outward unit normal vectors~$\Vec{n}$, where the smoothness of this field only depends on the mesh resolution. In contrast, the computation of outward unit normal vectors for arbitrarily oriented and shaped domain boundaries is far from trivial for the voxel-type grids used within {\LBMs}. In fact, due to the staircase approximation, {\LBMs} a priori lack a smooth outward unit normal vector field $\Vec{n}$ at boundary cells, which is necessary for the evaluation of key quantities in the assessment of the flow state in biomedical flows, for example the wall shear stress computed as $\vec{\textbf{WSS}}=\tens{\sigma}\vec{n}-(\Vec{n}^{\top}\tens{\sigma}\Vec{n})\Vec{n}$, respectively its norm $\textup{WSS}=\|\textbf{WSS}\|_2$ close to the vessel and aneurysm wall (see Section~\ref{sec:NumericalSimulations}).

\begin{figure}
\centering
\subfigure[Cartesian lattice used for an {\LBM}-simulation]{\label{fig:LBMMesh}\includegraphics[height=4.5cm]{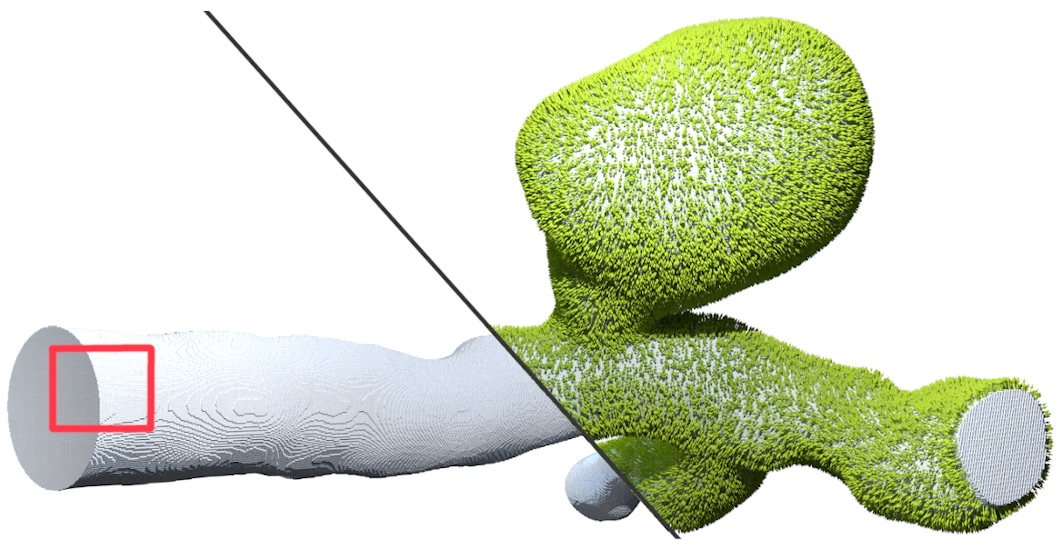}}
\subfigure[Close-up view into red rectangle with three sectors showing steps~I, II, and III for the computation of the outward unit normal vector field]{\label{fig:LBMNormalField}\includegraphics[height=4.5cm]{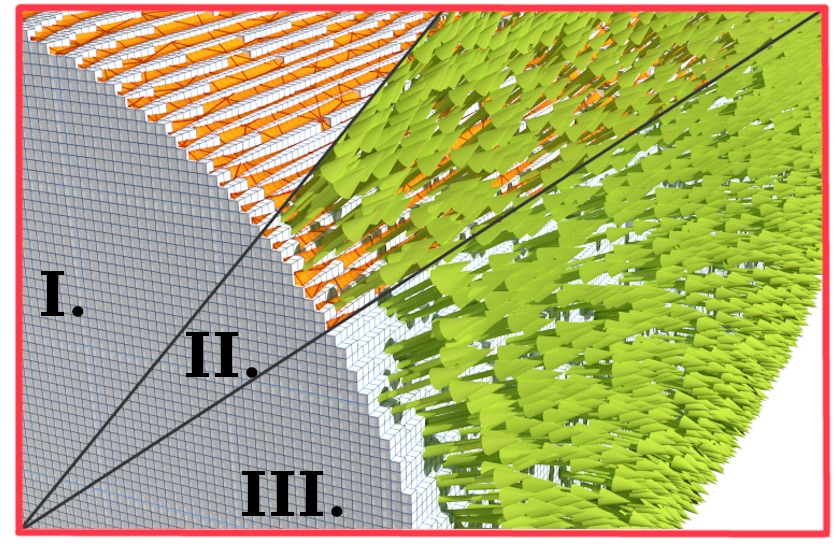}}
\caption{\emph{Left:} Cartesian grid used for an {\LBM}-simulation
(left sector: voxel-based geometry, right sector: outward unit normal field depicted by green arrows).
\emph{Right:} Zoom into the red marked region of \ref{fig:LBMMesh} to showcase the staircase approximation of curved boundaries.
Sector I shows the grid constructed \textit{inside}  the (orange) triangulation of the vessel surface.
Sector II. illustrates normal vectors originating from lattice-\ cells and surface mesh intersections,
while Sector III. depicts the (continuous to mesh resolution) field of normal vectors without the surface mesh.
\emph{Note:} For an improved visualization, only a random selection of cell normal vectors is shown.}
\label{fig:Grids}
\end{figure}

In order to reconstruct such a smooth outward unit normal field, purely grid-based approaches often result in only a finite number of different normal vectors, which often exhibit sharp changes in close proximity. As counter measure, we instead aim to employ the underlying surface mesh of the original vessel geometry,
that served as starting point for the distance octree algorithm to create the {\LBM}-grid (see Section~\ref{sec:Meshing}). We follow the approach from~\cite{normalCalc}, where an outward-facing normal field $\Vec{n}(\Vec{x})$ can be obtained from intersecting boundary cells of the {\LBM} grid with the original surface mesh's triangles and a subsequent averaging over these triangles individual normal values. Similar to the {\FEM}, the resulting field of normal vectors is again smooth up to the mesh resolution of the original vessel geometry. In Figure~\ref{fig:Grids}, the intersection situation is depicted together with the resulting outward unit normal vector field~$\Vec{n}(\Vec{x})$.

\section{Continuum-mechanics \& numerics of flow in aneurysms} \label{sec:ContinuumMechanicsAndNumericsOfFlowInAneurysms}

We now introduce the analytical models and numerical methods for the simulation of hemodynamics within the vessels.
In this contribution, we work with two different continuum models and approaches:
Section~\ref{subsec:FreeFlow} introduces the Navier-Stokes equations for incompressible fluid flow and their solution via Lattice Boltzmann methods (LBMs) on non-deforming domains,
that will be used to describe the blood flow if most geometrical features of medical devices are resolved.
Section~\ref{subsec:PoroElasticFlow} details the governing equations and finite element discretization of a poro-elastic model
to represent treated aneurysms on a deforming domain in a homogenized manner, i.e.\ without resolving geometrical details of medical devices.
Both formulations differ in the increased computational performance of the LBM
and in the ability of FEM to analyze deformable domains.


\subsection{Free-flow fluid model w/o fully resolved endovascular devices} \label{subsec:FreeFlow}

\subsubsection{Governing equations} \label{subsubsec:FreeFlowAnalysical}
Velocity~$\vel$ and pressure~$\p$ in the flow field are governed by the Navier-Stokes equations for incompressible fluid flow, reading
\begin{subequations}
\label{eq:NavierStokesSystem}
\begin{align}
    \frac{\partial\vel}{\partial t}+\left(\vel\cdot\grad\right)\vel+\frac{1}{\den}\left(\grad \p-\grad\cdot\left(2\vis\tens{\epsilon}(\vel)\right)\right)&=\Vec{f} \label{eq:NavierStokes1}\\
    \grad\cdot \vel&=0 \label{eq:NavierStokes3},
\end{align}
\end{subequations}
where~$\rho$, $\tens{\epsilon}(\vel)=\frac{1}{2}\left(\grad\vel+\grad\vel^{\top}\right)$, and~$\Vec{f}$ denote the density, strain rate tensor, and body force vector, respectively.
We note that~\eqref{eq:NavierStokes1} is formulated such that it is extensible to non-Newtonian fluid models (see Section~\ref{sec:NonNewtonian}).

To complete the system of equations in \eqref{eq:NavierStokesSystem},
we impose boundary conditions to the individual physical parts of the computational domain's boundary,
in particular at the inflow, the outflow, and the interface to the vessel wall.

\begin{figure}
\centering
\includegraphics[height=7cm]{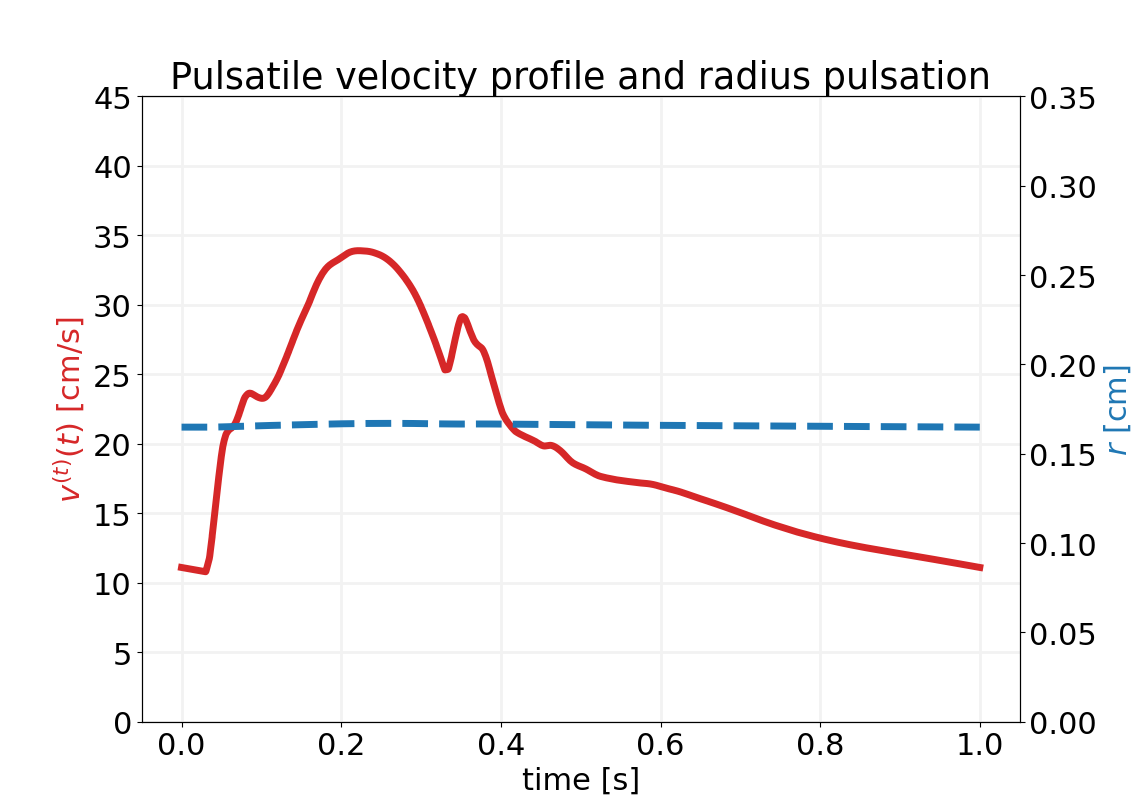}
\caption{Time-dependent pulsatile velocity amplitude profiles $v^{(t)}(t)$ and vessel radius pulsation $r(t)$
from the 0D-1D model~\cite{fritz20221d} over one heart cycle in the basilar artery (vessel \#\,22 in Figure~1 in~\cite{alastruey2007modelling}).
}
\label{fig:PulsationProfile}
\end{figure}

At the inlet cross section,
we impose a Dirichlet condition on the velocity field, which might depend on space $\Vec{x}$ and time $t$, hence reading~$\vel_{\mathrm{in}} = \vel_0(t,\vec{x})$.
In the special case of a \textit{circular} inflow boundary aligned with the $x$-coordinate axis as obtained from Section~\ref{subsec:GeometryPreprocessing},
we can directly prescribe a radially symmetric generalized Poiseuille velocity profile $\vel_0(t,r)$,
where $r$ denotes the radial distance from the inlet circle's center point.
Denoting the radius of the inlet-circle by $R$, the profile takes the form
\begin{equation}\label{eq:GammaPoiseuille}
\vel_0(t,r)=v^{(t)}(t)\cdot \left(\begin{array}{c}
v_{\perp}(r) \\
0 \\
0
\end{array}\right),\qquad \textup{with} \qquad v_{\perp}(r)=\frac{\gamma+2}{\gamma}\left(1-\left(\frac{r}{R}\right)^{\gamma}\right)
\end{equation}
where $\gamma$ is a shape parameter yielding the standard parabolic Poiseuille profile for the choice $\gamma=2$,
while larger values such as $\gamma=9$ are reported to fit well to experimental data~\cite{alastruey2008lumped, smith2002anatomically}
and to resemble a non-Newtonian, more plug-like, flow profile~\cite{vcanic2003mathematical} (also see Section~\ref{sec:NonNewtonian}).
Finally, $v^{(t)}(t)$ describes the time-dependent amplitude of the inflow velocity over the course of a heart beat.
We obtain such a temporal pulsatile profile from a desired position within the cardiac cycle from the 0D-1D model described in~\cite{fritz20221d}
with differences ranging from large oscillations in comparably also larger vessels up to a nearly temporally constant flow in the small capillaries.
In Figure~\ref{fig:PulsationProfile},
we plot such a velocity amplitude profile for the basilar artery (vessel \#\,22 in Figure~1 in~\cite{alastruey2007modelling}),
where the pulsatile nature in the velocity amplitude is still clearly visible,
while the vessel radius pulsation is essentially negligible.

At any outlet cross-section,
we impose a constant zero pressure boundary condition, meaning that the flow can freely leave the domain.

For this pure flow model, all vessel and aneurysm walls are assumed to be rigid and to carry no-slip conditions for the velocity, i.e. $\vel=\vec{0}$ at all walls.
While this rigidity assumptions poses a simplification,
it reduces the complexity of the model considerably and enables quick first insights into the flow in cerebral vessels and aneurysms.
We will lift the assumption of rigid vessel walls in Section~\ref{subsec:PoroElasticFlow}.
%
Similarly, medical treatment devices are represented by fully resolved, but static solid obstacles (see Sections~\ref{sec:DeviceModelling} and~\ref{sec:NumExCoiling} for details).
At their boundaries, we impose no-slip conditions as well.
Admittedly, the assumption of a static obstacle neglects effects such as coil compaction~\cite{kallmes1999histologic, reul1997long, sluzewski2004relation},
which deserve a detailed analysis in a separate work due to their medical relevance for the recanalization and eventually regrowth risk of aneurysms.

\subsubsection{Extension to non-Newtonian rheology modeling}
\label{sec:NonNewtonian}

In order to extend the constitutive behavior of the fluid model to account for the non-Newtonian shear-thinning nature of (human) blood,
we revisit the viscosity $\vis$ in~\eqref{eq:NavierStokes1} and formulate it as a dependent quantity $\vis=\vis(\dot{\gamma})$ for decreasing shear rates $\dot{\gamma}$.
Herein, $\dot{\gamma}$ is defined as $\dot{\gamma}:=2\sqrt{\tens{\epsilon}(\vel):\tens{\epsilon}(\vel)}$
with $\tens{\epsilon}(\vel):=\frac{1}{2}\left(\grad\vel+\grad\vel^{\top}\right)$ again being the symmetric strain rate tensor.
To recover the case of a viscous Newtonian fluid,
the stress-strain rate relation adopts the form $\tens{\sigma}=2\vis\tens{\epsilon}(\vel)$, taking the Frobenius-norm $\varsigma:=\sqrt{\tens{\sigma}:\tens{\sigma}}=\vis\dot{\gamma}$ where $\vis\in\R_+$ is just a constant. Starting from here, many different models for shear rate dependent viscosity functions $\vis(\dot{\gamma})$ do exists, yielding the nonlinear expression $\varsigma=\vis(\dot{\gamma})\dot{\gamma}$.

A common representative are power laws~\cite{PowerLawLBM} of the form~$\vis(\dot{\gamma})=\eta\dot{\gamma}^{n-1}, \eta\in\R_+, n\leq 1$,
where~$\eta$ again takes the role of the constant base viscosity in the Newtonian case $n=1$.
It should be noted that~$\vis_{\infty}:=\lim_{\dot{\gamma}\to\infty}\vis(\dot{\gamma})=0$ and~$\vis_0:=\lim_{\dot{\gamma}\to 0}\vis(\dot{\gamma})$ diverges towards~$\infty$ for such power law models.
Another model in this context is the Carreau-Yasuda model, which allows to smoothly attain finite values in these limit situations.
Here, the viscosity function is given by~$\vis(\dot{\gamma})=\vis_{\infty}+(\vis_0-\vis_{\infty})\left[1+(\lambda\dot{\gamma})^a\right]^{\frac{n-1}{a}}$~\cite{boyd2007analysis}.
A list of further classical non-Newtonian constitutive equations can be found in~\cite{cho1991effects}, from which we choose the Casson model for the present work.

The Casson model assumes the existence of a yield shear stress value~$\varsigma_Y$ of blood,
which is explained by interaction forces between red blood cells in~\cite{cho1991effects}.
In flow regimes with shear stress values below~$\varsigma_Y$,
these interaction forces are not overcome and the respective fluid parts behave solid-like with an apparent viscosity of~$\infty$.
Only once the threshold~$\varsigma_Y$ is exceeded,
deformation in form of strain occurs.
In the Casson model, the specific version of the stress-strain relation known as Casson's equation is given by:
\begin{align*}
	\sqrt{\eta\dot{\gamma}}=\begin{cases}
		\sqrt{\varsigma}-\sqrt{\varsigma_Y},&\varsigma>\varsigma_Y\\
		0,&\textup{else.}
	\end{cases}
\end{align*}
or, for a direct expression of the apparent viscosity $\vis$ or the strain rate $\dot{\gamma}$:
\begin{equation}\label{eq:Casson2}
	\vis(\dot{\gamma})=\begin{cases}
	\eta\cdot\frac{\varsigma}{(\sqrt{\varsigma_Y}+\sqrt{\eta\dot{\gamma}})^2}&\varsigma>\varsigma_Y\\
	\infty,&\textup{else.}
	\end{cases},\qquad\qquad
    \dot{\gamma}=\begin{cases}
	\frac{\varsigma_Y}{(\sqrt{\vis(\dot{\gamma})}-\sqrt{\eta})^2}&\varsigma>\varsigma_Y\\
		0,&\textup{else,}
	\end{cases}
\end{equation}
A detailed account including the derivation of an analytical solution to a Casson model equivalent of a Hagen-Poiseuille flow can be found in~\cite{fung2013biomechanics},
being the aforementioned profile that \eqref{eq:GammaPoiseuille} resembles for values of e.g. $\gamma=9$.

\subsubsection{Numerical treatment via the Lattice Boltzmann method}\label{subsubsec:FreeFlowNumericalTreatmentLBM}
Within this sub-section we will give a very brief account on the {\LBM} to be applied to the presented model-problem. In addition to the references given in this sub-section, we refer the generally interested reader to the textbook~\cite{kruger2017lattice} for an in-depth introduction to a variety of aspects of the {\LBM}.

\paragraph{Gasdynamical background}
The first of our numerical methods to solve the hemodynamical fluid model equations \eqref{eq:NavierStokesSystem} is the {\LBM}.
In contrast to classical numerical schemes such as {\FEM} or finite differences being based directly on the macroscopic quantities like velocity $\vel$ or density $\rho$ as well as the continuum-level equations like Navier-Stokes system (in weak or strong form),
the {\LBM} operates on a mesoscopic level.
While it is not yet dealing with individual particles as in molecular dynamics,
the central term of {\LBM} are so called particle distribution functions $f(t,\vec{x},\Vec{\xi})$ depending on time $t$, spatial coordinates $\Vec{x}$ and \textit{particle} velocity $\Vec{\xi}$ representing the probability to encounter a particle at time $t$ at spatial location $\Vec{x}$ with velocity $\Vec{\xi}$.
Boltzmann's equation
\begin{equation}\label{eq:Boltzmann}
    \frac{\textup{d}f}{\textup{d}t}=\frac{\partial f}{\partial t} + \vec{\xi}\cdot\nabla_{\vec{x}}f + \frac{\Vec{G}}{\rho}\cdot\nabla_{\vec{\xi}}f=\Omega(f)
\end{equation}
now describes the total change $\textup{d}f/\textup{d}t$ within these populations due to external forces $\Vec{G}$ as well as particle collisions that are modeled by a highly complex collision operator $\Omega(f)$,
see e.g.~\cite{hanel2006molekulare} for a derivation.
Due to its complexity in~\cite{bhatnagar1954model}, a simplification $J(f)$ for $\Omega(f)$ was proposed reading
\begin{equation*}
    \Omega(f)\approx J(f)=-\frac{1}{\tau}\left(f-M\right),\qquad\textup{with}\qquad M=M(\Vec{\xi};\rho,\vel,T)=\rho\left(\frac{1}{2\pi R_sT}\right)^{3/2}\exp\left(-\frac{\|\Vec{\xi}-\vel\|_2^2}{2R_sT}\right)
\end{equation*}
being the so called gas-dynamical \textit{Maxwell-equilibrium} depending on absolute temperature $T$, the specific gas constant $R_s$ of the fluid
as well as the macroscopic quantities $\rho$ and $\vel$ which are related to the particle distribution functions via:
\begin{equation}\label{eq:MacroQuant}
    \rho(t,\Vec{x})=\int_{\mathds{R}^3}f(t,\Vec{x},\Vec{\xi})~\textup{d}\Vec{\xi},\qquad \vel(t,\Vec{x})=\frac{1}{\rho(t,\Vec{x})}\int_{\mathds{R}^3}\vec{\xi}\,f(t,\Vec{x},\Vec{\xi})~\textup{d}\Vec{\xi}.
\end{equation}
Hence, the simplified collision operator $J(f)$ corresponds to a linearization around the equilibrium distribution $M$ with $\tau$ being a relaxation time later to be related to the fluid's viscosity $\vis$.

\paragraph{Discretization in LBM}
Discretization within the {\LBM} is applied in ($d$-dimensional) physical space
by means of the aforementioned Cartesian cube-lattice on the computational domain
and in velocity space by means of a finite, discrete set of lattice velocities $\Vec{c}_i,i=0,1,\dots,n-1$ per grid-cell,
which then forms the so-called D$d$Q$n$-stencil.
Most common representatives are the D$2$Q$9$ and D$3$Q$27$ stencils in 2D and 3D, respectively.
Neglecting force contributions (i.e. $\Vec{G}=\vec{0}$) for the moment,
applying a forward Euler discretization to the temporal derivative in \eqref{eq:Boltzmann}
and a finite-difference approximation in direction $\Vec{\xi}$ for the spatial gradient therein
as well as restricting the velocity variable $\Vec{\xi}$ to the discrete set of lattice-velocities,
equation \eqref{eq:Boltzmann} -- with the simplified collision operator and evaluated at time $t$ in cell $\Vec{x}$ -- becomes:
\begin{equation}\label{eq:LBMUpdate}
    f_i(t+\Delta t,\vec{x}+\vec{c}_i\Delta t) = f_i(t,\vec{x}) - \frac{\Delta t}{\tau}\left(f_i(t,\Vec{x})-f_i^{\textup{(eq)}}(t,\Vec{x})\right)
\end{equation}
Herein, $f_i$ is the particle distribution associated to the $i$-th lattice velocity and $f_i^{\textup{(eq)}}$ the following, truncated after quadratic order, approximation to the Maxwell-equilibrium from above:
\begin{equation}\label{eq:feq}
    f_i^{\textup{(eq)}} = w_i\rho\left(1+\frac{\vec{c}_i\cdot\vel}{c_s^2}+\frac{(\Vec{c}_i\cdot\vel)^2}{2c_s^4}-\frac{\vel\cdot\vel}{2c_s^2}\right)
\end{equation}
where $c_s$ is the - in lattice units - speed of sound (in case of D2Q9 and D3Q27 it is $1/\sqrt{3}$, $w_i\in\R, i=0,1,\dots n-1$ are method-specific weights associated to the lattice-velocities and furthermore the continuous expressions for density $\rho$ and macroscopic velocity $\vel$ from \eqref{eq:MacroQuant} are replaced by the discrete versions:
\begin{equation}\label{eq:DiscrMacroQuant}
    \rho=\sum_{i=0}^{n-1}f_i,\qquad \vel=\frac{1}{\rho}\sum_{i=0}^{n-1}\vec{c}_if_i
\end{equation}
The three equations \eqref{eq:DiscrMacroQuant}, \eqref{eq:feq} and \eqref{eq:LBMUpdate}, inserted into each other in this order, then form one single time step and hence the very core of an {\LBM}-scheme,
consisting of the \textit{collision} step being the purely cell-local right hand side of~\eqref{eq:LBMUpdate}
and the \textit{streaming} step, which involves only the direct neighborhood of a given cell and stands for the propagation of population data,
i.e. the left hand side of~\eqref{eq:LBMUpdate}.
For the high-performant execution of these steps as well as their handling in a parallel distributed memory and multi-processor architecture,
we employ the software framework \textit{waLBerla}\footnote{waLBerla software webpage: \url{https://walberla.net/}} \cite{bauer2021lbmpy, bauer2019code, bauer2021walberla} developed at the FAU Erlangen.
Finally, the link to the macroscopic Navier-Stokes system \eqref{eq:NavierStokesSystem}
is established via a so called \textit{Chapman-Enskog analysis}~\cite{chapman1990mathematical, chopard2002cellular, silva2014truncation}
which also relates the relaxation time parameter~$\tau$ within~\eqref{eq:LBMUpdate} to the macroscopic viscosity of the fluid~$\vis$ within~\eqref{eq:NavierStokes1} via
\begin{equation}\label{eq:CEResult}
    \vis=\frac{\rho}{3}\Delta t\left(\tau-\frac{1}{2}\right).
\end{equation}

\paragraph{Boundary conditions in LBM}
The boundary conditions specified in Section~\ref{subsubsec:FreeFlowAnalysical} are realized in the LBM framework as follows. For the no-slip boundary conditions $\vel=\Vec{0}$ at the vessel walls, as well as on all fully resolved medical device-surfaces are resolved by means of a (simple i.e. with wall-velocity $\vec{0}$) \textit{bounce-back} operator, i.e. total reflection of all populations $f_i$ that would leave the computational fluid domain. The bounce-back operator in a generalized form~\cite{ladd1994numerical, ladd2001lattice, yin2012improved, zhang2012general} can also be used at the inflow-boundary to incorporate the non-zero, e.g. generalized Poiseuille velocity profile there, while at the outlet for the pressure-condition \textit{anti-bounce-back} scheme~\cite{ginzburg2008antiBounceBack, ginzburg2008antiBounceBackPressure, Izquierdo2008antiBounceBack} is employed. For the generalized bounce-back rule, an additional term is introduced with the prescribed velocity value to ensure the correct application of the boundary condition. While the \textit{bounce-back} operator reflects back the population leaving the domain, the \textit{anti-bounce-back} operator introduces antireflection by changing the sign for the corresponding populations. In order to incorporate the given condition properly, as in the case of the generalized \textit{bounce-back} scheme, extra terms are present in terms of the distribution function. 

\paragraph{Non Newtonian extension of the LBM}
Within our short review about the {\LBM} so far,
we have assumed that the viscosity $\vis$, e.g. in \eqref{eq:NavierStokes1} and \eqref{eq:CEResult} is constant.
Switching from a Newtonian to a non-Newtonian model now changes this from a constant to an actually velocity-dependent viscosity,
where \eqref{eq:Casson2} gives an expression for the shear rate $\dot{\gamma}$ evaluated from the Casson model.
However, following~\cite{chopard2002cellular, ouared2005lattice}, the current strain rate can also be directly computed from within the {\LBM} via
\begin{align*}
    \dot{\gamma}_{\textup{LBM}}&=2\sqrt{\tens{\epsilon(\vel)}:\tens{\epsilon(\vel)}},\\
    \tens{\epsilon(\vel}) &= -\frac{3}{2\rho\tau}\,\tens{\Pi}(\rho,\vel),\\
    \tens{\Pi}(\rho,\vel)&=\sum_{i=0}^{n-1}\left(f_i-f_i^{\textup{(eq)}}(\rho,\vel)\right)\,\Vec{c}_i\Vec{c}_i^{\top}.
\end{align*}
Equating these two expressions for $\dot{\gamma}$ into one equation - due to its quadratic nature - allows to solve it~\cite{ouared2005lattice} for the apparent viscosity $\vis$,
which then by \eqref{eq:CEResult} allows to set the relaxation time parameter accordingly for the next time step of the {\LBM},
hence, step by step, adapting viscosity to the current flow situation as described by Casson's law.
A similar approach but with a different Non-Newtonian model, i.e. the Carreau–Yasuda model instead of the Casson-model, is given in~\cite{horvat2024lattice}.

%
%

\subsection{Poro-elastic flow}
\label{subsec:PoroElasticFlow}

To overcome the bottleneck of geometry (re-)construction of the devices
--- either being impossible due to a lack of high-resolution medical images or being tedious to generate and mesh ---
 the theory of deforming {\PM} is applied to analyze the flow behavior and deformation of the {\CA}, the devices and the parent vessel in a homogenized sense.

{\PMs} typically consist of a porous skeleton phase, where the intermediate spaces are filled with one or more fluid phases.
Additionally, one assumes that all pores are interconnected and no enclosed pores exist.
For a fully saturated {\PM}, only a solid phase and one fluid phase is considered. 

The exact and rather complex geometry of the fine pores on the microscale are generally hard to resolve and often unknown.
Luckily, such detailed knowledge at the pore scale is often not required to answer physically relevant questions.
The focus of the theory of deformable {\PM} is to provide an estimate on the macro or continuum scale.
Several techniques as \textit{homogenization}, \textit{volume averaging} or \textit{asymptotic expansion} exists to provide a relation between the different scales and can be applied to obtain a formulation of the quantities of interest on the respective scale. 

In the following a homogenization approach is considered, which smears the fluid and solid phase over the domain and provides averaged quantities.
This results in an overlapped continuum with the domain $\porodomain$,
where fluid and solid phases consist simultaneously with an initial domain $\domaininitial$ and deformed domain $\domaindeformed$.
In fact, some quantities are defined within both phases.
To indicate whether a variable is attributed to the solid or fluid phase of the {\PM},
we use the superscripts $(\bullet)^\solid$ and $(\bullet)^\fluid$, respectively.
Consequently, the solid volume fraction is denoted with $\soliddomain$.
The fluid volume fraction of a domain $\fluiddomain$ is combined via the porosity $\porosity$,
to relate it with the deformed volume: $\porosity~\infinetismald\domaindeformed=\infinetismald\domaindeformedfluid$.

\subsubsection{Governing equations} \label{subsubsec:poroelasteq}

While a detailed derivation of the theory of {\PM} can be found in~\cite{coussyPoromechanics2004} for example,
we just briefly summarize the fundamental system of equations within an isothermal state
and recite the balance equations with a constitutive law for {\PMs}.
Thereby, we base our formulation on~\cite{chapelleGeneralCouplingPorous2014,vuongGeneralApproachModeling2015},
however keep a convective term of the fluid phase.

In the following, we employ the shorthand notation~$\dt{(\bullet)} = \partial (\bullet) / \partial \mytime$ to abbreviate time derivatives.
Moreover, $\dis$, $\vs = \dt{\dis}$, $\defgrad = \defgradterm$ with~$\detdefgrad = \det(\defgrad)$, $\Pktwostress$, $\bodyforcesolid$, and~$\rhos_0$
denote the solid's displacement, velocity, deformation gradient, second Piola-Kirchhoff stress tensor, body force vector and density, respectively.
Accordingly, $\vf$, $\pf$, $\rhof$, $\bodyforcefluid$, $\visfluid$, $\permea$, and~$\fluidstressvisc$
refer to the fluid's velocity, pressure, density, body force, kinematic viscosity, permeability tensor and stress tensor, respectively.
The porosity in the material configuration is denoted by~$\porosity_0$ with its push forward to the spatial configuration reading~$\porosity = \detdefgrad\porosity_0$.
Then, poro-elastic medium is governed by the following coupled system of equations:
\begin{subequations}
\begin{align}
\pcttrans + \pctsolid + \pctfluid = 0  \quad \text{in}~\domaindeformed \times  \kle{\timeinterval}\\
\fpmttrans+\fpmtconv+\fpmtpress-\fpmtbody+\fpmtdarcy-\fpmtvisc=\vecnull \quad \text{in}~\domaindeformed \times  \kle{\timeinterval} \label{eq:poro_fluid_momentum_strong}\\
\spmttrans - \spmtstress - \spmtbody - \spmtpress - \spmtdarcy = \vecnull \quad \text{in}~\domaininitial \times  \kle{\timeinterval} 
\end{align}
accompanied by appropriate initial conditions
\begin{align}
\us =\dsini \quad \text{in}~\domaininitial\times \kle{\initial}   \quad  \frac{\partial \us}{\partial \mytime}= \hatvs_0\quad \text{in}~\domaininitial \times \kle{\initial} \\
\vf=\vfini \quad \text{in}~\domaindeformed \times \kle{\initial}\quad  \porosity=\hatporosityinitial  \quad \text{in}~\domaininitial\times \kle{\initial}  
\end{align} 
and boundary conditions
\begin{align}
\vf\cdot \normalvectorcurr = \hat{v}^\fluid_n  \quad \text{on}~ \boundaryfluiddirichtime \times \kle{\timeinterval} \\
-\pf\Identity\cdot \normalvectorcurr =\tractionfluidneuman  \quad \text{on}~ \boundaryfluidneumantime \times \kle{\timeinterval} \label{eq:PressureBC}\\
\klr{\defgrad\Pktwostress}\cdot \normalvectorinitial = \tractionsolid \quad \text{on}~\boundarysolidneumaninitial \times \kle{\timeinterval} \\
\us = \hat{\us} \quad \text{in}~\boundarysoliddirichinitial \times \kle{\timeinterval}\\
\porosity = \hatporosityinitial \quad \text{in}~\boundaryporositydirichinitial \times \kle{\timeinterval}
\end{align}
on the respective Dirichlet boundary~$\boundarydirichtime = \boundarysoliddirichtime \cup \boundaryfluiddirichtime \cup \boundaryconstrainttime$ and Neumann boundary~$\bcsurfneu = \boundarysolidneumantime \cup \boundaryfluidneumaninitial \cup \boundaryconstrainttime$
with~$\bcsurfdiri \cap \bcsurfneu = \emptyset$ for every time.

\end{subequations}

The deforming {\PM} assumes an incompressible Newtonian fluid,
where the viscous stress tensor~$\fluidstressvisc= 2\visfluid (\grad \vf +\transpose{\klr{\grad \vf}})$ relates the kinematic viscosity~$\visfluid$ and the fluid velocity.
The Darcian contribution within \eqref{eq:poro_fluid_momentum_strong} with the current permeability tensor $\permea$
can be interpreted as a linear resistance acting against the convective fluid flow and causing dissipation.
The influence of viscous dissipation within the fluid through a {\PM} is accounted for by a Brinkman-type contribution.
The difference in comparison to~\cite{vuongGeneralApproachModeling2015} resides in the presence of the convective term.
By including the convective term, the fluid balance equation may simplify to the Navier-Stokes equation for a specific choice of parameters.

The constitutive model of the skeleton is adopted from~\cite{vuongGeneralApproachModeling2015} and depends on the porosity~$\porosity$,
the deformation in form of the Green-Lagrange strain tensor $\GreenLagrangestrain=\greenlagrangeterm$
and the deformation gradient~$\defgrad$.
An additive split of the solid's strain energy function $\strainenergysolidEJporosity$ yields
\begin{align}
\label{eq:PoroStrainEnergyFunctions}
    \strainenergysolidEJporosity = \macrostrainskel +\macrostrainvol + \macrostrainpen.
\end{align}
The macroscopic strain energy function $\macrostrainskel$ governs the constitutive behavior of the porous skeleton.
The contribution~$\macrostrainvol$ accounts for changes due to the compression or expansion of the solid volume
due to the fluid pressure within the pores compared to the initial porosity~$\porosityinitial$, reading
\begin{align}
\macrostrainvol=\macrostrainvolterm.
\end{align}
For a comparison to the linear Biot theory, see~\cite{chapelleGeneralCouplingPorous2014}.
The penalty strain energy function~$\macrostrainpen$ ensures the maintenance of a physical reasonable porosity
\begin{align}
\macrostrainpen=\macrostrainpenterm \text{.}
\end{align}
Consequently, the second Piola-Kirchhoff stress tensor~$\Pktwostress$ can be expressed as
\begin{align}
\Pktwostress= \Pktwostressterm  \quad \text{in}~\domaininitial.
\end{align}
The fluid pressure~$\pf$ is related to the strain energy of the {\PM} via
\begin{align}
\pf=\frac{\partial \strainenergysolidEJporosity}{\partial  \klr{\detdefgrad \porosity}}\label{eq:pf-strain}.
\end{align}
Equation~\eqref{eq:pf-strain} is used as an additional strong equation.

\subsubsection{Weak formulation}
\label{sec:WeakForm}

Following the formulation and implementation of~\cite{vuongGeneralApproachModeling2015},
we use finite elements for the spatial discretization with the function spaces $\mathcal{S}$ for the primary variables $\us$,$\vf$,$\pf$ and $\porosity$ with the according weighting function spaces $\mathcal{V}$, in particular reading:
\begin{align*}
    \mathcal{S}_{\us} &= \left\{ \us \in \mathcal{H}^1\klr{\domain}^{\nsd}  \big| \quad \us = \hat{\us} \quad \text{on}~\boundarysoliddirichinitial \right\} \\
    \mathcal{V}_{\us} &= \left\{ \virtualus \in \mathcal{H}^1\klr{\domaininitial}^{\nsd} \big | \virtualus =\vecnull ~\text{on}~\boundarysoliddirichinitial \right \}\\
    \mathcal{S}_{\vf}&= \left\{ \vf \in \mathcal{H}^1\klr{\domaindeformed}^{\nsd}  \big| \quad \vf \cdot \normalvectorcurr = \hat{v}^\fluid_n \quad \text{on}~\boundaryfluiddirichtime \right\} \\
    \mathcal{V}_{\vf}&= \left\{ \virtualvf \in \mathcal{H}^1\klr{\domaindeformed}^{\nsd} \big | \virtualvf \cdot \normalvectorcurr =0 ~\text{on}~\boundaryfluiddirichtime \right \} \\
    \mathcal{S}_{\pf}&=\left\{ \pf \in  \mathcal{H}^1 \left( \domaindeformed \right) \right\}\\
    \mathcal{V}_{\pf}&=\left\{ \qf \in  \mathcal{H}^1 \left( \domaindeformed \right) \right\}\quad \\
    \mathcal{S}_{\porosity}&= \left\{ \porosity \in \mathcal{H}^1\klr{\domaindeformed} \big| \quad \porosity = \hatporosityinitial \quad \text{on}~\boundaryporositydirichinitial \right\}  \\
    \mathcal{V}_{\porosity}&= \left\{ \virtualporosity \in \mathcal{H}^1\klr{\domaindeformed}^{\nsd} \big | \virtualporosity =0 ~\text{on}~\boundaryporositydirichinitial \right \}
\end{align*}
The according Sobolev spaces are denoted with $\mathcal{H}^1$ with the number of spatial dimensions $\nsd$.
The porosity is chosen as an additional degree of freedom.
By doing this, the higher continuity requirements for the porosity gradient from \eqref{eq:poro_fluid_momentum_strong} is lowered.
Applying now the standard procedure within the finite element method, the variational formulation of the poro-elastic problem can be obtained:
\begin{align*}
\bilinear{\virtualporosity}{\partialfrac{\porosity}{\mytime}}{\domaindeformed} +
\bilinear{\virtualporosity}{\porosity \nabla \cdot \vf }{\domaindeformed} +
\bilinear{\virtualporosity}{ \vf \nabla \porosity }{\domaindeformed} -
\bilinear{\virtualporosity}{\vs \nabla \porosity}{\domaindeformed} = 0 \\
\bilinear {\virtualvf}{\partialfrac{\rhof\vf}{\mytime}}{\domaindeformedfluid} +
\bilinear {\virtualvf}{\fpmtconv}{\domaindeformedfluid} 
- \bilinear {\grad \virtualvf}{\pf}{\domaindeformedfluid}  
-\bilinear {\virtualvf}{\fpmtbody}{\domaindeformedfluid} \\
- \bilinear{\virtualvf}{\tractionfluid}{\boundaryfluidneumantime \cup \boundaryconstraint  } +
\bilinear {\virtualvf}{\fpmtdarcy}{\domaindeformedfluid}
+\bilinear{\grad \virtualvf}{\sigmafvis}{\domaindeformedfluid}
-\bilinear{\virtualvf}{\frac{1}{\porosity}\nabla \porosity \cdot \sigmafvis}{\domaindeformedfluid}
=0
\\
\bilinear{\virtualus}{\rhos_0(1-\porosityinitial) \partialfrac{\vs}{\mytime}}{\domaininitial}
+\bilinear{\virtualE}{\Pktwostress}{\domaininitial}
-\bilinear{\virtualus}{\spmtbody}{\domaininitial}
-\bilinear{\virtualus}{\spmtpress}{\domaininitial} \\
-\bilinear{\virtualus}{\spmtdarcy}{\domaininitial}
-\bilinear{\virtualus}{\tractionsolid}{\boundarysolidneumaninitial}
=0
\\
\bilinear{\qf}{\pf}{\domaininitial}- \bilinear{\qf}{\frac{\partial \strainenergysolidEJporosity}{\partial  \klr{\detdefgrad \porosity}}}{\domaininitial}=0
\end{align*}

\subsection{Discretization and residual-based stabilization}
For the spatial discretization, we use linear finite elements for all fields of primary unknowns,
whereas discretization in time is done via finite differences and employs One-Step-Theta time integration.
To satisfy the inf-sup condition for~$\vf$ and~$\pf$,
we employ residual-based stabilization terms,
namely {\PSPG}, {\SUPG}, and a stabilization for the Darcy term.
The stabilization for the {\PSPG} and Darcy term are motivated from a residual-based variational multi-scale decomposition from~\cite{badiaStabilizedContinuousDiscontinuous2010} and can be stated as
\begin{align*}
    \PSPGterm
\end{align*}
with~$\discreteresidualPFM$ denoting the discrete strong form residual of \eqref{eq:poro_fluid_momentum_strong}.
and~$\taumstab$ being a user-chosen stabilization parameter.
The stabilization for to the Darcy term is given as
\begin{align*}
    \PSPGDarcyterm.
\end{align*}
To stabilize the nonlinear convective term, we add a {\SUPG} contribution reading
\begin{align*}
\sum^{\nele}_{e=1}\bilinear{\klr{\vf \cdot \grad} \virtualvf}{\taumstab\discreteresidualPFM}{\Omega_e}.
\end{align*}
The necessary stabilization parameters are determined by a combination of ideas in~\cite{Taylor1998} and~\cite{badiaStabilizedContinuousDiscontinuous2010}.
We note that such stabilization techniques have been around for several years and proved to be very effective.
Yet, this particular combination of stabilization methods and terms poses a novelty.
The final implementation of the numerical approach just presented was accomplished in the collaborative software project \textit{BACI}~\cite{baci}\footnote{An LBM implementation of a similar, even though only porous not poro-elastic approach is discussed in~\cite{horvat2024lattice} in terms of the volume-averaged Navier-Stokes equations (VANSE) approach.}.

\section{Modeling of endovascular devices} \label{sec:DeviceModelling}

\begin{figure}
\centering
\subfigure[Simulated coiling device model within narrow-necked {\CA}, see also Figure~\ref{fig:CoilModel}.]{\includegraphics[width=4.5cm]{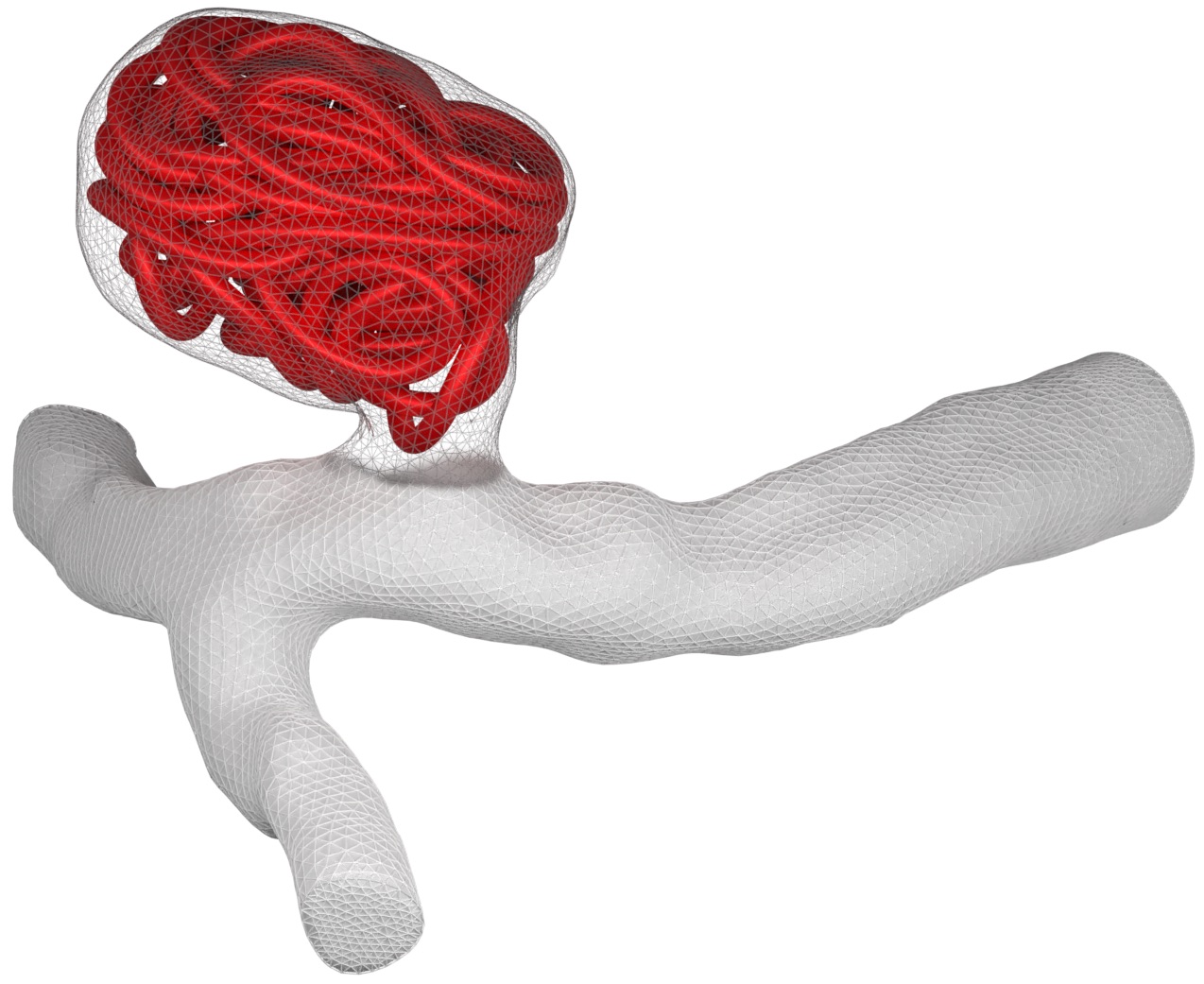}}\hspace*{5mm}
\subfigure[Geometric {\web} device model within wide-necked {\CA}, see also Figure~\ref{fig:WebInsertion}.]{\includegraphics[width=5.5cm]{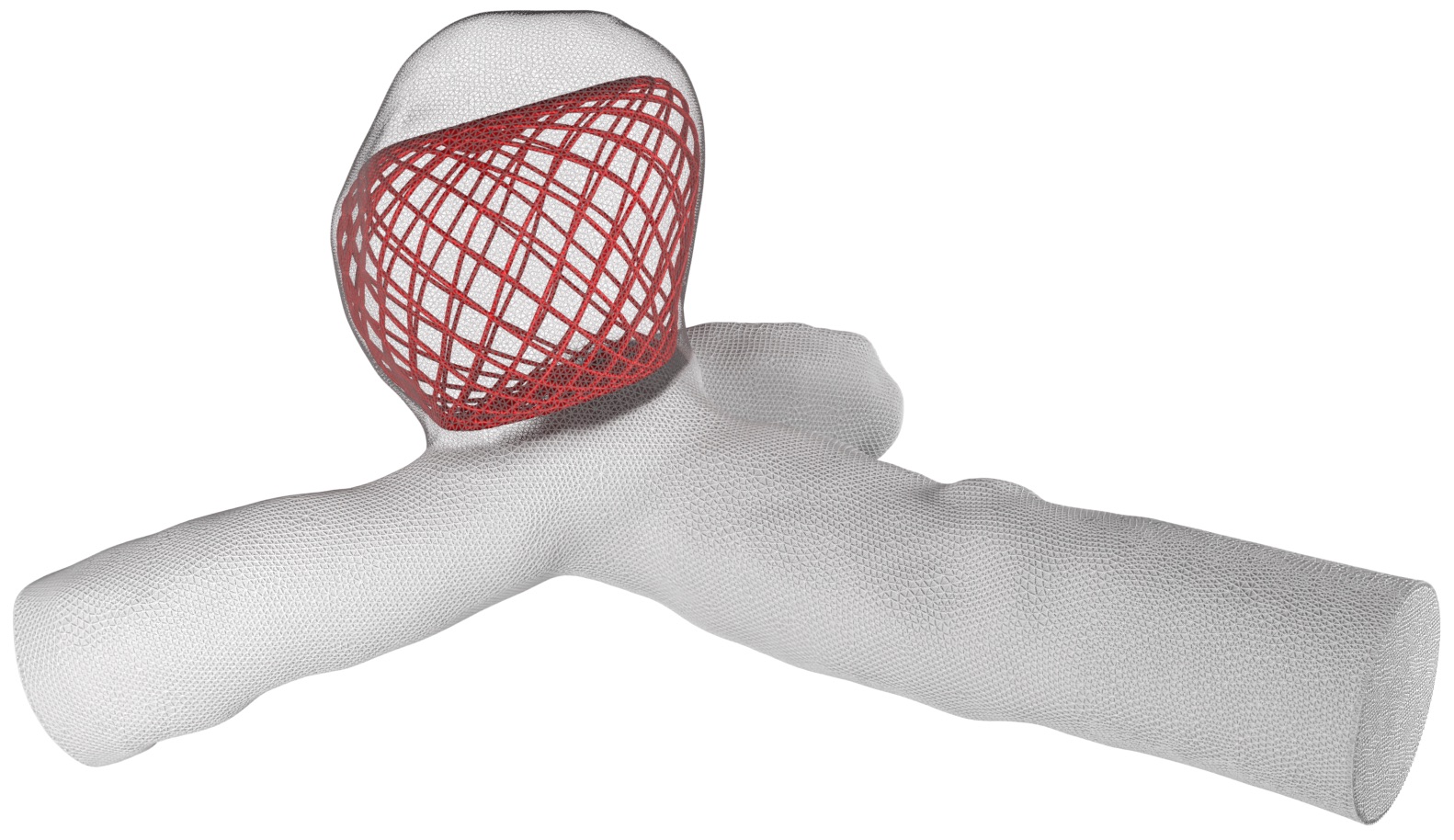}}\hspace*{5mm}
\subfigure[Geometric stent model even covering two {\CAs} directly above it.]{\includegraphics[width=5.5cm]{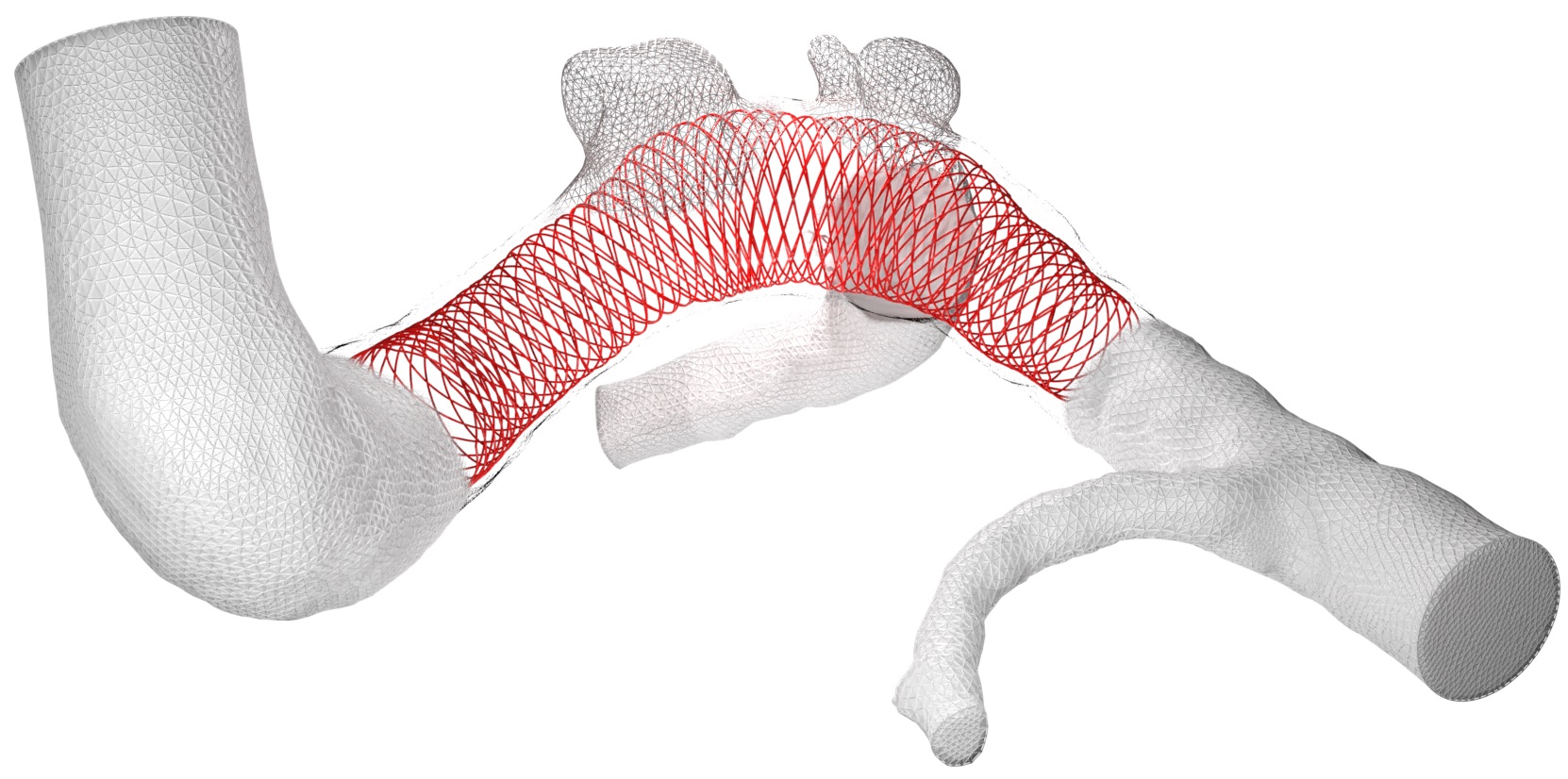}}
\caption{Virtual model examples for the three classes of endovascular devices\label{fig:Devices}}
\end{figure}

In this section we will discuss mathematical models, both mechanical as well as purely geometrical ones, for different kinds of aneurysm treatment devices that are used to model their accurate shape and placement within the vessel geometries obtained in Section~\ref{subsec:GeometryPreprocessing}. As already mentioned in Section~\ref{subsec:FreeFlow} they can then also be included in the hemodynamical simulations, see Section~\ref{sec:NumericalSimulations}, as fully resolved flow-obstacles to analyze their effect on the blood-flow \textit{in} resp. \textit{into} the aneurysm. It should be noted that the upcoming device-(placement) simulations are all still conducted \textit{in absence} of the blood-stream hence neglecting its effects on the device shapes but rather have to be seen as separate simulations or procedures simply to \textit{generate} realistic device shapes in preparation of the final hemodynamics simulations, where these devices are then used as (still static) obstacles.\\
\indent Based on the medical background from Section~\ref{sec:MedicalRelevanceAndChallenges} the three kinds of devices we will consider are endovascular coils in Section~\ref{subsec:Coils}, {\web} devices in Section~\ref{subsec:WebDevices} and stents/flow diverters in Section~\ref{subsec:Stents}.
Figure~\ref{fig:Devices} anticipates images of virtual representation of these devices.

\subsection{Coiling devices} \label{subsec:Coils}
\begin{figure}
\centering
\hfill
\subfigure[Schematic micro structure of a coiling wire: $D_3$ denotes the approximate diameter of the coil-node forming within the aneurysm,
$D_2$ the apparent diameter of the (macroscopic) coiling wire
and~$D_1$ the diameter of the actually microscopic wire wrapped in a helical way to form the macroscopic~($D_2$) wire.]{\label{fig:CoilPhotosSketch}\hspace{3em}\includegraphics[height=4cm]{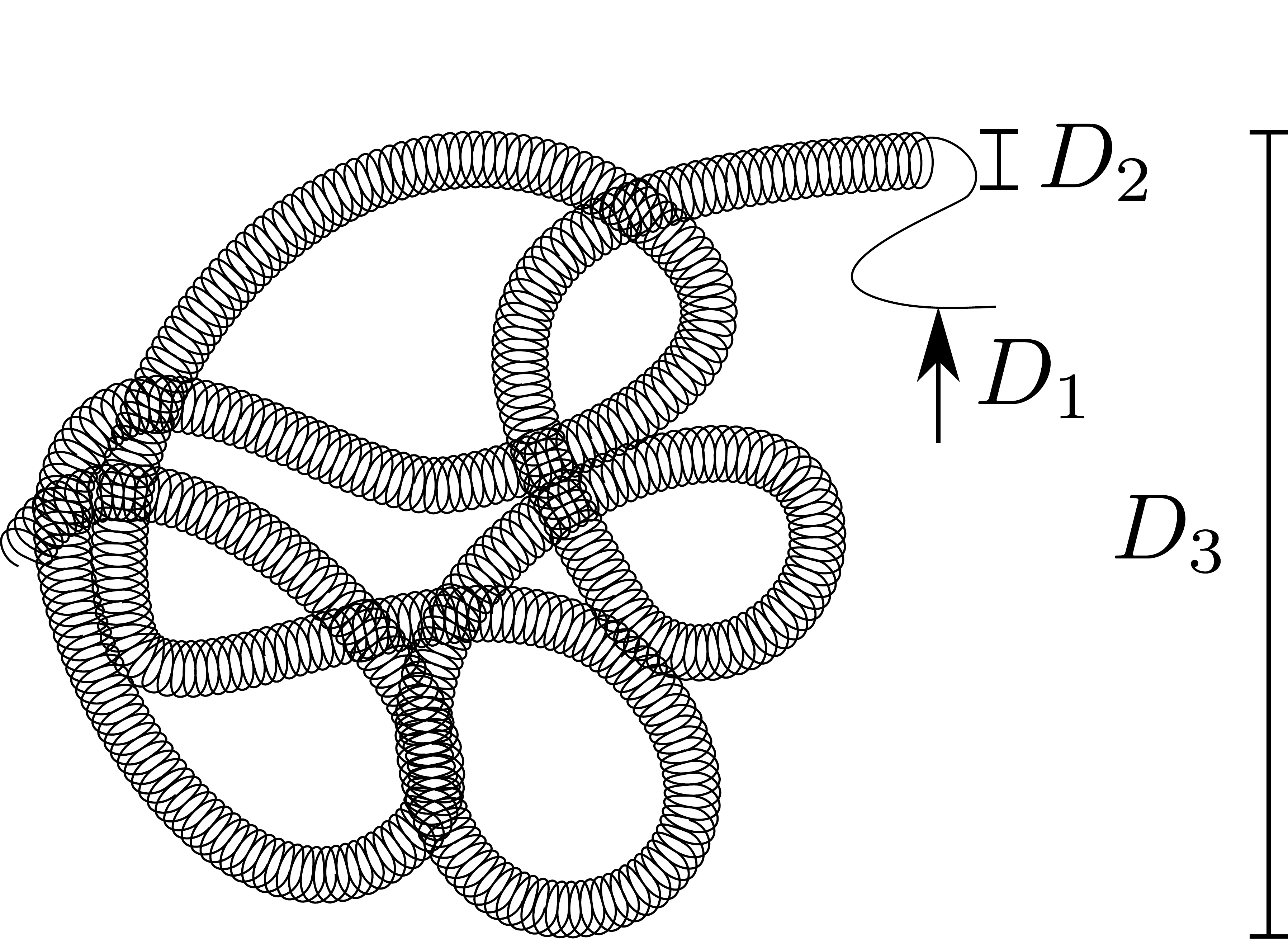}\hspace{3em}}
\hfill
\subfigure[Expulsion of an actual coil into regular surrounding (air) with contactless loop formation and curling]{\label{fig:CoilPhotosExpulsion}\includegraphics[height=3cm]{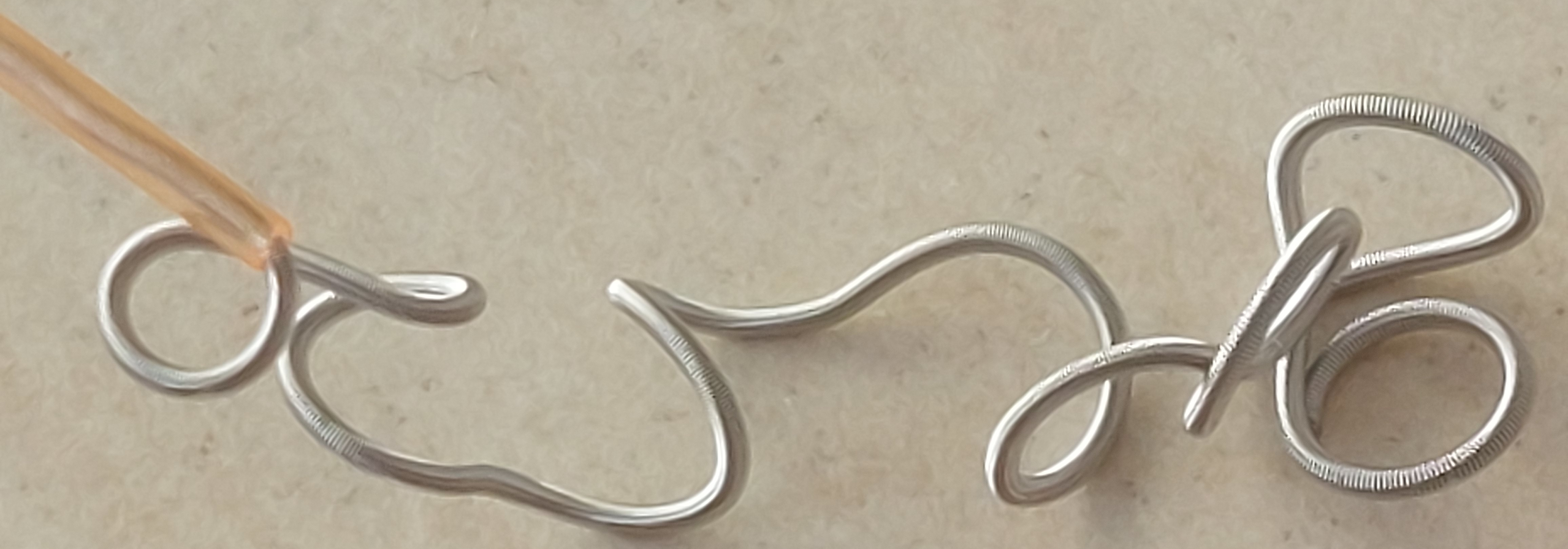}}
\hfill~
\caption{Sketch and photo of coiling wires}
\label{fig:CoilPhotos}
\end{figure}

First, we focus on endovascular coils~\cite{eddleman2013endovascular, holzberger2024comprehensive, hui2014history}.
Figure~\ref{fig:CoilPhotosExpulsion} shows a coil emerging from a micro-catheter onto a stable surface in unconfined air.
For coil deployment into the {\CA} cavity, the micro-catheter is first placed into the aneurysm's parent artery.
Then, coils can then be pushed through the catheter and into the aneurysm with the goal of filling it and occluding it from the blood flow in the parent vessel.
The initial phase of this process involves the insertion of so-called ``framing coils''.
These have an intrinsic spherical shape to form a basket-like structure, when being placed inside the aneurysm.
Then, the aneurysm is treated through the placement of ``filling coils'', featuring helical loops,
or ``finishing coils'' which consist of irregular wire loops.
Both types share the attribute of being much softer than the farming coil~\cite{ito2018experimental, kanenaka2016comparative}. 
In this way,
the procedure ensures full volumetric occlusion up to a desired global packing density $c_g=V_{\textup{coil}}/V_{\textup{aneurysm}}$.
In clinical practice, the achieved packing density usually varies between 30-40\,\%.
It can be computed using the measures obtained in Section~\ref{subsec:Measures}. 
Figure~\ref{fig:CoilPhotosSketch} depicts the main schematics of a coiling wire.
The stock wire with radius~$D_1$ undergoes a sequence of transformations.
Initially, it adapts into a secondary helical structure with radius~$D_2$.
Then, this secondary structure is transformed into a tertiary structure of radius~$D_3$.
The latter two structures are imprinted into the stock wire by the manufacturing process.
Depending on the manufacturer, the dimensions of~$D_1$ approximately exhibit a range from~$0.050 \,\si{\milli\meter}$ to~$0.076 \,\si{\milli\meter}$. Correspondingly, the secondary structure diameter~$D_2$ might vary between~$0.254 \,\si{\milli\meter}$ and~$0.381 \si{\milli\meter}$~\cite{white2008coils}.
When considering tertiary structures, framing coils can differ from filling coils significantly.
For the former, $D_3$ is often approximately the size of the aneurysm, while a much smaller size is desired for the latter.
Therefore, the diameter~$D_3$ can vary from~$2 \,\si{\milli\meter}$ to~$20 \,\si{\milli\meter}$~\cite{white2008coils}.
The stock wire of the majority of contemporary coil designs is composed of a platinum alloy (see Figure~\ref{fig:CoilPhotos}). 
Arguably one of the most significant parameters in coil design is the bending stiffness often relabeled as the aforementioned softness,
given by~\cite{white2008coils}
\begin{align}
    k_b= \frac{D_1^4 G}{8 D_2^3 n_p},
\end{align}
where~$G$ is the shear modulus of the stock wire (which is approximately~$82\,\si{\giga\pascal}$ for platinum), and~$n_p$ represents the pitch, denoting the space between two stock wire loops within the secondary helical structure~$D_2$, often being negligible or very small. Framing coils generally have a higher~$k_b$ value to retain their spherical structure, while filling coils possess a lower $k_b$ value. This bending stiffness together with the imprinted shape migrates the risk of poking the aneurysm wall, preventing rupture. The impact of the bending stiffness and imprinted shape can be seen very well, when the coil is for example protruded just into regular surrounding, when it already starts to curl up on its own -- even without contact (see Figure~\ref{fig:CoilPhotosExpulsion}).
This behavior is accounted for in the mechanical coiling model by letting the ``natural'' shape of the coil be curved.
Meaning that when no strain due to stretching, bending and torsion is present in the coil,
it attains a curved shape that is prescribed by the manufacturing process,
as for example shown in~\cite{wallace2001coildesigns}.
For the framing coil, this would be the aforementioned basket-like wall-adhering structure which in reality and numerics -- due to (blunt) contact with the wall and the wire itself -- of course might deviate from the idealistic shape.
For a mathematical description of the coil insertion process, we base our work on the \textit{Discrete Elastic Rods} formulation from~\cite{bergou2008discrete, bergou2010discrete}, which is a method for discretizing Kirchhoff rods. Please refer to \cite{holzberger2024comprehensive} for a more detailed account on the implementation as well as validation and application of the coiling model together with its mathematical derivation and backgrounds. In summary, within the model the coil is considered as a discrete one-dimensional space curve consisting of the~$n$ material points~$\Vec{x}_0,...,\Vec{x}_{n-1}\in\mathds{R}^3$ that are connected by the~$n-1$ edges~$\Vec{e}^j=\Vec{x}_{j+1}-\Vec{x}_j$. Each edge is equipped with an orthonormal basis of material directors~$\big[ \Vec{d}_1^j, \Vec{d}_2^j, \Vec{d}_3^j \big]\in\operatorname{SO}(3)$, indicating the orientation of the edges cross section. It is required that the third material director fulfills~$\Vec{d}_3^j=\Vec{e}^j/\| \Vec{e}^j \|$ for~$j=0,...,n-2$, meaning that it is adapted to the tangent of the curve. These assumptions let us formulate the elastic strain energy of the curve
\begin{align}
    E = \underbrace{\frac{1}{2}\sum\limits_{j=0}^{n-2} \alpha_s(\varepsilon^j)^2\|\bar{\Vec{e}}^j\|}_{\text{stretching energy}} + 
    \underbrace{\frac{1}{2}\sum_{i=1}^{n-2}\frac{\beta_t}{\bar{l}_i}(\vartheta_i-\bar{\vartheta}_i)^2}_{\text{torsion energy}} + 
    \underbrace{\frac{1}{2} \sum_{i=1}^{n-2} \frac{\beta_b}{\bar{l}_i}(\Vec{\kappa}_i-\bar{\Vec{\kappa}}_i)^2}_{\text{bending energy}},
\end{align}
where the first term refers to the axial strain energy due to the relative stretch~$\varepsilon^j= \|\Vec{e}^j\| / \|\bar{\Vec{e}}^j\|-1$ on the edge~$j$, the second term refers to the torsion energy containing the discrete integrated twist~$\vartheta_i$ at the material point~$\Vec{x}_i$ and the last term refers to the bending energy containing the discrete integrated curvature vector~$\Vec{\kappa}_i$ measured as well on the material point $\Vec{x}_i$. The quantities~$\vartheta_i$, $\Vec{\kappa}_i$ can be directly obtained from the curve as described in~\cite{bergou2008discrete}. Bared quantities~$\bar{\Vec{e}}^j, \bar{\vartheta}_i, \bar{\Vec{\kappa}}_i$ refer to the naturally curved shape of the rod and are therefore constant and~$\bar{l}_i =(\| \bar{\Vec{e}}^i \| + \| \bar{\Vec{e}}^{i+1} \|)/2$. The internal forces acting on a material point~$\Vec{x}_i$ can then be written as $\Vec{P}_{\textup{int},i} = -\partial E / \partial \Vec{x}_i$.  
Having an expression for the internal forces~$\Vec{P}_{\textup{int},i}$ due to stretching and letting the external forces due to collisions be~$\Vec{P}_{\textup{ext},i}$ (that are described in the following) acting on node~$\Vec{x}_i$, by Newton's second law,
one can describe the dynamics of the coil as the second order system of ordinary differential equations
\begin{align*}
m_i\Ddot{\vec{x}}_i+\nu\dot{\Vec{x}}_i=\Vec{P}_{\textup{int},i} + \Vec{P}_{\textup{ext},i},\quad \forall i=0,1,\dots,n-1,
\end{align*}
that we solve numerically by the semi-implicit Euler method. Therein, $\nu$ refers to a global damping parameter, which can be interpreted to depend on the surrounding fluid within the aneurysm damping the insertion process. For the  external forces~$\Vec{P}_{ext,i}$ due to coil-coil and coil-wall contacts,
we adapt the formulation from~\cite{gazzola2018forward} with slight simplifications.
Coil-coil contacts are modeled by the repulsive force on the material point~$\Vec{x}_i$ as
\begin{align}
    \Vec{P}_{cc,i} = \sum\limits_{j=0}^{n-2} H(\epsilon_{ij})\left(-k_{cc}\epsilon_{ij}-\gamma_{cc}\left(\dot{\Vec{x}}_i - \frac{\dot{\Vec{x}}_{j}+\dot{\Vec{x}}_{j+1}}{2}\right)\cdot \Vec{d}^{ij}_{min}\right)\Vec{d}^{ij}_{min},
\end{align}
that is activated by the Heaviside function~$H$ when the scalar overlap~$\epsilon_{ij}$ between material point~$\Vec{x}_i$ and edge~$\Vec{e}^j$ is positive. One defines the scalar overlap with help of the minimum distance vector~$\Vec{d}^{ij}_{min}$ from~$\Vec{x}_i$ to edge~$\Vec{e}^j$ as~$\epsilon_{ij}=D_2 - \|\Vec{d}^{ij}_{min}\|$. For the coil-wall contact, a Coulomb friction model is assumed. The wall force on the material point~$\Vec{x}_i$ is decomposed into tangential and normal component~$\Vec{P}_{cw,i}(\Vec{P}_i)=\Vec{P}_{cw,i}(\Vec{P}_i)_{\perp} \oplus \Vec{P}_{cw,i}(\Vec{P}_i)_{||}$, that both depend on the force~$\Vec{P}_i$ acting on the material point. Starting with the tangential component one has
\begin{align}
    \Vec{P}_{cw,i}(\Vec{P}_i)_{||} = -\mu_w \|\Vec{P}_{i,\perp}\|  \frac{\dot{\Vec{x}}_{i,||}}{\| \dot{\Vec{x}}_{i,||} \|},
\end{align}
for that we decompose the nodal forces and velocities into their tangent and normal components with respect to the wall normal~$\Vec{P}_i = \Vec{P}_{i,\perp} \oplus \Vec{P}_{i,||}$ and~$\dot{\Vec{x}}_i = \dot{\Vec{x}}_{i,\perp} \oplus \dot{\Vec{x}}_{i,||}$. The parameter~$\mu_w$ is the friction coefficient for sliding coil-wall contact. Accordingly, the normal contact force is  
\begin{align}
    \Vec{P}_{cw,i}(\Vec{P}_i)_{\perp}  = H(\epsilon) (-\|\Vec{P}_{\perp,i}\| + k_{cw}\epsilon - \gamma_w \dot{\Vec{x}}_i \cdot \Vec{n}_w)\Vec{n}_w.
\end{align}
Similar as for the coil-coil contacts,
the activation of the normal force is triggered by the Heaviside function with~$\epsilon=D_2/2-d_{min,w}^i$ and~$d_{min,w}^i$ denoting the smallest distance from~$\Vec{x}_i$ to the wall.
The parameters~$k_{cw}$ and~$\gamma_w$ refer to the wall stiffness and dissipation coefficients, respectively.
$\Vec{n}_w$ denotes the wall outward unit normal at the position that is closest to~$\Vec{x}_i$.
Finally, the external forces on the coil are given by~$\Vec{P}_{\textup{ext},i} = \Vec{P}_{cc,i} + \Vec{P}_{cw,i}(\Vec{P}_{int,i} + \Vec{P}_{cc,i})$,
which completes the discussion of the contact model. 
Since not all nodes are within the aneurysm from the beginning of the simulation,
but are rather inserted one after another,
we supply a ``flagging'' data structure to the wire nodes, such that all the previously mentioned computations are only conducted for those already inserted, while ``exterior'' nodes are equipped with ``boundary-conditions'' of a predefined ``insertion velocity''~$\dot{\Vec{x}}_i=\vel_{\textup{ins}}$.
Note that the model stated here neglects the rotational moment acting on the coil,
which allows us to drop the rotational degree of freedom and, thus, making the systems dynamics depend only on the material point forces.
\begin{figure}
\centering
\begin{minipage}{\textwidth}
    \includegraphics[width=0.24\textwidth]{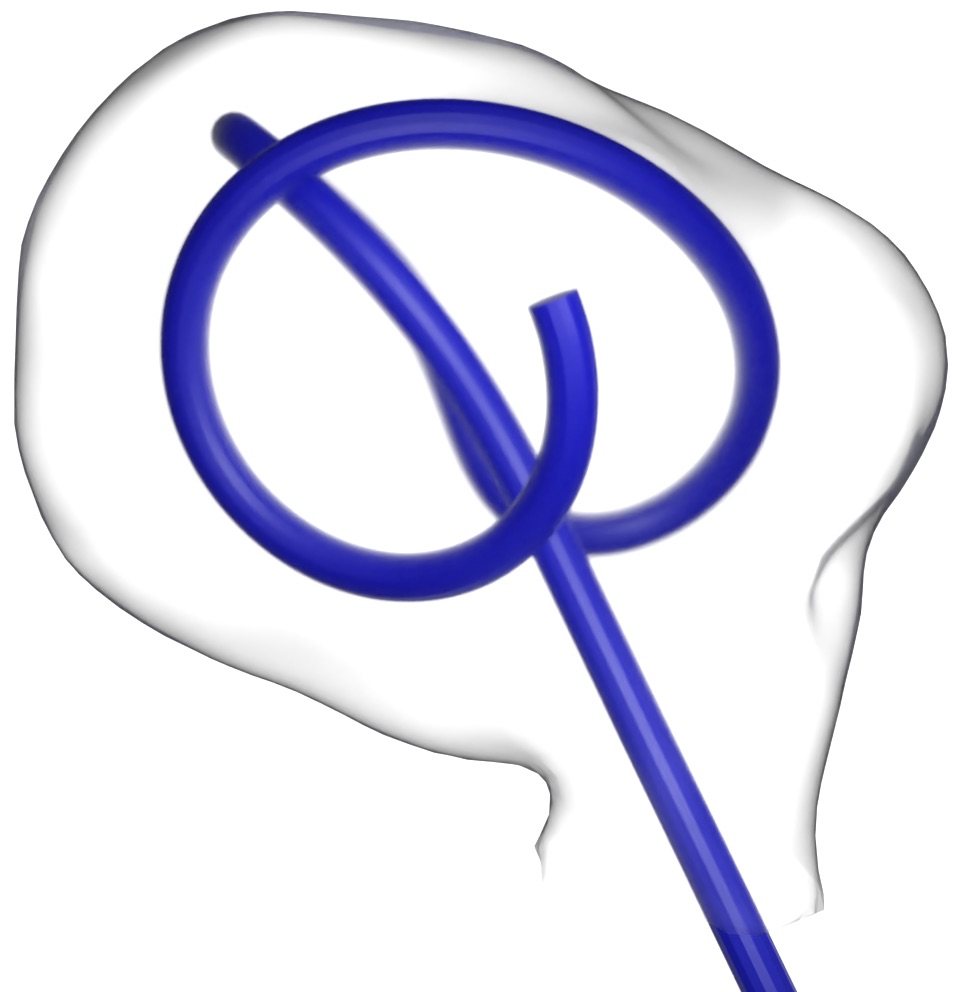}
    \includegraphics[width=0.24\textwidth]{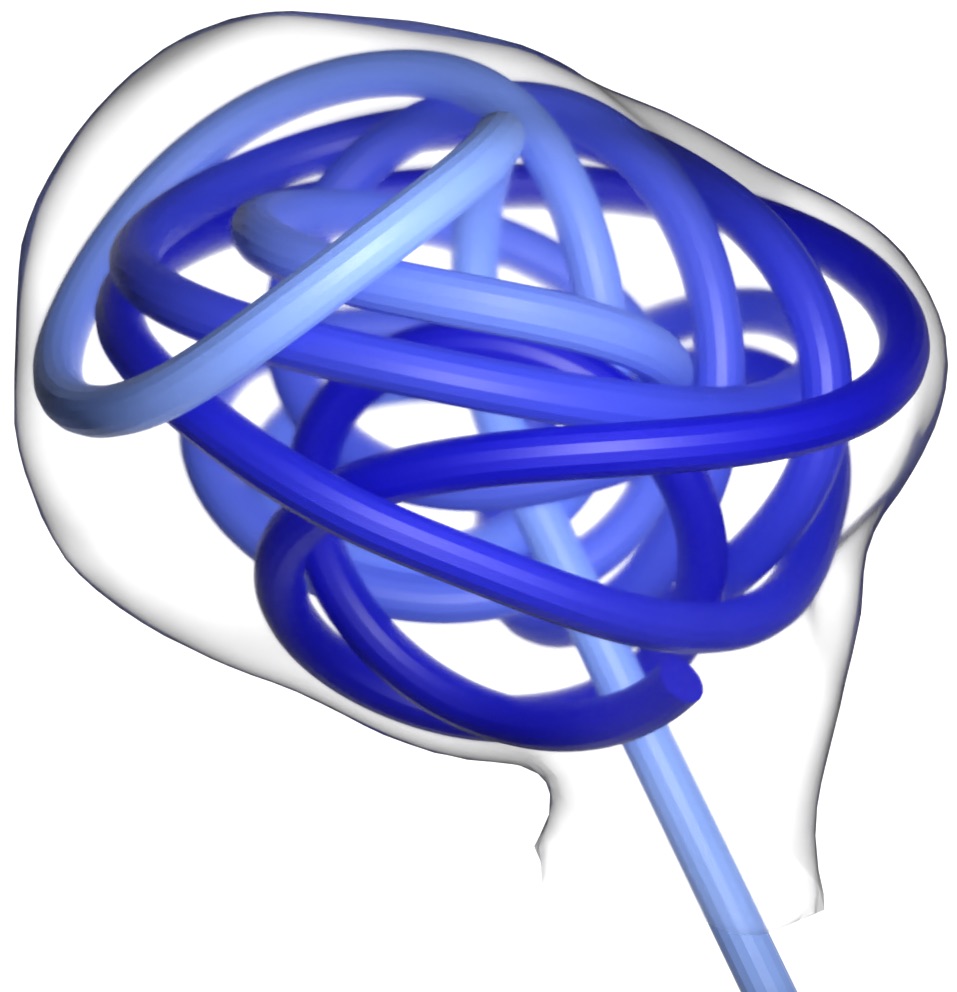}
    \includegraphics[width=0.24\textwidth]{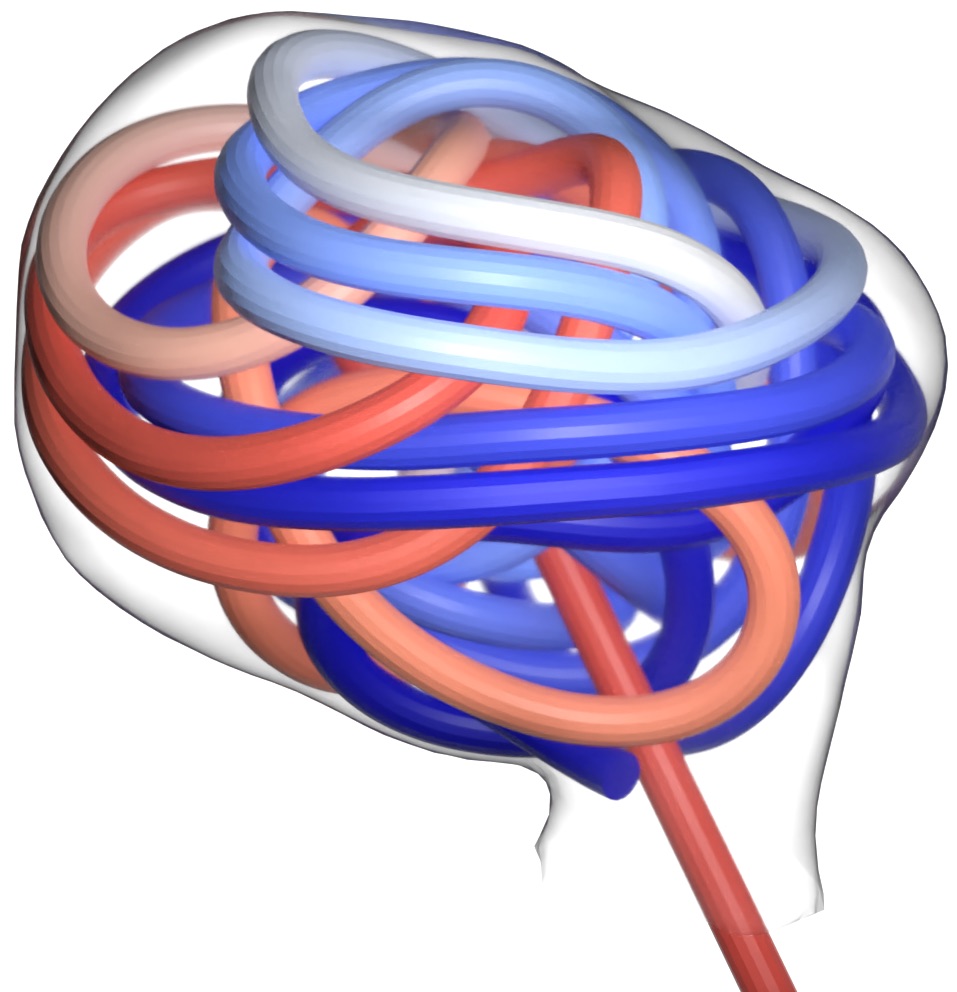}
    \includegraphics[width=0.24\textwidth]{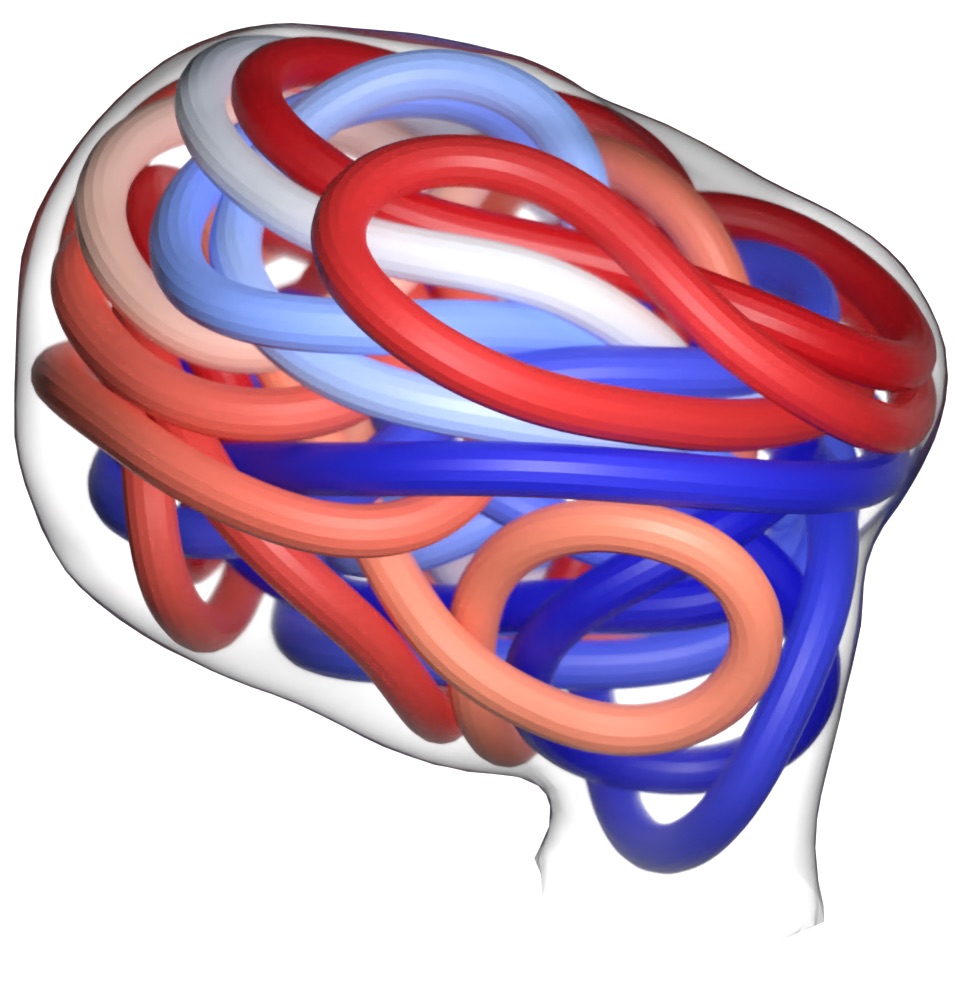}
\end{minipage}\\
\begin{minipage}{\textwidth}
\includegraphics[width=0.58\textwidth]{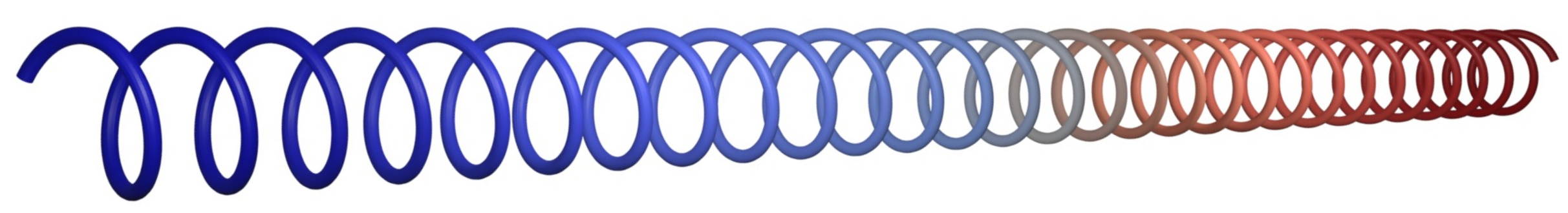}
\includegraphics[width=0.38\textwidth]{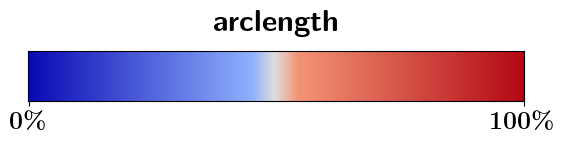}
\end{minipage}
\caption{\emph{From left to right:} Insertion of a coil into an intracranial aneurysm. \emph{Bottom left:} The natural curvature of the coil is extracted from a helix. Color coding is with respect to the coil length. The final volumetric packing density in this example is approximately $33.7\,\%$.}
\label{fig:CoilModel}
\end{figure}
In this work, a helix coil is used (see Figure~\ref{fig:CoilModel} bottom left for the naturally curved shape).
Parameters and boundary conditions for the simulation of the coil insertion process can be found in Table~\ref{fig:CoilPramameters}.
A sequence of snapshots of the coil insertion process can be seen in Figure~\ref{fig:CoilModel}.
As expected, the coil is forming into a curved shape since we have imposed a helix onto its natural shape.
The more coil is inserted, the more densely it is packed inside the aneurysm.
To model a realistic occlusion from the blood flow,
up to~$30-40\,\si{\percent}$ of the aneurysm volume should be be occupied by the coil, which in our case is $33.7\,\si{\percent}$.
After having placed the coil into the aneurysm,
the time step size is decreased to~$10\,\si{\percent}$ of the original time step size and the micro-catheter is removed,
allowing the coil to settle in its final equilibrium-position.

\begin{table}
\centering
\caption{Parameters used in the simulation of the helix coil.}
\label{fig:CoilPramameters}
\begin{tabular}{|c|c|c|}
\hline
Parameter & Value & Unit \\
\hline
Coil diameter $D_2$ & 0.45 & \si{\milli\meter} \\
Helix diameter $D_3$ & 4 & \si{\milli\meter} \\
Helix wavelength $P_H$ & $2\pi D_2$ & \si{\milli\meter} \\
Coil total length $L$ & 37.8 & \si{\centi\meter} \\
\hline
Stretching constant $\alpha_s$ & 18 & \si{\joule\per\meter} \\
Bending constant $\beta_b$ & $1\cdot 10^{-3}$ & \si{\joule} \\
Torsion constant $\beta_t$ & $73\cdot 10^{-5}$ & \si{\joule} \\
\hline
Coil-coil contact stiffness $k_{cc}$ & $1\cdot10^{4}$ & \si{\kilogram\per\squared\second} \\
Coil-coil contact dissipation $\gamma_{cc}$ & 1 & \si{\kilogram\per\second} \\
Coil-wall contact stiffness $k_{wc}$ & $2\cdot 10^{4}$ & \si{\kilo\newton\per \square\meter} \\
Coil-wall contact dissipation $\gamma_{wc}$ & $1\cdot 10^{-4}$ & \si{\kilo\newton\per\meter} \\
Coil-wall slip coefficient $\mu_{w}$ & 0.2 &  \\
\hline
Material points $n$ & 1682 &  \\
Time step-size $\Delta t$ & 1 & \si{\micro\second} \\
Velocity dissipation parameter $\nu$ & 0.2 & \si{\kilogram\per\second}  \\
Insertion velocity $\vel_{ins}$ & 3 & \si{\centi\meter\per\second} \\
\hline
\end{tabular}
\end{table}

\subsection{Woven Endo-Bridge (WEB) devices} \label{subsec:WebDevices}
{\web} devices~\cite{goyal2020web, pierot2015web} are also deployed to the aneurysm via a catheter,
from which they are protruded into the sac of the aneurysm.
In contrast to the comparably simple one-dimensional structure of the coiling wire,
the {\web} device unfolds into an intricate net shape on a two-dimensional (up to thread thickness) spherical or box-like manifold,
see Figure~\ref{fig:WebModel}, that, if sized properly, fits and gets stuck in the aneurysm by adapting to its shape and with slight tension.
The device is supposed to serve the same purpose as a coiling wire, i.e. to occlude the aneurysm sac from further blood flow.
However, it reaches this goal not by a full \textit{volumetric} occlusion but rather by hindering inflow to the aneurysm by the devices very narrow wire-spacing acting like a fluid cage.
Since a full mechanical simulation of the {\web} device unfolding process lies beyond the scope of this publication,
we pursue a purely geometric, parameterized shape model for the {\web} device.
This device model can then be deformed in order to fit as precise as possible into a given aneurysm geometry to become an obstacle for the hemodynamical simulation,
like the final shape of the coil from before, in order to analyze its effect on the blood flow behavior.

\subsubsection{Mathematical shape design of {\web} devices}
\begin{figure}
\begin{center}
\subfigure[Influence of different shape-parameters on the design function $z_d(\cdot)$]{\includegraphics[width=10cm]{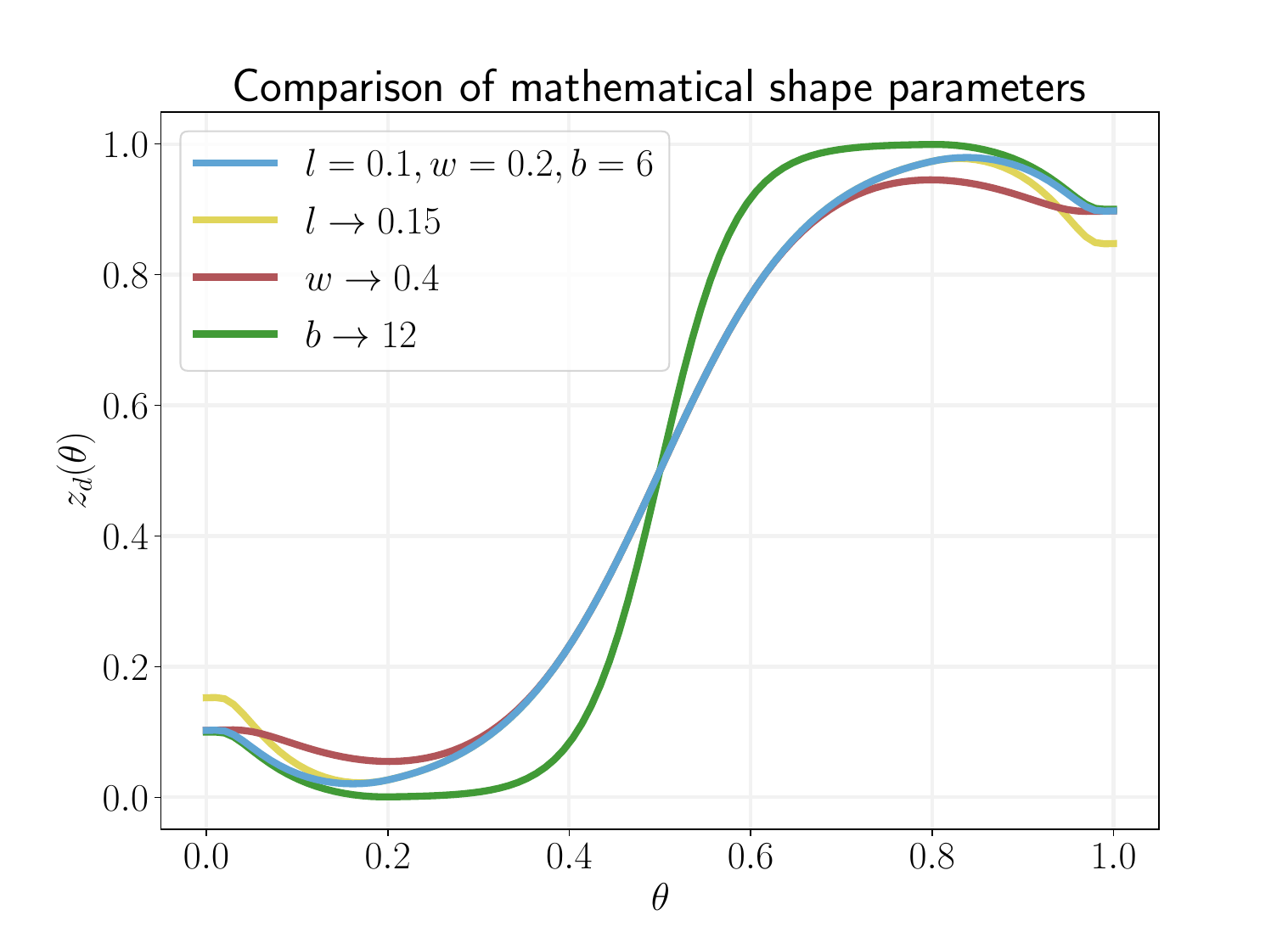}}\subfigure[A more realistic {\web} device w.r.t. thread-count $n_t=70$ and -thickness $r_t=0.0055$, which, however, requires an extremely high mesh resolution. The remaining parameters are the same as in the last (green) case of subfigure (c).]{\includegraphics[width=7cm]{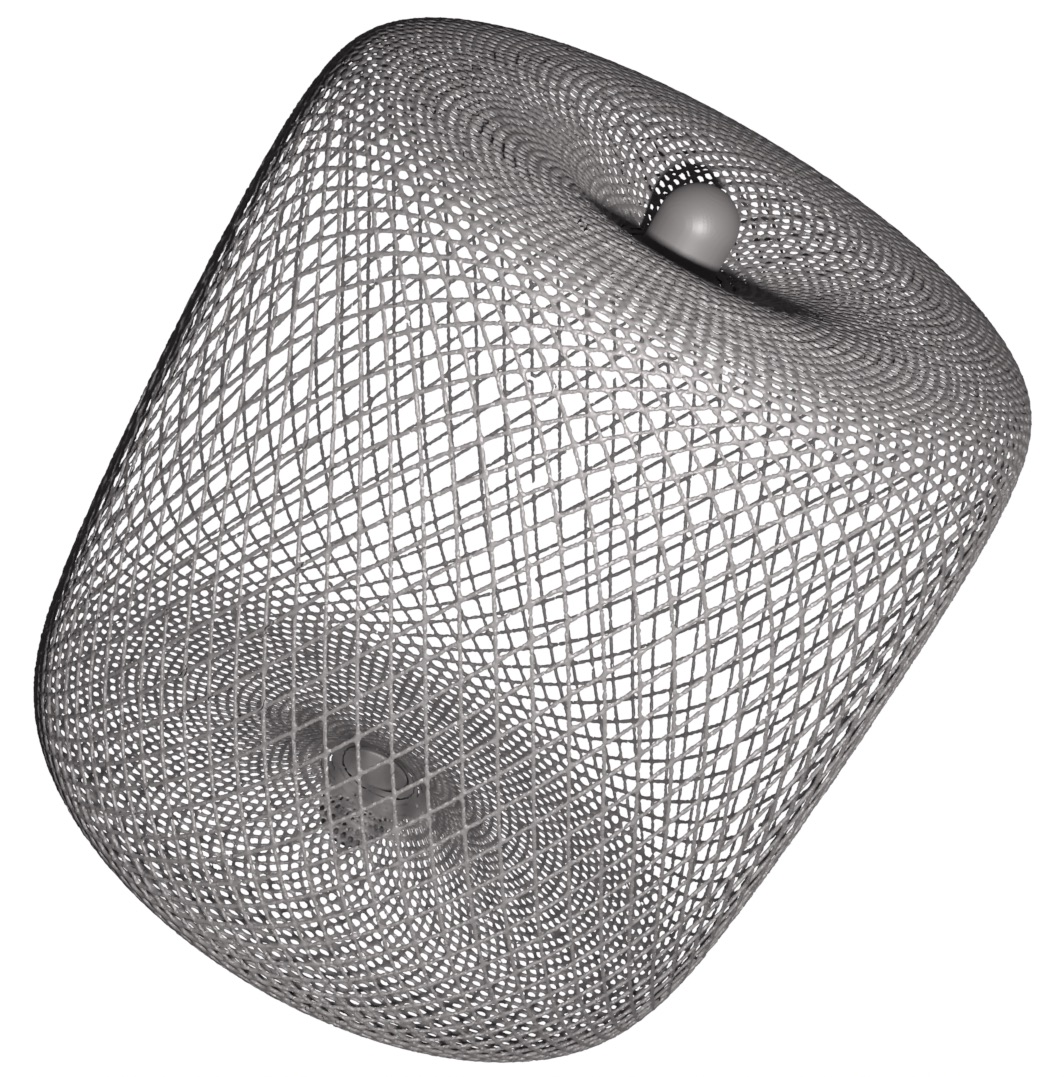}}\\
\subfigure[Differently shaped {\web} devices according to the mathematical design parameters used in subfigure (a) (left to right matches the legend entries from top to bottom in subfigure (a)). They all use $n_t=20$ thread-pairs, $n_w=1, h_{d,\max}=2, r_{d,\max}=1$ and $r_t=0.02$.]{\includegraphics[width=4.25cm]{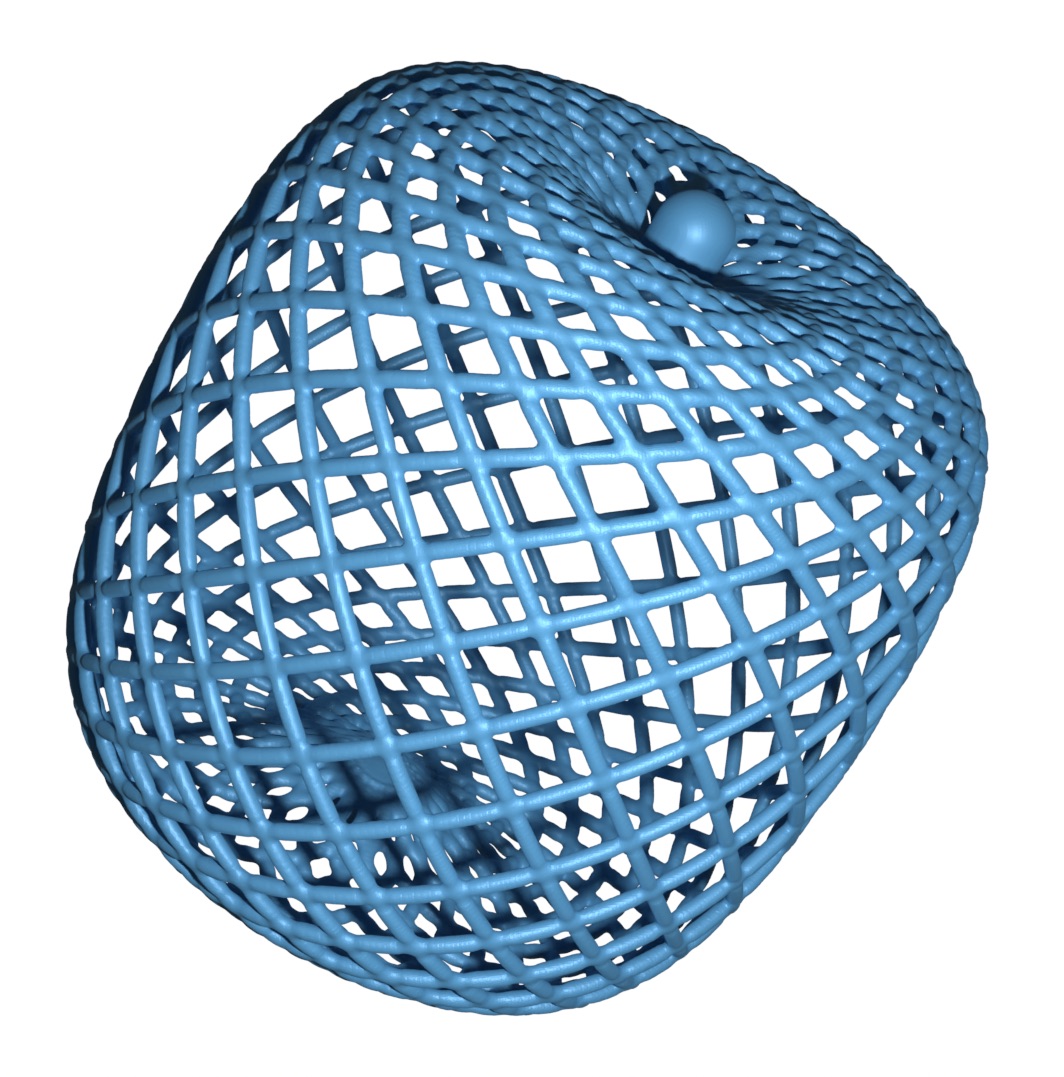}\includegraphics[width=4.25cm]{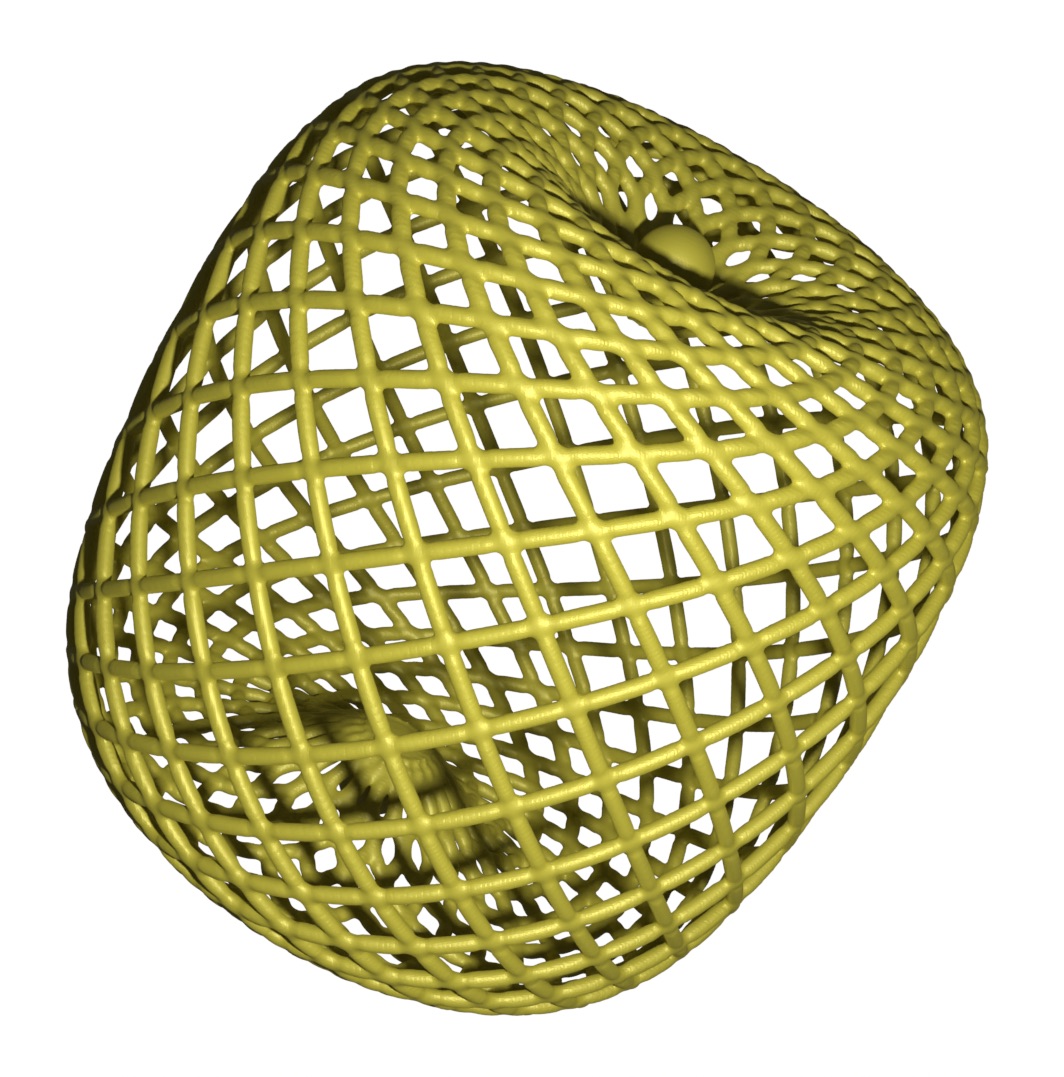}\includegraphics[width=4.25cm]{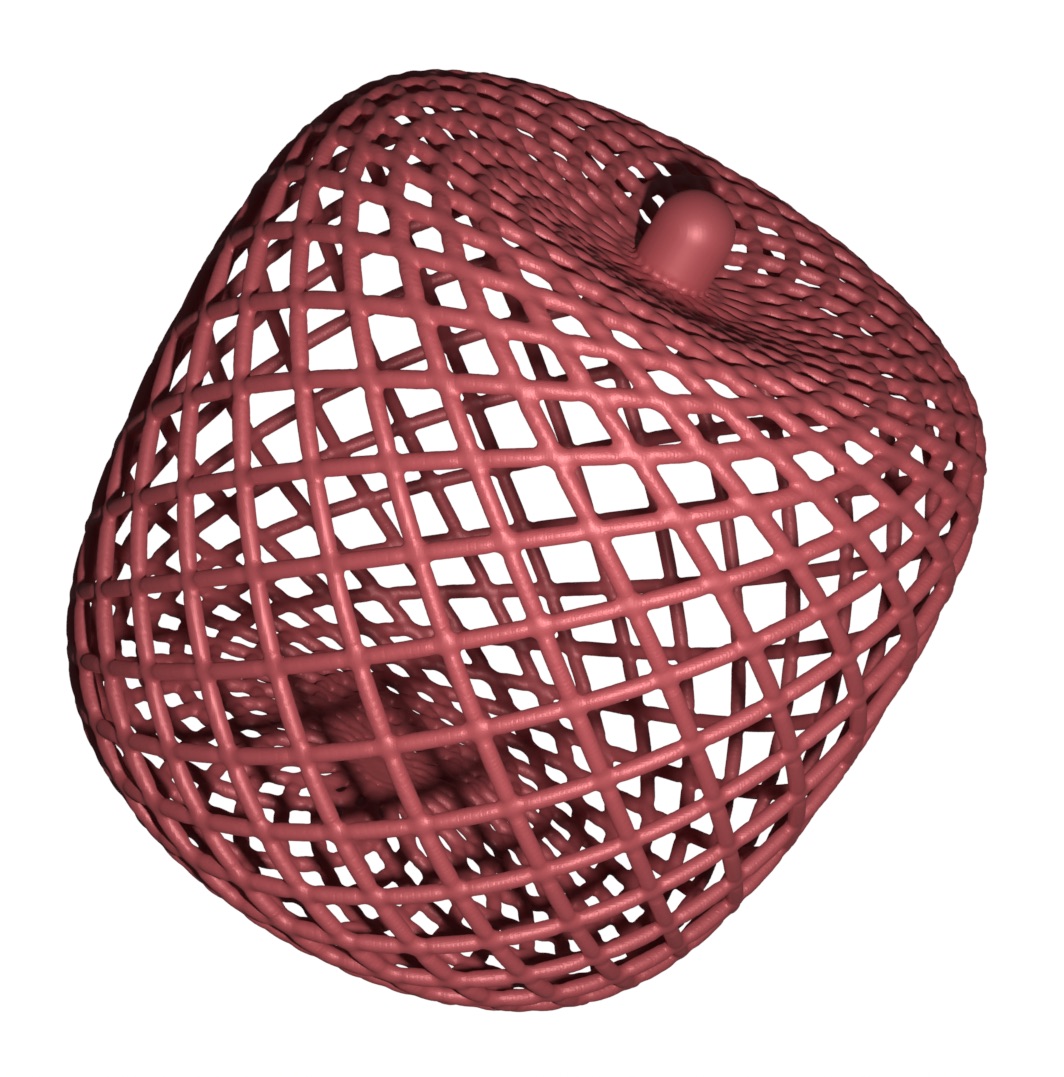}\includegraphics[width=4.25cm]{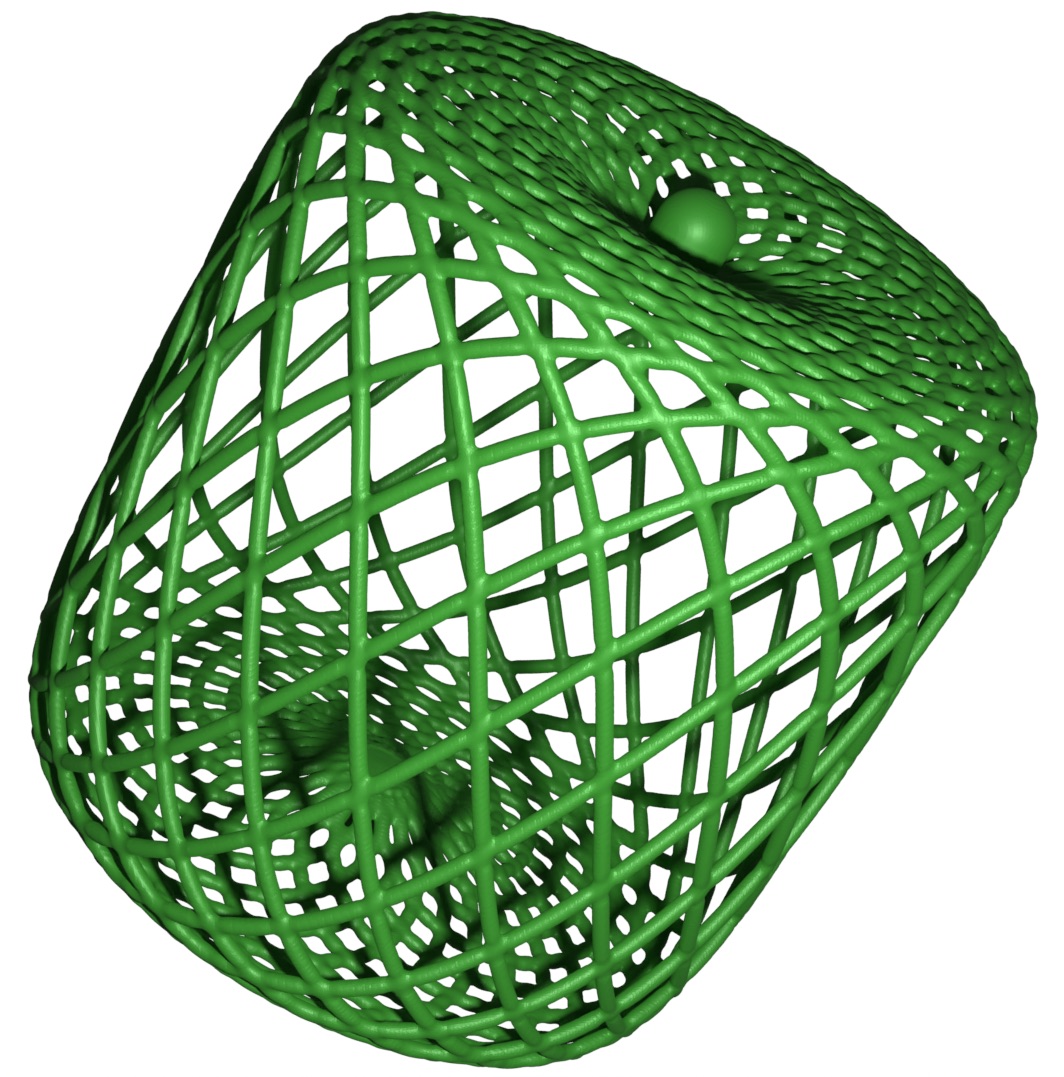}}
\caption{Parametrized {\web} device models using different design parameters. The parameter $l$ (compare blue vs. yellow) steers the depth of the bundling points at the device's poles. The parameter $w$ (compare blue vs. red) steers the width of the bundling point sink and $b$ (compare blue vs. green) probably takes the most important influence on the devices's shape transitioning between a more spherical one (small $b$) and a more cylindrical one (large $b$). All devices have radiopaque markers attached at the ``poles'' in a post-processing step.\label{fig:WebModel}}
\end{center}
\end{figure}

The creation of {\web} device shapes relies on the boolean union of a rotationally symmetric arranged collection of individual parameterized threads, which at the ``north-'' and ``south-pole'' of the device can even be equipped with the half-``pill''-shaped radiopaque markers,
i.e. the ``tips'' of the device.
The main manufacturer shape parameters used for each individual thread as well as for the thread-collection and hence complete web-device are
the threads radius $r_t~\left[\si{\milli\meter}\right]$,
the (maximal) radius~$r_{d,\max}~\left[\si{\milli\meter}\right]$ and height~$h_{d,\max}~\left[\si{\milli\meter}\right]$
of the {\web} device base-cylinder,
the number of winding rounds per thread $n_w~\left[1\right]$,
and the number of thread-pairs $n_{t}~\left[1\right]$, where each of them consists out of a clockwise and a counter-clockwise tortuous thread.
These parameters are now used to spin the threads, parameterized via $\theta\in[0,1]$, across the surface of a cylindrical shape with varying radius $r_d(\theta)$ which becomes zero at the ``poles'' of the shape in order to close it at the top and bottom.
Assuming the radial direction of the cylinder always resides within the $x$-$y$-plane,
we equip the threads with a parameterized $z$-coordinate~$z_d(\theta)$, which will allow also for re-entrant parts such as at the bundling points of the radiopaque markers.
Taking inspiration from the (level set) parametrization of a general stent graft in~\cite{sogn2014stabilized},
the parametrization of the $i$-th, $i=1,2,...,n_t$, thread pair $\vec{\gamma}_i$ with one thread in clockwise ($\vec{\gamma}_{i,\textup{c}}$) and one in counter-clockwise direction ($\vec{\gamma}_{i,\textup{cc}}$) -- distinguished by the different signs $\pm$ within the trigonometric function -- reads:
\begin{align*}
    \vec{\gamma}_{i,\textup{c}/\textup{cc}}:[0,1]\rightarrow\mathds{R}^3, \quad \theta\mapsto \vec{\gamma}_{i,\textup{c}/\textup{cc}}(\theta)=\left(\begin{array}{c}
         r_{d,\max}\cdot r_d(\theta)\cdot \cos(\pm 2\pi n_w\theta + \varphi_i) \\
         r_{d,\max}\cdot r_d(\theta)\cdot \sin(\pm 2\pi n_w\theta + \varphi_i) \\
         h_{d,\max}\cdot z_d(\theta)
    \end{array}\right),\qquad \varphi_i = \frac{2\pi(i-1)}{n_t}.
\end{align*}
For our concrete shapes, we are using the following radius and $z$-coordinate/height functions (see Figure~\ref{fig:WebModel} for a general overview over the function $z_d(\theta)$ as well as the influence of its parameters):
\begin{align*}
    r_d(\theta) &= 1-(2\theta -1)^{2}\\
    z_d(\theta) &= \left[\frac{1}{2}\tanh\left(b\cdot\left(\theta-\frac{1}{2}\right)\right)+\frac{1}{2}\right] + \left[l\cdot \textup{sign}\left(\theta-\frac{1}{2}\right)\cdot \large\big(m(\theta; w)-1\big)\right]
\end{align*}
Regarding $z_d(\cdot)$ the expression in the \textit{first} bracket determines the overall shape of the device while the expression in the \textit{second} bracket is responsible for the (inward pointing) tip parts of the device at its ``poles''. Here we have chosen the following further \textit{mathematical} - and purely {\web} device related - shape parameters, in contrast to the aforementioned \textit{manufacturer's} parameters which will be reused in the upcoming flow diverter model.
\begin{itemize}
    \item $b\in\mathds{N}$ steering the ``boxiness'' of the {\web}-device
    \item $l\in[0,1/2]$ steering the \textit{length} of the (inward pointing) tips of the device
    \item $w\in[0,1/2]$ steering the \textit{width} of the (inward pointing) tips of the device
\end{itemize}
Finally, the function~$m(\cdot; w)$, that appeared as part of the design function $z_d(\cdot)$, is a symmetric standard mollifier bridge function on $[0,1]$ with ascending flank from $0$ to $w$ and descending one from $1-w$ to $1$, which pronounces the region around the ``poles'' of the device for the aforementioned re-entrant tip modification.
It reads:
\begin{equation*}
    m(\cdot; w):[0,1]\rightarrow[0,1],\quad \theta\mapsto m(\theta; w)=\begin{cases}
        0, & \theta = 0\\
        \exp\left(1-\left(1-\left(\frac{x-w}{w}\right)^2\right)^{-1}\right), & \theta\in \left(0,w\right]\\
        1, & \theta\in (w,1-w)\\
        \exp\left(1-\left(1-\left(\frac{x-(1-w)}{w}\right)^2\right)^{-1}\right), & \theta\in[1-w,1)\\
        0, & \theta = 1
    \end{cases}
\end{equation*}
The individual parameterized curves~$\vec{\gamma}_{i,\textup{c,cc}}$ are then inflated to becomes tubes of radius~$r_t$,
then combined to the whole {\web} device and surface meshed using \textit{gmsh}\footnote{gmsh software webpage: \url{https://gmsh.info/}} \cite{gmsh} and \textit{MeshLab}\footnote{MeshLab software webpage: \url{https://www.meshlab.net/}} \cite{LocalChapterEvents:ItalChap:ItalianChapConf2008:129-136}. In Figure~\ref{fig:WebModel}, the influence of a few if these shape parameters, both manufacturers as well as mathematical ones, on the {\web} devices appearance is shown.

\subsubsection{Placement of {\web} devices}
\begin{figure}
\begin{center}
\subfigure[\textit{Un}deformed ``design-shape'' {\web} device inserted along predefined axis still leaving large inflow zones due to insufficient wall contact. ]{\label{fig:WebInsertionGreenUnDeforemd}\includegraphics[width=7cm]{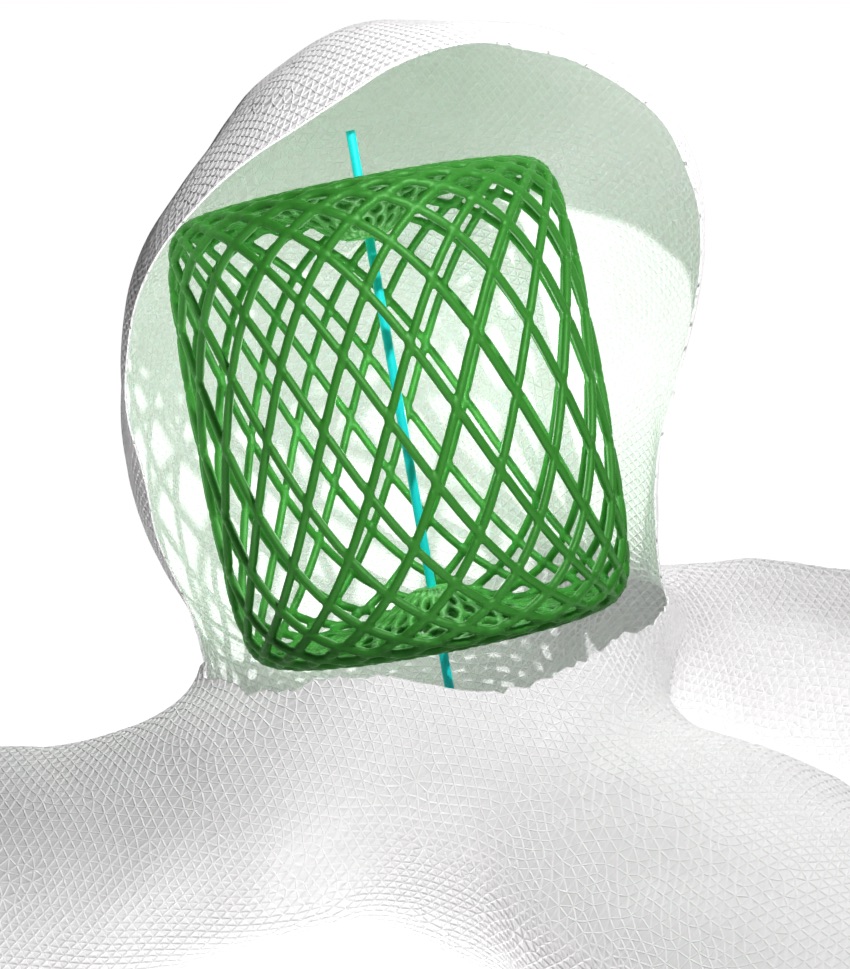}}
\hspace{1cm}
\subfigure[Deformed {\web} device aligning with the aneurysm surface reaching much higher occlusion.]{\label{fig:WebInsertionDeforemd}\includegraphics[width=7cm]{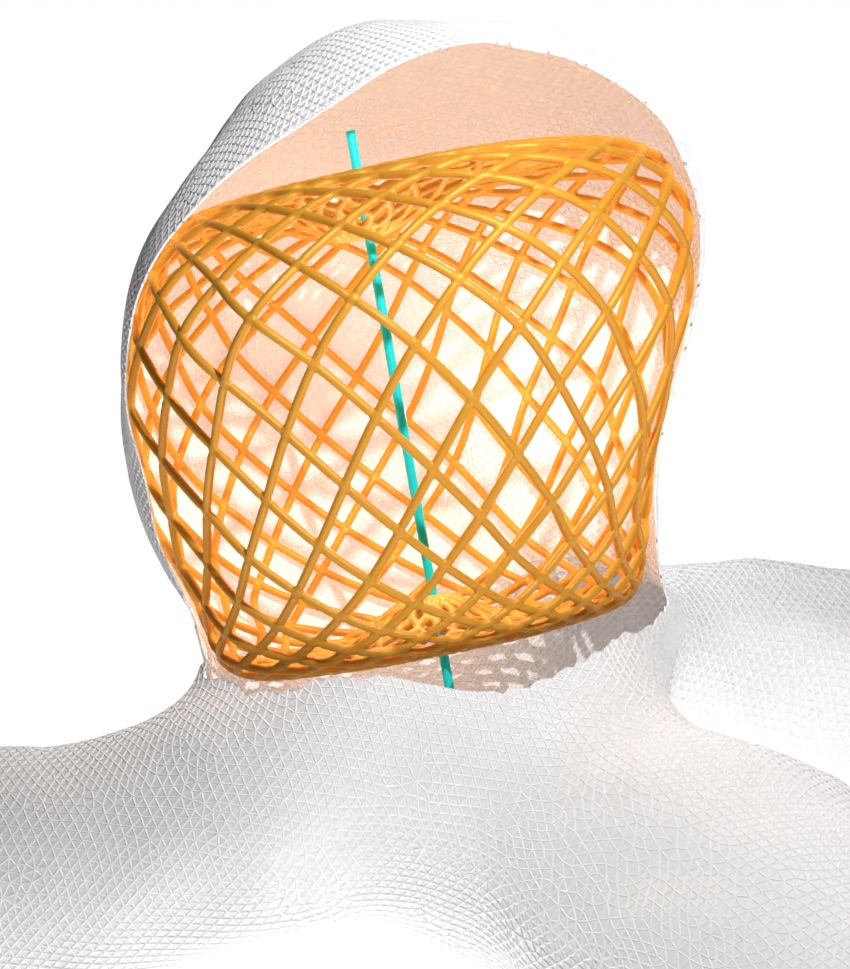}}\\
\subfigure[View from inside the vessel onto the \textit{un}deformed device's bottom surface]{\includegraphics[width=7cm]{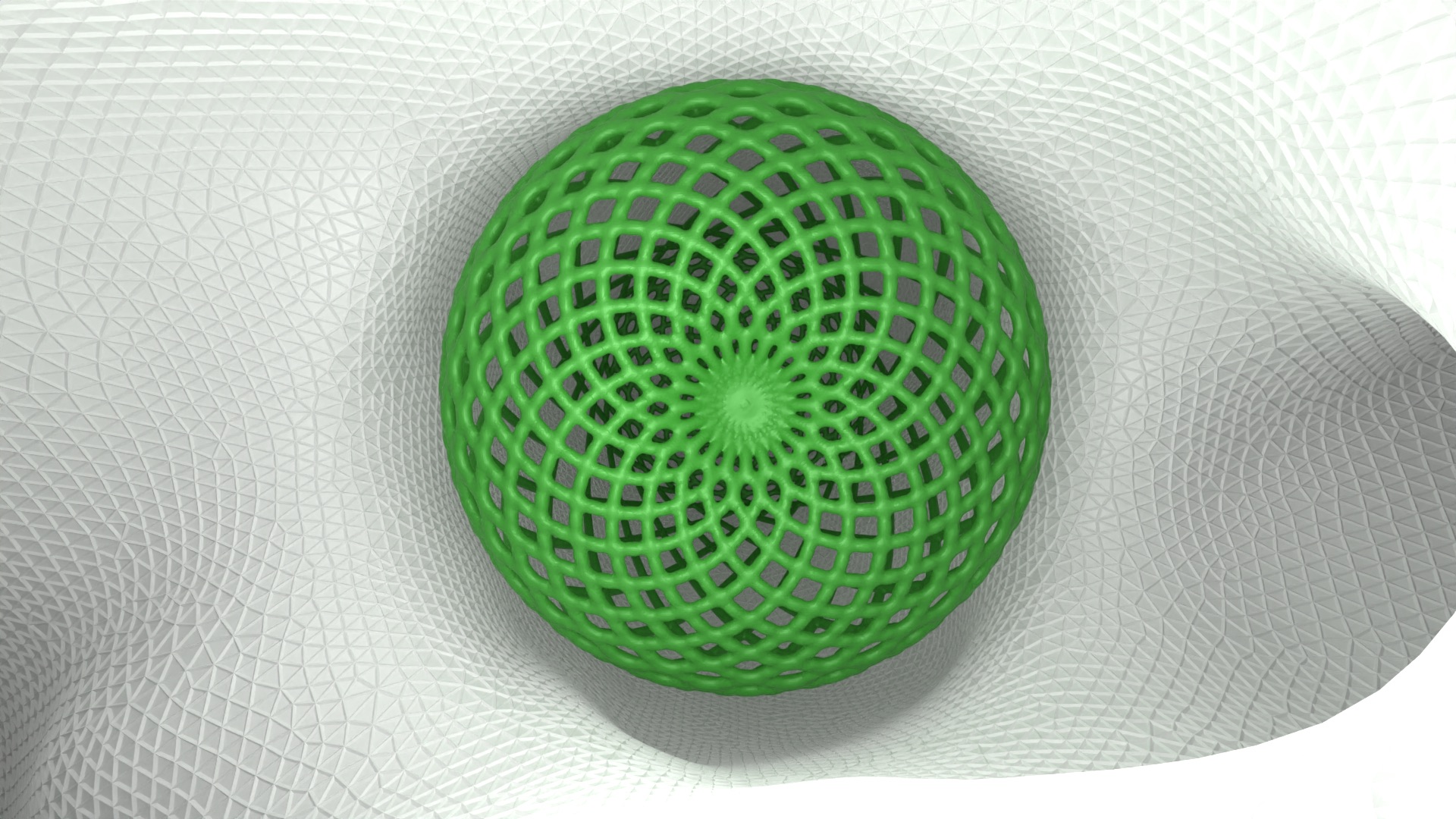}}
\hspace{1cm}
\subfigure[View from inside the vessel onto the expanded device's bottom surface]{\includegraphics[width=7cm]{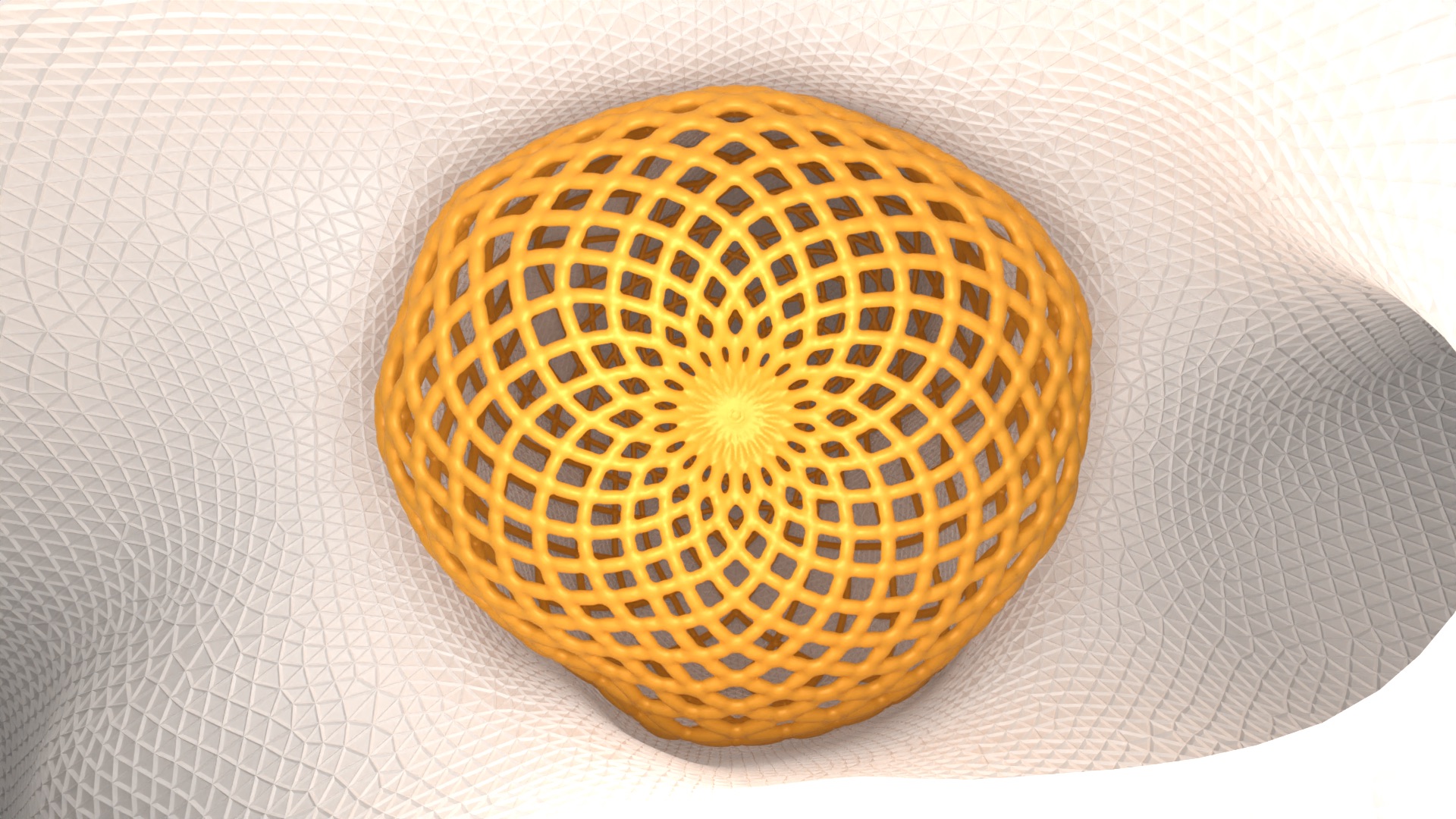}}
\caption{\emph{Left column with green devices:} Insertion and deformation process of {\web} devices along a predefined insertion direction axis (cyan).
\emph{Right column with yellow devices:} Improved occlusion due to the device deformation (viewed from within the vessel).}
\label{fig:WebInsertion}
\end{center}
\end{figure}

The shape design of {\web} devices so far always dealt with their unobstructed, contact less shape that one would expect when protruding the device into void space or centrally into a perfectly cylindrical aneurysm of exactly matching radius $r_{d,\max}-r_t$.
Such an idealistic perfect aneurysm would be filled exactly up to its surface by the {\web} device and hence could be fully occluded from the blood-stream, assuming the device is inserted in a way that its bottom, for large values of $b$ (nearly) planar surface aligns with the ostium of the aneurysm. Unfortunately realistic aneurysms are seldomly shaped that perfectly but usually exhibit some ``organic'', saccular surface shape with potentially further outpouchings, or in the extreme even smaller aneurysms on top of an aneurysm. Filling such general shaped aneurysm with an \textit{un}deformed {\web} device by simple e.g. central placement would always result in imperfect or even non-sufficient occlusion, even risking that the device could ``fall'' out of the aneurysm into the adjacent vessel posing an even additional medical risk (see Figure~\ref{fig:WebInsertionGreenUnDeforemd}).

It is hence our next step to \textit{deform} a designed {\web} device in  such a way that it aligns as good as possible with a given aneurysm shape, resembling the realistic, clinical situation, where after measuring the concrete patient specific aneurysm the surgeon chooses a {\web} device that is slightly larger in diameter than the aneurysm such that during insertion it can gently press itself onto the aneurysms surface and hence attach to that. Realistically this deformation process involves complex mechanics on both sides, i.e. the {\web} device as well as the dome of the aneurysm, which might deform, even though not stretch, as well. As before, for the purpose of this work, we will not involve any mechanical model to resolve these processes, but instead employ a purely geometric model, that deforms a given {\web} device, inserted with a given axis into the aneurysm, simply by expanding it in radial direction until it touches the aneurysm walls. By allowing angular as well as axial height dependent variation of that expansion factor, we can bring all the device's threads into persistent contact with the aneurysm's wall and again fully occlude the aneurysm (see Figure~\ref{fig:WebInsertionDeforemd}). The depth of insertion, that is responsible for the amount of empty volume \textit{above} the {\web} device towards the \textit{upper} end of the insertion axis as well as the amount of device that protrudes into the vessel after insertion at the \textit{lower} end of the axis, can be steered as well simply by placing the undeformed aneurysm prior to the expansion deeper or more superficial inside the aneurysm, i.e. let it slide across the insertion axis until the desired depth is reached. In the example in Figure~\ref{fig:WebInsertion}, we have chosen the depth in a way that at the parent vessel end the device closes the aneurysms without any larger protrusion hence connecting the vessel branching to the left and right from the aneurysm relatively straight.

\subsection{Stents and flow diverters}\label{subsec:Stents}

In contrast to coils or {\web} devices, stents and flow diverters are usually placed across the aneurysm's opening in the parent vessel.
Flow diverters aim to redirect the blood stream across the aneurysm's opening, reducing the inflow through the narrow windings of the flow diverter drastically. Stents come with a less denser mesh and are intended to intra-aneurysmal devices such as coils to fall out from the aneurysm sac into the vessel, hence stabilizing its position and effect.

For the geometric modeling of such stents, we follow very similar steps as for the {\web} device.
Again, the complete device comprises of several parametrized threads being spun across a cylindrical shape.
In its straight reference shape,
the parametrization of its threads is actually even simpler than in the {\web} device case,
since the radial shape function~$r_d(\theta)$ can just be chosen as constant~$r_d(\theta)\equiv 1$ and~$z_d(\theta)=\theta$ as linear,
giving the stent the shape of a straight cylinder.
Hence, with the same manufacturer shape parameters as in the {\web} device case,
the $i$-th, $i=1,2,\dots,n_t$, individual thread-pair's parametrization -- again one clockwise, one counter-clockwise -- for the stent reads:
\begin{align}\label{eq:StentRefparametrization}
    \vec{\gamma}_{i,\textup{c}/\textup{cc}}:[0,1]\rightarrow\mathds{R}^3, \quad \theta\mapsto \vec{\gamma}_{i,\textup{c}/\textup{cc}}(\theta)=\left(\begin{array}{c}
         \cos(\pm 2\pi n_w\theta + \varphi_i) \\
         \sin(\pm 2\pi n_w\theta + \varphi_i) \\
         h_{d,\max}\cdot \theta
    \end{array}\right),\qquad \varphi_i = \frac{2\pi(i-1)}{n_t}
\end{align}
For the maximal radius in the $x$-$y$-plane, we have deliberately chosen~$r_{d,\max}=1$ to actually obtain a ``normalized'' straight reference shape of the stent.
Once we start to insert the stent into a vessel, it not only needs to be deformed to follow the vessel's contours, 
but also can be scaled in radial direction to fit to the gauge of the vessel, either its minimal gauge (to again have a device of constant radius) or even adaptively to the local gauge along the vessel's arc length.

\begin{figure}
\begin{minipage}{0.15\textwidth}
    \subfigure[Three thread-pairs in straight reference configuration,i.e. the graphs of the curves $\Vec{\gamma}_{i,\textup{c/cc}}, i=1,2,3$, each with $n_w=4$ windings.]{\label{fig:StentModelReferenceState}\includegraphics[height=6cm]{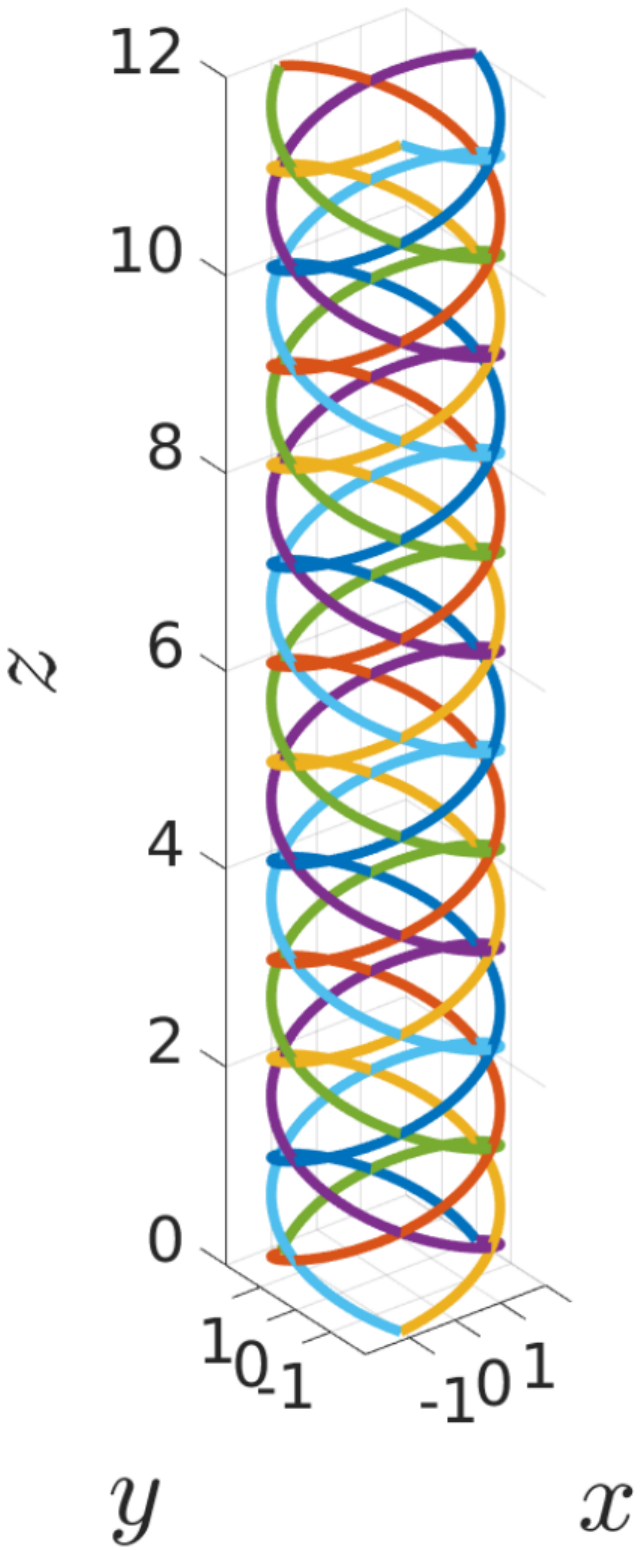}}
\end{minipage}\hspace{5mm}\begin{minipage}{0.5\textwidth}
\subfigure[Centerline of the blood-vessel with color-coded arclength. In blue the parallel-transported vector-field $\Vec{u}$, in orange the vector-field $\Vec{v}$.]{\label{fig:StentModelCenterline}\includegraphics[height=4.25cm]{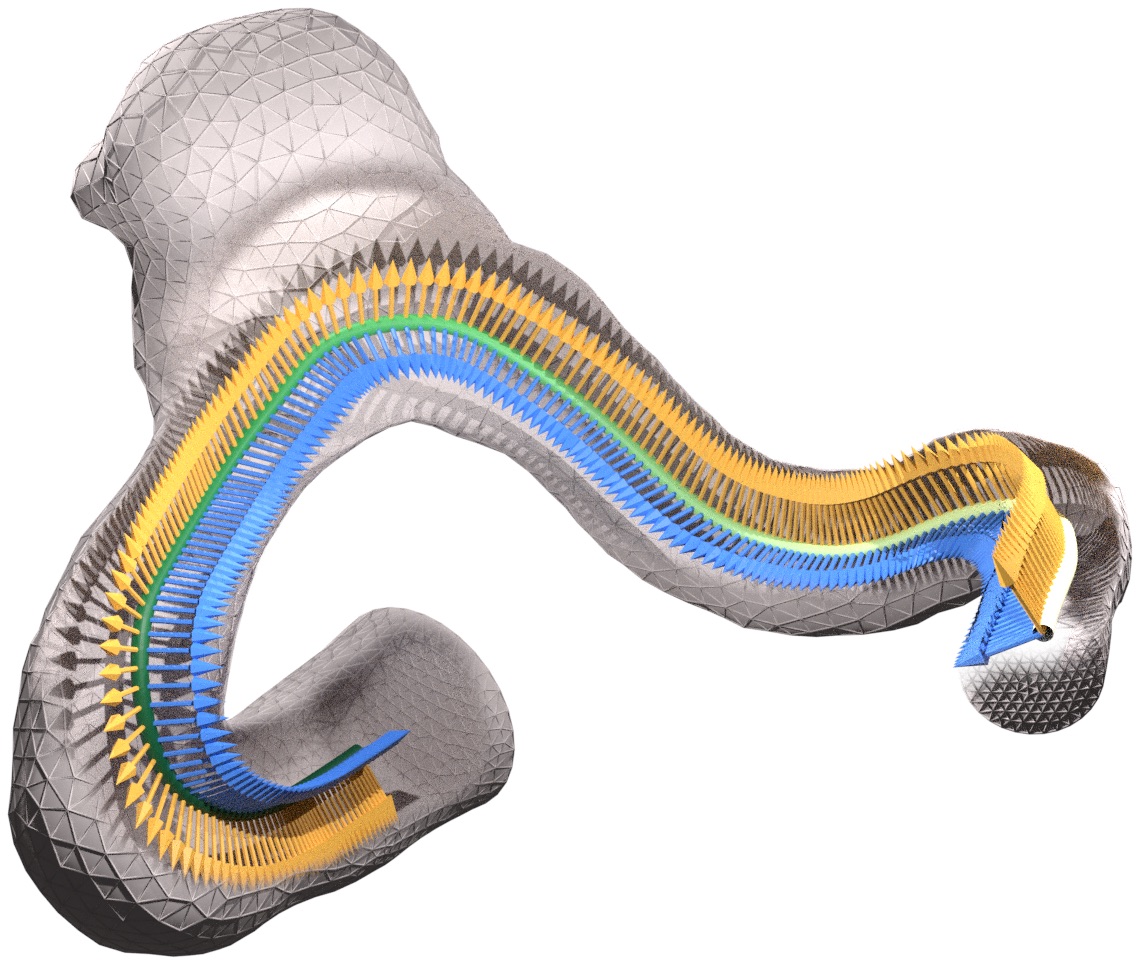}\includegraphics[height=5cm]{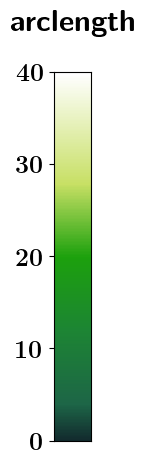}}\hspace*{5mm}\subfigure[The truncated centerline portion in front of the aneurysm with color coded vessel-radius information together with the physical device configuration, i.e. the graphs of the curves $\Vec{\Gamma}_{i,\textup{c/cc}}$. In transparent: The (curved) cylinder-surface on which the threads are spun, i.e. the image of the original reference cylinder.]{\label{fig:StentModelTruncatedCenterline}\includegraphics[height=4.25cm]{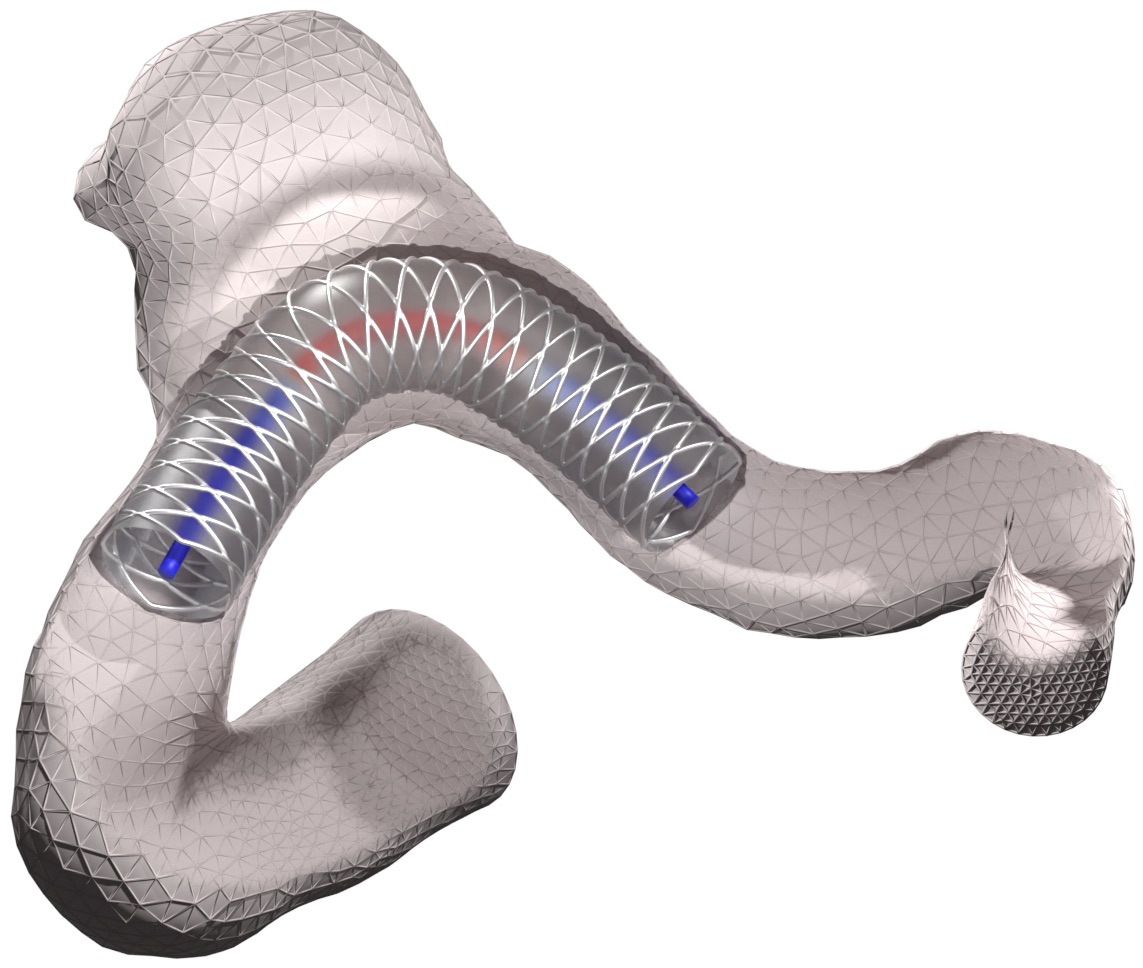}\includegraphics[height=5cm]{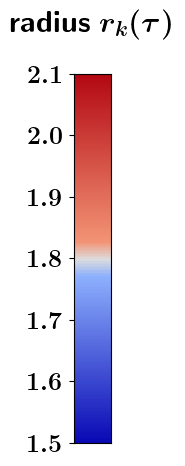}}\\
\subfigure[Different parameter variations of the same stent, all with $n_w=6$ windings; \emph{Left:} $n_t=4, r_t=0.05 \,\si{\milli\meter}$; \emph{Middle:} $n_t=8, r_t=0.025 \,\si{\milli\meter}$; \emph{Right:} same as in  left picture, but with local radius $r_d$ adapted to centerline radius information $r_k$ instead of a constant radius.]{\label{fig:StentModelRadii}\includegraphics[height=2.8cm]{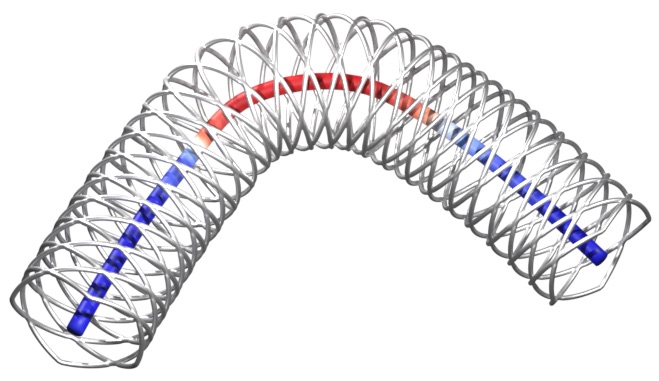}\includegraphics[height=2.8cm]{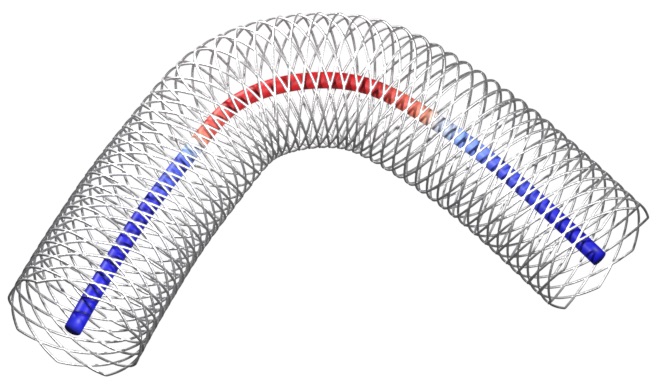}\includegraphics[height=2.8cm]{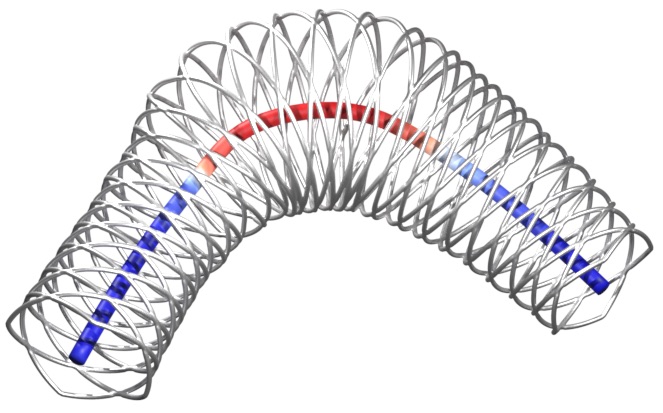}}
\end{minipage}
    \caption{Visualization of the stent creation process from reference configuration over the setup of an orthonormal, parallel-transported centerline frame, up to the finished stent-graft wrapped around the centerline with different manufacturer's parameters. The vessel geometry used in this case is \#AN166 from~\cite{IntrA}.}
\label{fig:StentModel}
\end{figure}

Having the parametrization of the reference stent's threads available, we now transform them such that the resulting device adapts to a given vessel geometry by curvilinear bending to take its place in front of the aneurysm.
This is achieved by transforming the reference thread's parametrizations to not spin around the straight vertical $z$-axis of the reference configuration anymore, but actually around the centerline of the vessel-geometry with a radial distance depending on the gauge of the vessel.
To this end, we proceed as follows (see Figure~\ref{fig:StentModel} for a graphical depiction of the upcoming process):
\begin{enumerate}
    \item \textbf{Reading centerline information:} Starting from a preprocessed vessel geometry, we approximate its centerline by linear interpolation between discrete curve points~$\Vec{x}_i, i=1,2,\dots,n$ computed by \textit{vmtk} and denote it by~$\vec{k}:[0,L]\rightarrow\mathds{R}^3,~\tau\mapsto \vec{k}(\tau)$.
    Furthermore by $r_k:[0,L]\rightarrow\mathds{R},~\tau\mapsto r_k(\tau)$, we denote the radius of the maximal vessel inscribed sphere around the point at parametric coordinate~$\tau$ on the centerline, hence its minimal distance from the vessel walls, also obtained by linear interpolation between the discrete curve points. By~$L$ we refer to the overall arc length of the centerline, see the two color-bared quantities in Figures~\ref{fig:StentModelCenterline} and~\ref{fig:StentModelTruncatedCenterline}.
    \item \textbf{Truncate the centerline:} We truncate the centerline to the user-chosen parametric interval $[\tau_{\textup{s}},\tau_{\textup{e}}]\subseteq[0,L]$ that covers the region of the vessel, where the stent is supposed to be placed, i.e. in front of the aneurysm. The length~$\tau_{\textup{e}}-\tau_{\textup{s}}$ is hence the arc length~$h_{d,\max}$ of the stent to be inserted.
    By abuse of notation from here on we will re-index the discrete curve points~$\Vec{x}_i$ that constitute the truncated portion of the centerline as well by $i=1,2,\dots,n$ and re-parameterize the interval~$[\tau_{\textup{s}},\tau_{\textup{e}}]$ to $[0,h_{d,\max}]$.
    Furthermore, the name $\vec{k}$ and the radius-information $r_k$ are carried over to the truncated curve as well, see Figure~\ref{fig:StentModelTruncatedCenterline}.
    
    \item \textbf{Smooth the centerline:} To denoise the centerline obtained from \textit{vmtk},
    a one-dimensional version of a Laplace-Taubin smoother~\cite{taubin1995signal} listed in Algorithm~\ref{alg1} is applied to smooth the centerline representation.
    We denote by $\Vec{x}_i^{(0)}, i=1,2,\dots, n$ the initial discrete curve points of the (truncated) centerline stemming from the \textit{vmtk} output.
    We perform $n_{\textup{it}}$ smoothing iterations with smoothing parameters~$\lambda$ and~$\mu$.

    \begin{algorithm}
        \caption{Laplace-Taubin smoother, also see~\cite{taubin1995signal}}\label{alg1}
        \begin{algorithmic}
        \For {$k=0,1,...,n_{\textup{it}}-1$}
        \State $\vec{\Delta}_i^{(k)}= \frac{1}{2}\left({\Vec{x}_{i+1}^{(k)}-2\vec{x}_i^{(k)}+\Vec{x}_{i+1}^{(k)}}\right),\quad \forall i=2,...,n-1$
        \If {$k$ \textup{~is even}}
        \State $\Vec{x}_i^{(k+1)}=\Vec{x}_i^{(k)} + \lambda \vec{\Delta}_i^{(k)},\quad \forall i=2,...,n-1$
        \Else {\textbf{if~} $k$ \textup{~is odd~}\textbf{then}}
        \State $\Vec{x}_i^{(k+1)}=\Vec{x}_i^{(k)} - \mu \vec{\Delta}_i^{(k)},\quad \forall i=2,...,n-1$
        \EndIf
        \EndFor
        \end{algorithmic}
    \end{algorithm}
    To be on the safe side, we usually choose $n_{\textup{it}}=512$, $\lambda=0.5$, and $\mu=0.25$

    By abuse of notation, we again keep the name $\vec{k}$ for the now smoothed centerline curve, $\Vec{x}_i$ for its discrete curve-points, $r_k$ for its radius information as well as the interval naming $[0,h_{d,\max}]$ for its truncation to the desired portion of the complete centerline,
    where the concrete polygonal segment's arc lengths might have changed due to the smoothing operation.

    \item \textbf{Compute centerline tangents:} For the polygonal curve $\vec{k}$, we can define a discrete tangent vector~$\Vec{t}_i$ at every interior point~$\Vec{x}_i, i=2,\dots, n-1$ as the length-aware average of the adjacent edge directions, evaluating to 
    \begin{equation*}
        \Vec{t}_i=\frac{\Vec{x}_{i+1}-\Vec{x}_{i-1}}{\|\Vec{x}_{i+1}-\Vec{x}_{i}\|_2+\|\Vec{x}_{i}-\Vec{x}_{i-1}\|_2},
    \end{equation*}
    while we use $\Vec{t}_1=\Vec{x_2}-\Vec{x}_1$ as well as $\Vec{t}_{n}=\Vec{x}_n-\vec{x}_{n-1}$ at the boundary points, then normalizing each of them to length one.

    \item \textbf{Create an orthonormal frame at the beginning of the centerline:} At the beginning of the truncated centerline, i.e. at $\Vec{x}_1$,
    we use $\Vec{t}_1$, choose a second unit-vector $\Vec{u}_1$, which is arbitrary, but orthogonal to $\Vec{t}_1$ (e.g. a random vector not parallel to $\Vec{t}_1$ then taking out its projection onto $\Vec{t}_1$ and normalizing) and form a third vector $\Vec{v}_1=\Vec{t}_1\times \Vec{u}_1$. The triad $(\Vec{t}_1,\Vec{u}_1,\Vec{v}_1)$ then forms an orthonormal frame associated with the curve-point $\Vec{x}_1$, i.e. so far at the beginning of the centerline only!

    \item \textbf{Parallel transport of initial frame:} We now parallel-transport the orthonormal frame $(\Vec{t}_1,\Vec{u}_1,\Vec{v}_1)$ along the discrete centerline curve, successively from point $\Vec{x}_i$ to point $\Vec{x}_{i+1},~i=1,2,\dots,n-1$. Each such step of the (space) parallel-transport is facilitated by determining the rotation that maps the current points tangent $\Vec{t}_i$ onto the next point's tangent $\Vec{t}_{i+1}$ and then apply its representing rotation matrix $\tens{R}$ to the whole frame $(\Vec{t}_i,\Vec{u}_i,\Vec{v}_i)$ to obtain the new frame $(\Vec{t}_{i+1},\Vec{u}_{i+1},\Vec{v}_{i+1})=(\tens{R}\Vec{t}_i,\tens{R}\Vec{u}_i,\tens{R}\Vec{v}_i)$. By means of remark~\ref{rem:RodriguesFormula}, this rotation matrix is well defined and can be explicitly computed using formula~\eqref{eq:Rotation} by choosing $\Vec{a}=\Vec{t}_i$ and $\Vec{b}=\Vec{t}_{i+1}$. See Figure~\ref{fig:StentModelCenterline}) for a graphical representation of the orthonormal frame vector fields~$\vec{u}_i$ and~$\vec{v}_i, i=1,2,\dots,n$.
    
    We employ specifically the described way of generating an orthonormal frame for the centerline, since this parallel transported frame, which is also refer to as Bishop frame in~\cite{jawed2018primer}, is twist-free. Hence, if we base the mapping from the reference configuration to the physical configuration on this transport, it will also generate a bent, but twist-free stent. We note that currently the frame is only defined at discrete points $\Vec{x}_i$. However, as with the centerline curve~$\Vec{k}$ itself, we interpret $(\Vec{t},\Vec{u},\Vec{v})$ as the continuous linear spline interpolation between the discrete frame-points $(\Vec{t}_i,\Vec{u}_i,\Vec{v}_i)$ again parameterized over the arc length interval $[0,h_{d,\max}]$ for a subsequent, continuous frame evaluation.

    \item \textbf{The stent-graft mapping:} We now combine the thread-wise reference configuration parametrization $\vec{\gamma}_{i,\textup{c/cc}}$ of the stent~\eqref{eq:StentRefparametrization} with the parallel frame to construct the following thread-wise parametrization of the physical configuration of the stent~$\vec{\Gamma}_{i,\textup{c/cc}}$, where we are using the short-hand notation~$\vec{\gamma}_{i,\textup{c/cc}} = (\gamma_x,\gamma_y,\gamma_z)^{\top}$ for the three spatial components of~$\vec{\gamma}_{i,\textup{c/cc}}$:
    \begin{align*}
        \vec{\Gamma}_{i,\textup{c/cc}}:[0,1]\rightarrow\mathds{R}^3,\quad \theta\mapsto \Vec{\Gamma}_{i,\textup{c/cc}}(\theta) = (\vec{k}\,\circ\,\gamma_z)(\theta) + (r_d\,\circ\,\gamma_z)(\theta)\cdot\left[\gamma_x(\theta)\cdot(\vec{u}\,\circ\, \gamma_z)(\theta) + \gamma_y(\theta)\cdot (\Vec{v}\,\circ\, \gamma_z)(\theta)\right]
    \end{align*}
    Within this formula, the first additive term drives along the centerline to the position corresponding to arc length~$\gamma_z(\theta)$.
    The two terms in brackets then map the radial components of the individual thread curves from their $x$/$y$-axis reference values~$\gamma_x(\theta)$ and~$\gamma_y(\theta)$ into the bent, local $\Vec{u}/\Vec{v}$-coordinate frame.
    In front of the bracket, there is again a term denoted by~$r_d$ which allows to radially scale the stent. Here, $r_d$ could either be chosen as $r_d(\tau)\equiv\min_{\tau\in [0,h_{d,\max}]}r_k(\tau)$, i.e. to represent the minimal vessel gauge, or alternatively $r_d(\tau)=r_k(\tau),~\forall\tau\in[0,h_{d,\max}]$ for the local vessel gauge.

    \item \textbf{Thread inflation and meshing:} Similar to the {\web} device, once all threads are placed as desired, they are inflated to become tubes of radius $r_t$ (see manufacturer shape parameters in Section~\ref{subsec:WebDevices}), then combined and surface meshed to yield the final stent-graft object.
\end{enumerate}

\section{Numerical experiments} \label{sec:NumericalSimulations}

Now, we will employ the numerical models introduced in Section~\ref{sec:ContinuumMechanicsAndNumericsOfFlowInAneurysms}
to showcase some concrete use cases.
In a first example, we will study the modeling approach for a {\FD} by a homogenized deformable porous medium from Section~\ref{subsec:PoroElasticFlow}
on a parametric aneurysm geometry and discuss the influence of modeling parameters.
In a second experiment, we will assess the influence of a medical treatment device on quantities of interest
by a fully resolved simulation of the fluid flow field, once done with and once without the device being inserted.

\subsection{Proof of concept: poro-elastic modeling of an aneurysm treated with a flow diverter} \label{subsec:poroexample}

We will now outline a numerical example to discuss the modeling approach of elastic porous media
to represent the fine wire structures of the medical devices in a homogenized fashion.
Exemplarily, we will therefore consider a {\FD}, where the wire structure is represented as a homogenized porous medium.

\begin{figure}
\begin{minipage}{0.3\linewidth}
 \includegraphics[scale=0.75]{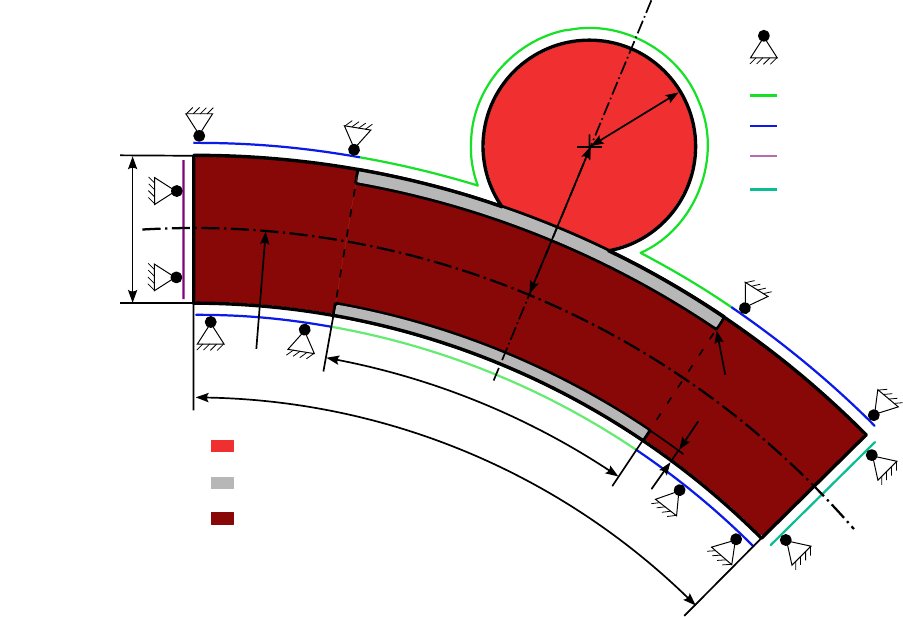}
\scalebox{0.750}{
\begin{tikzpicture}[overlay]

\node[text width=1cm] at (10.9 ,9.3) {$r_{an}$};
\node[text width=1cm] at (9.4 ,7.5) {$h_{an}$};

\node[text width=1cm] at (8,4.6) {$l_{FD}$};
\node[text width=1cm] at (12.85,4.5) {$r_{FD}$};
\node[text width=1cm] at (12.3,3.5) {$d_{FD}$};

\node[text width=1cm] at (2.1,7.0) {$d_{A}$};
\node[text width=1cm] at (4.8,4.7) {$d_{a}$};
\node[text width=1cm] at (8.3,3.6) {$l_{A}$};

\node[text width=1cm] at (4.7,3.5) {$\Omega_{CA}$};
\node[text width=1cm] at (4.7,2.85) {$\Omega_{FD}$};
\node[text width=1cm] at (4.7,2.2) {$\Omega_{A}$};

\node[text width=4cm] at (15.5,10.3) {fixed displacement};
\node[text width=5cm] at (16,9.45) {impermeable boundary for fluid};
\node[text width=4cm] at (15.5,8.9) {no-slip condition};
\node[text width=4cm] at (15.5,8.4) {arterial inflow};
\node[text width=4cm] at (15.5,7.85) {arterial outflow};

\node[text width=1cm] at (5.5,7.5) {};

\node[text width=4cm] at (19,5.5) {$r_{an} = 2.25 \si{\milli\meter}$}; 
\node[text width=4cm] at (19,5.0) {$h_{an} = 3.0 \si{\milli\meter}$};
\node[text width=4cm] at (19,4.5) {$l_{FD} = 22.5 \si{\degree}$};
\node[text width=4cm] at (19,4.0) {$r_{FD} = 0.105 \si{\milli\meter}$};
\node[text width=4cm] at (19,3.5) {$d_{FD} = 0.15 \si{\milli\meter}$};
\node[text width=4cm] at (19,3.0) {$d_{A} = 3 \si{\milli\meter}$};
\node[text width=4cm] at (19,2.5) {$d_{a} = 12 \si{\milli\meter}$};
\node[text width=4cm] at (19,2.0) {$l_{A} = 45^\circ$};

\end{tikzpicture} }
\end {minipage}
    \caption{Two-dimensional sketch of the parametric geometry for the three different domains with dimensions and boundary conditions. The colored lines represent the different boundary conditions of the fluid phase.}
    \label{fig:ToyFD}
\end{figure}

For now,
we will restrict ourselves to a simplified artificial geometry. 
The suggested parametric geometry of aneurysm is an idealization of a saccular aneurysm,
which is composed of a torus and a sphere segment, respectively.
Figure~\ref{fig:ToyFD} sketches the geometry with the placed device and all associated geometric dimensions.
The sphere is positioned at a distance of~$h_{an}$ between the centerline and the center of the sphere.
For the in-silico treatment of this saccular aneurysm, a {\FD} is inserted into the artery.
The final configuration of the {\FD} is assumed to be ideal,
such that it attaches perfectly to the aneurysm and surrounding artery.
Additionally, the position of the {\FD} is symmetrical with respect to the rest of the aneurysm.
In addition, a flow extension of length~$3d_A$ is added at the outflow cross section
to stabilize the flow and help with backflow stabilization~\cite{Bertoglio2017}.

The proposed modeling approach uses the equations from Section~\ref{subsubsec:poroelasteq} for the three types of domains,
namely the artery~$\Omega_A$, the {\FD}~$\Omega_{FD}$ and the {\CA}~$\Omega_{CA}$.
These domains may be interpreted as a single computing domain with varying material parameters,
so that no additional coupling conditions are necessary to relate the device with the surroundings.
The underlying idea is to use high values for permeability and porosity in areas with expectedly high fluid flow
(such as the artery or {\CA} cavity),
whereas the material parameters within the medical device are chosen closer to a less permeable porous medium.

The constitutive behavior of the porous skeleton governed by~$\macrostrainskel$ is assumed to follow a hyperelastic isotropic Neo-Hookean law~\cite{Holzapfel2000a}
with Young's modulus~$E_{NH}$ and Poisson's ratio~$\nu_{NH}$ described by the strain energy function
\begin{equation*} 
    \Psi_{NH}=
 \frac {E_{NH}}{4(1+\nu_{NH})}(\text{tr}~\boldsymbol{C}-3)+\frac {E_{NH} \nu_{NH}}{4(1+\nu_{NH})(1-2\nu_{NH})} (J^{-2\frac {\nu_{NH}}{1-2\nu_{NH}}}-1)
\end{equation*}
with~$\rightcauchgreendeftensor$ and~$J$ denoting the right Cauchy-Green tensor and the determinant of the deformation gradient, respectively.

The mesh is created by the preprocessing pipeline outlined in \ref{subsec:GeometryPreprocessing}
and consists of $59769$ nodes in $59360$ quadrilateral linear finite elements
with approximately $358614$ degrees of freedom in total. For the One-Step-Theta time integration scheme, a constant time step size of $\Delta t=10^{-3} \si{\second}$ with a $\Theta=0.6$ is used.
The boundary conditions of the mixture are assigned with respect to the different phases,
which can be identified from Figure~\ref{fig:ToyFD}. 
For the solid phase,
we set the displacements in~$\Omega_A$ to zero, but leave $\Omega_{FD}$ and~$\Omega_{CA}$ free to deform.
For the fluid phase, we prescribe a pulsatile velocity profile at the inflow cross section
as outlined in Section~\ref{subsubsec:FreeFlowAnalysical}.
At the outflow cross section, a time-dependent pressure profile is imposed via a Neumann boundary condition.
To prevent potential instabilities due to backflow through the Neumann boundary,
we employ an additional backflow stabilization from~\cite{Bertoglio2017}.
At the interfaces to the arterial walls, we impose a no-slip condition on the velocity field.
Along the outer boundary of~$\Omega_{FD}$ and~$\Omega_{CA}$,
an impermeability condition is enforced along~$\boundaryconstraint$ via a Lagrange multiplier field~$\pmb{\lambda}$
to ensure that the velocities of both phases coincide at the boundary.
This condition can be stated as
\begin{equation*}
    \int \klr{\vf -\vs} \virtual \pmb{\lambda}\,\textup{d}\boundaryconstraint= 0,
\end{equation*}
resulting in a vanishing fluid flow across the boundary.
We finally prescribe the porosity~$\porosity$ everywhere in the domain:
For ~$\Omega_A$ and~$\Omega_{FD}$, the porosity is fixed at $0.999$ (as close as numerically possible to $1$),
whereas~$\Omega_{FD}$ adopts a constant porosity of $0.5$.

We will vary the permeability tensor~$\permea$ in the spatial configuration
using the Carman-Kozeny relation~\cite{Carman1997}
\begin{align}
  \permea = l_0^2\frac{\detdefgrad\porosity^3}{\left(1-\detdefgrad\porosity\right)^2} \Identity,  \qquad \permea= \pushforward{\permeainitial}, \label{eq:permea-carman}
\end{align}
where~$\Identity$ denotes the identity tensor,
$\permea$ is governed by a reference length with a geometric constant~$l_0$, that is related to the shape of the skeleton,
and $\permeainitial$ refers to the permeability tensor in the material configuration.
While the permeability is a tensor quantity, which may also account for anisotropic effects,
we restrict ourselves to an isotropic permeability here.
Aiming at a homogenized representation of a deforming device,
the idea is that the permeability acts as a resistance against the convective velocity
and therefore decreases on the one hand the velocity within $\Omega_{FD}$
and on the other hand increases the coupling stress within the solid balance of the mixture.

We will now assess the capability of the {\PM} modeling approach to represent a flow diverter placed across the opening of a {\CA}.
Therefore, we compare three instances of a flow diverter modeled with different values for the reference length of the permeability from~\eqref{eq:permea-carman} 
to a CFD reference solution of an untreated {\CA}.
In particular, we study cases for~$l_0^2\in\{1.0\cdot 10^{-5}, 6.0\cdot 10^{-6}, 5.0\cdot 10^{-7}\} \si{\milli\square\meter}$.
Further parameters are summarized in Table~\ref{tab:paraporosimu}.

\begin{table}
    \centering
    \caption{Selected parameters for the different simulation domains }
    \label{tab:paraporosimu}
    \subfigure[ Parameters for artery domain~$\Omega_A$ and aneurysm domain~$\Omega_{CA}$]{   
            \begin{tabular}{|c|c|c|} 
                \hline
                Parameter & Value & Unit \\ \hline 
                $\visfluid$ & $4$ & \si{\pascal \second }\\
                $\rhof$& $1.0\cdot 10^{-6}$ & \si{\kilogram \per \milli\cubic\meter}\\  
                \hline
                $\rhos$& $0.94\cdot 10^{-6}$ & \si{\kilogram \per \milli\cubic\meter}\\
                $E_{NH}$ & $1$ & \si{\newton \per \milli\square\meter}\\
                $\nu_{NH}$& $0.3$ & \\
                $\bulkmod$ & $1.0$ & \si{\newton \per \milli\square\meter} \\
                $\poromatpenalty$& $0.1$ & \\ \hline
                $\porosity$& $\approx 1$& \\
                $l_0^2$ & $1.0\cdot 10^{20}$ & \si{\milli\square\meter}\\
                \hline 
            \end{tabular}
            
        }
        \hspace{2cm}
    \subfigure[Parameters for {\FD} domain~$\Omega_{FD}$]{    
        \begin{tabular}{|c|c|c|} 
            \hline
            Parameter & Value & Unit \\ \hline 
            $\visfluid$ & $4$ & \si{\pascal \second }\\
            $\rhof$& $1\cdot 10^{-6}$ & \si{\kilogram \per \milli\cubic\meter}\\  
            \hline
            $\rhos$& $4.5\cdot 10^{-5}$ & \si{\kilogram \per \milli\cubic\meter}\\
            $E_{NH}$ & $2000$ & \si{\newton \per \milli\square\meter}\\
            $\nu_{NH}$& $0.3$ & \\
            $\bulkmod$ & $1100$ & \si{\newton \per \milli\square\meter} \\
            $\poromatpenalty$& $0.1$ & \\ \hline
            $\porosity$& $0.5$& \\
            $l_0^2$ & variable & \si{\milli\square\meter}\\
            \hline
        \end{tabular}
        }
\end{table}

\begin{figure}
\begin{center}
    \includegraphics[scale=0.4]{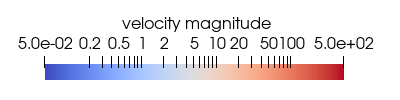}\\
    \subfigure[Reference solution of a CFD simulation ]
        {\label{fig:FDCFDReference}\includegraphics[width=7cm,trim=5cm 15cm 5cm 0cm,clip]{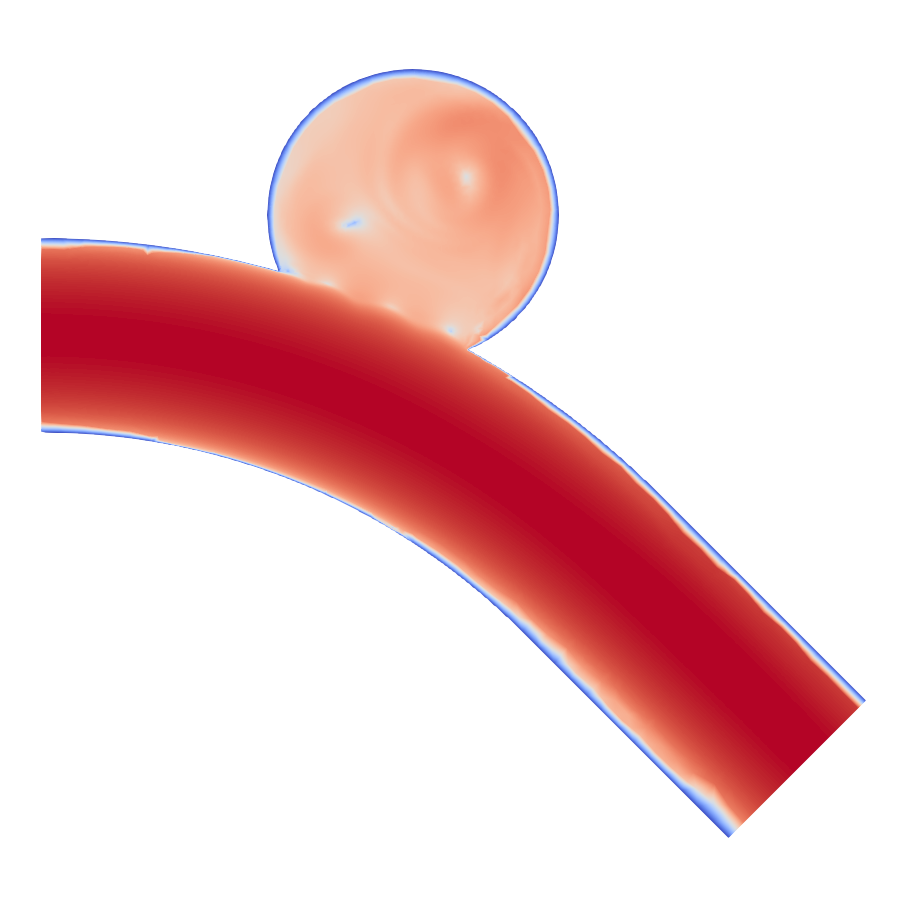}} \hspace*{1cm}
    \subfigure[Velocity solution for $l_0^2=1.0\cdot 10^{-5} \si{\milli\square\meter}$]
        {\label{fig:FDCFDk1e-5}\includegraphics[width=7cm,trim=5cm 15cm 5cm 0cm,clip]{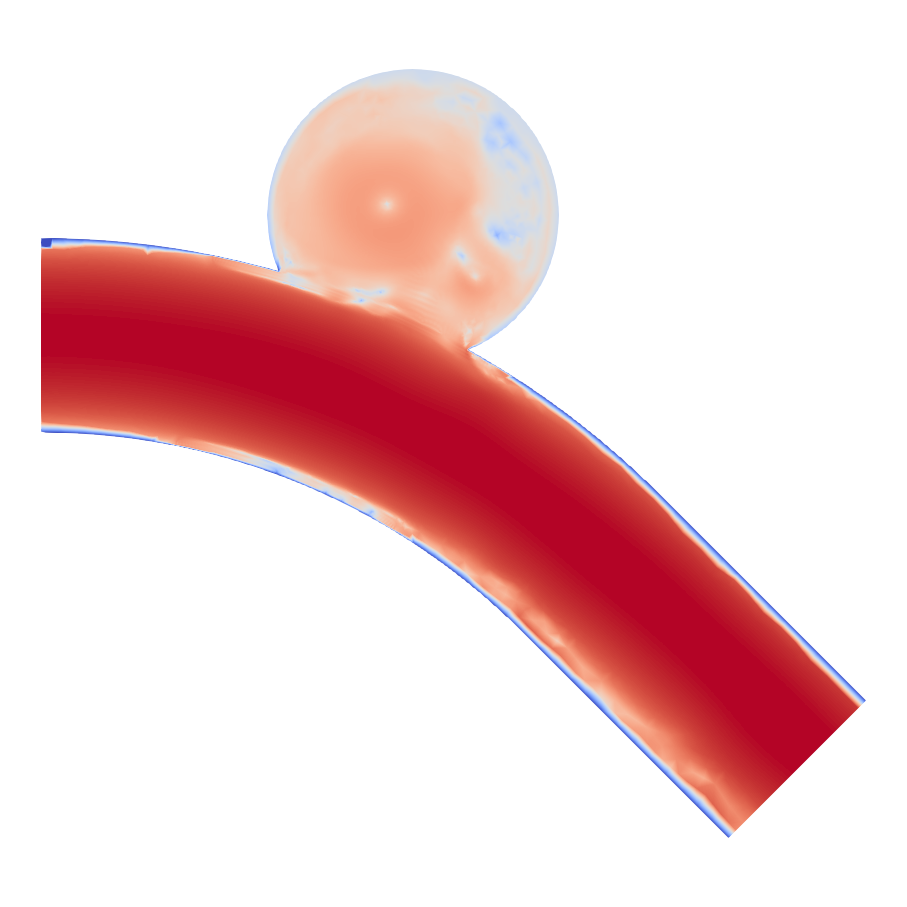}}  \\
    \subfigure[ Velocity solution for $l_0^2=6\cdot 10^{-6} \si{\milli\square\meter}$ ]
    {\label{fig:FDCFDk0.6e-6}\includegraphics[width=7cm,trim=5cm 15cm 5cm 0cm,clip]{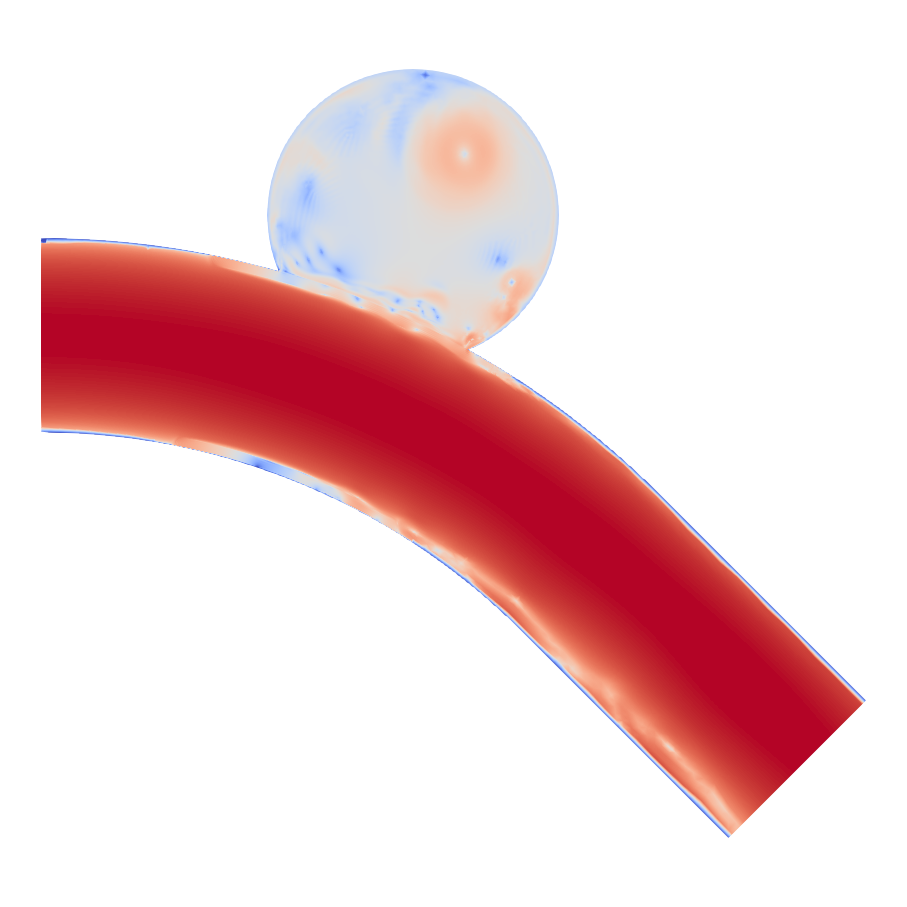}} \hspace*{1cm}
        \subfigure[ Velocity solution for $l_0^2=5.0\cdot 10^{-7} \si{\milli\square\meter}$ ]
    {\label{fig:FDCFDk0.5e-6}\includegraphics[width=7cm,trim=5cm 15cm 5cm 0cm,clip]{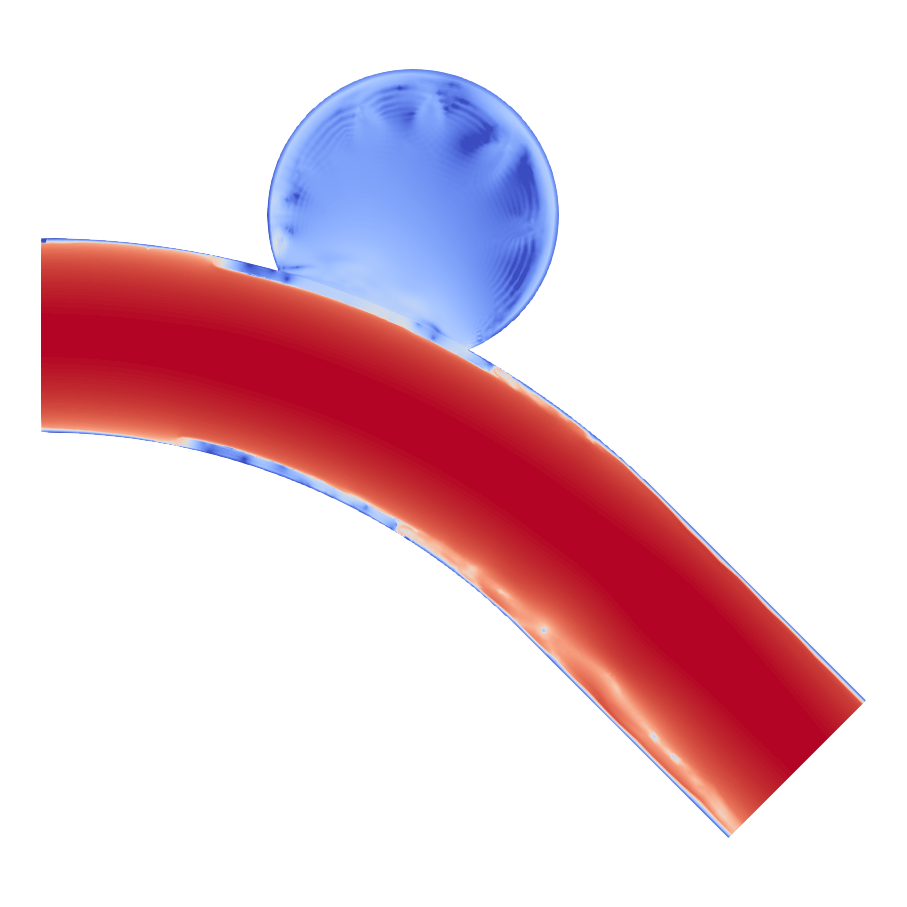}}    
    \caption{Comparison of the fluid flow in the {\CA} at the systolic peak for three different permeabilities with CFD simulation of an untreated {\CA}.
    The different permeability values decrease the magnitude of the flow. For the lowest permeability, the fluid may barely pass the {\FD},
    whereas  fluid flow can still be observed for the other configurations. The color coding uses a logarithmic scale.}
    \label{fig:fdvelocitycomparison}
\end{center}
\end{figure}

Figure~\ref{fig:fdvelocitycomparison} shows a comparison of the fluid velocity fields for the different configurations
with coloring based on a logarithmic scale.
As a reference solution, Figure~\ref{fig:FDCFDReference} depicts a pure CFD solution on a rigid domain of an untreated {\CA}.
While the highest flow velocities naturally occur in the artery~$\Omega_A$,
almost the entire aneurysm cavity~$\Omega_{CA}$ exhibits a flow field with non-negligible velocities.
Due to the inserted {\FD} and changed permeability, the respective magnitude of the fluid flow within the {\CA} diminishes.
Naturally, a less permeable {\FD} (e.g. as in Figure~\ref{fig:FDCFDk0.6e-6}) reduces the inflow in the aneurysm cavity
and, thus, the velocity magnitude inside the cavity more than {\FD} models with larger permeability.
In addition, the center of the vertex can be seen as an indicator for the convective velocity transport.
For Figure~\ref{fig:FDCFDReference}, the center of vertex has already moved further
compared to the treated aneurysm in Figures~\ref{fig:FDCFDk1e-5}, \ref{fig:FDCFDk0.6e-6}, and~\ref{fig:FDCFDk0.5e-6}. 
From the decreasing velocity magnitude within the aneurysm domains of Figures~\ref{fig:FDCFDk1e-5}, \ref{fig:FDCFDk0.6e-6}, and~\ref{fig:FDCFDk0.5e-6},
the proof of concept can be concluded.
The velocity field in Figure~\ref{fig:FDCFDk0.5e-6} exhibits a small magnitude. 
So far, the solution in Figure~\ref{fig:FDCFDk0.5e-6} still shows grid dependencies,
which will be addressed with additional mesh refinement and adaptation in future work.

\begin{figure}
\centering
\includegraphics[width=8cm,trim=5cm 15cm 4cm 0cm,clip]{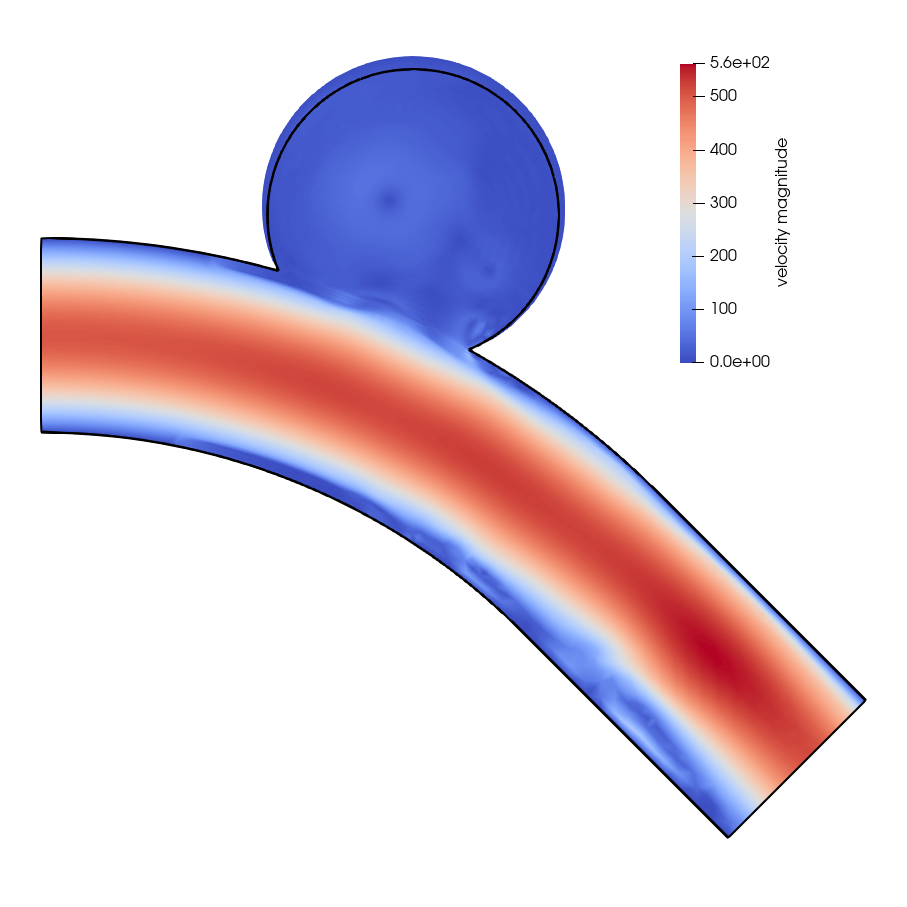}
\caption{Deformed configuration of the {\CA} cavity at the systolic peak with a contour plot of the velocity field:
The black line indicates the undeformed configuration.
The deformation is scaled by a factor of 2 for better visibility.
}
\label{fig:CAwithDeformation}
\end{figure}

To assess the deformation of the {\CA} cavity, Figure~\ref{fig:CAwithDeformation} shows the deformed configuration at the systolic peak.
Naturally, displacements increases towards the tip of the aneurysm.
In addition, the aneurysm is slightly shifted against the fluid flowing passing through the {\FD}.
The area of the two-dimensional artificial aneurysm increases by $\approx 9.5\%$ at the systolic peak compared to the initial state.
Due to the stiffness of the {\FD}, its deformation is small and, thus, cannot be seen in~Figure~\ref{fig:CAwithDeformation} with the naked eye.

\subsection{Impact of coiling on the flow field inside a cerebral aneurysm}
\label{sec:NumExCoiling}

We will compare simulation results for the flow patterns, velocity magnitude as well as the wall shear stress as quantities of interest
for a patient-specific {\CA} geometry before and after a coiling procedure.
While the untreated {\CA} just requires the vessel geometry,
the coiled {\CA} simulation utilizes the virtual coiling models from section \ref{subsec:Coils}
to obtain a coil geometry and use it as a fully resolved obstacle within the aneurysm sac.
We have chosen coiling over stenting or {\web} device treatment for this example,
since the coil model is the most advanced of the three so far, since it is supplemented by a mechanical model as well.
The aneurysm and vessel geometry employed for the experiment is again \#C0074b from~\cite{aneurisk}
with bounding box measures $18.38\,\textup{mm}\times 15.5\,\textup{mm}\times 12.15\,\textup{mm}$, for which the preprocessing steps were already outlined in Figure~\ref{fig:CutNExtendGeometry}.
The final geometry after preprocessing is digitally available in~\cite{Frank2024a}.
In both the empty and the coiled simulation, we study the bloodflow over a time period of roughly 2.7 heartbeats,
ending within the mid of the diastole phase,
where also all subsequent time snapshot images are taken.
The inflow boundary conditions are adopted from~\eqref{eq:GammaPoiseuille} using $\gamma=2$
and the time amplitude profile from Figure~\ref{fig:PulsationProfile},
which we have continued periodically after the first full heartbeat and for stability reasons the values were decreased by 100. At the outlet, for the pressure boundary condition, the constant pressure $p=0$ is applied. For the simulations below, we assume that the fluid is Newtonian which is valid for the flow in the coiled aneurysms\cite{MORALES2013}. 


For discretization, we employ the {\LBM} as introduced in section \ref{subsubsec:FreeFlowNumericalTreatmentLBM}
on a lattice comprising of $N_{\textup{cells}}=\num[group-separator={,}]{15323938}$ cells,
each of size $\Delta x=2.75\times10^{-4}\,\textup{mm}$ and $N_{\textup{d}t}=637113$ time steps with a time step size $\Delta t=4.24\times10^{-6}\,\textup{s}$.
Hence, we run approximately $235967$ time steps for a single heartbeat with a duration of roughly one second.
For the inserted coil,
material and discretization parameters used for its insertion simulation are summarized in Table~\ref{fig:CoilPramameters}.
Snapshots of the virtual coiling procedure are depicted in Figure~\ref{fig:CoilModel}.
Please note again, that the time step size of the coil formation simulation given in Table~\ref{fig:CoilPramameters}
does not have to match the time step size of the fluid simulation,
as they can be regarded as completely independent and sequential computational tasks.
The same holds true for the distance of the material points in the 1D coil vs. the lattice spacing~$\Delta x$ of the fluid simulation.
Considering the coil diameter~$D_2=0.45\,\textup{mm}$ and the given \LBM lattice spacing, the coil's wire is resolved by approximately $16$ cells in diameter, justifying the term ``fully resolved'' for the device simulation.
The simulation was then conducted on an Asus ESC8000A-E11 Server using an AMD Epyc 7713 processor with a core count of $64\times 2.0\,$GHz and 2\,TB of DDR4-RAM ({from which roughly 256\,GB were actually being required).

We now visualize and discuss results at a point in time in the mid of the diastole phase of the third simulated heart beat.
Visualizations have been gerenated using \textit{blender} in combination with the \textit{BVtkNodes}\footnote{BVtkNodes software webpage: \url{https://bvtknodes.readthedocs.io/en/latest/BVTKNodes.html}} add-on~\cite{BVTK, BVTK2}.

First, we study the volumetric cell-wise velocity field~$\vel(\Vec{x})$ for both scenarios, i.e. before and after coiling.
Using \textit{ParaView}\footnote{ParaView software webpage: \url{https://www.paraview.org/}},
we compute streamlines in each case and use them to compare the respective flow patterns in Figure~\ref{fig:Streamlines}.
\begin{figure}
\begin{center}
    \subfigure[Streamline pattern within the complete computational domain, including the \textit{un}treated aneurysm and the adjacent vessel where a pipe-flow pattern develops, while a vortex is formed inside the aneurysm.]{\includegraphics[width=10cm]{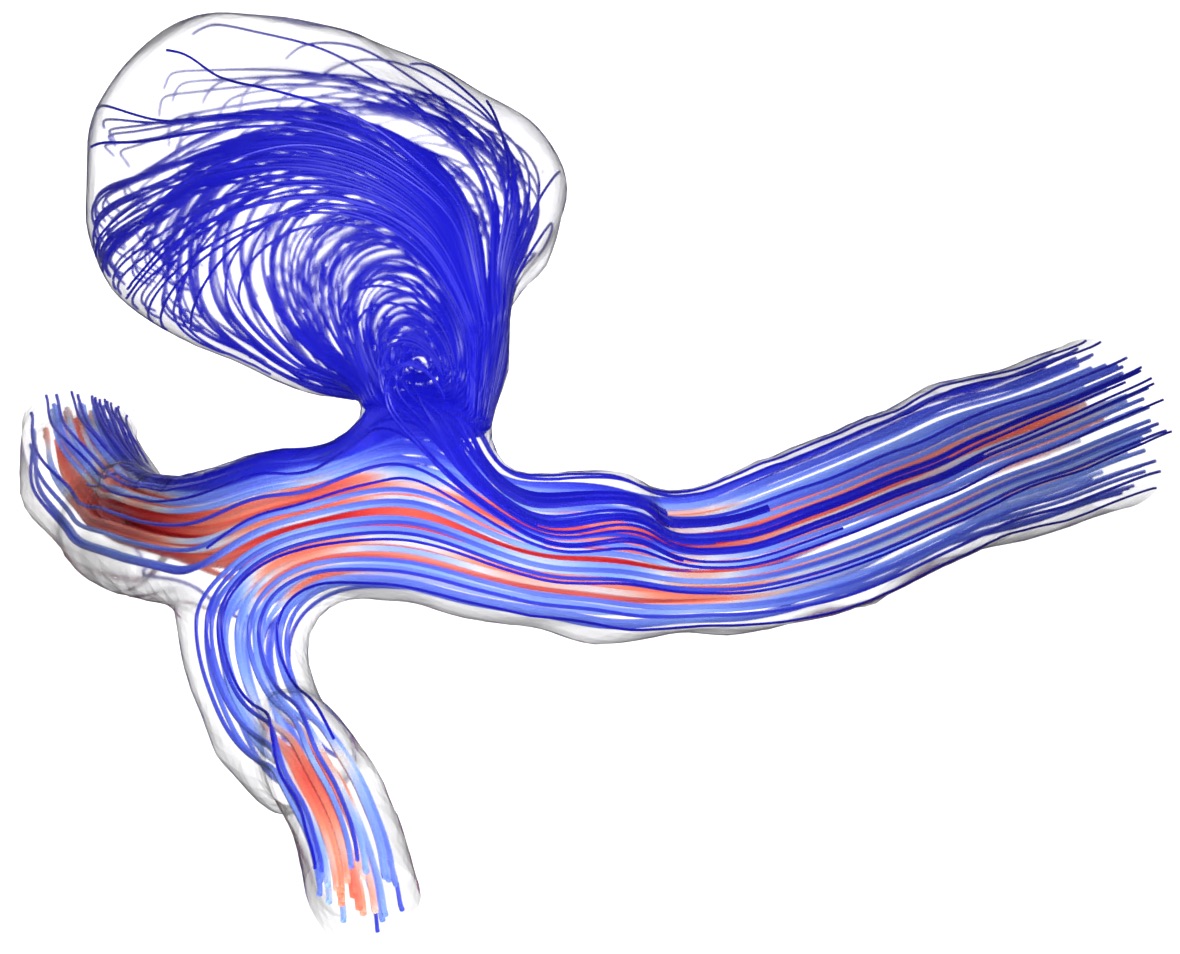}\hspace*{1cm}\raisebox{1cm}{\includegraphics[width=2.75cm]{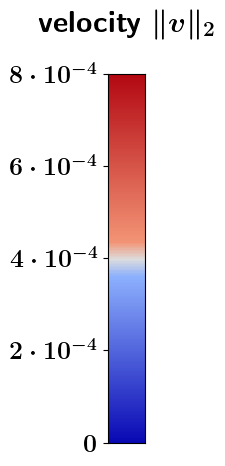}}}\\
    \subfigure[Zoom into the untreated aneurysm part of the domain to highlight the vortex pattern inside the {\CA}.]{\includegraphics[width=4cm]{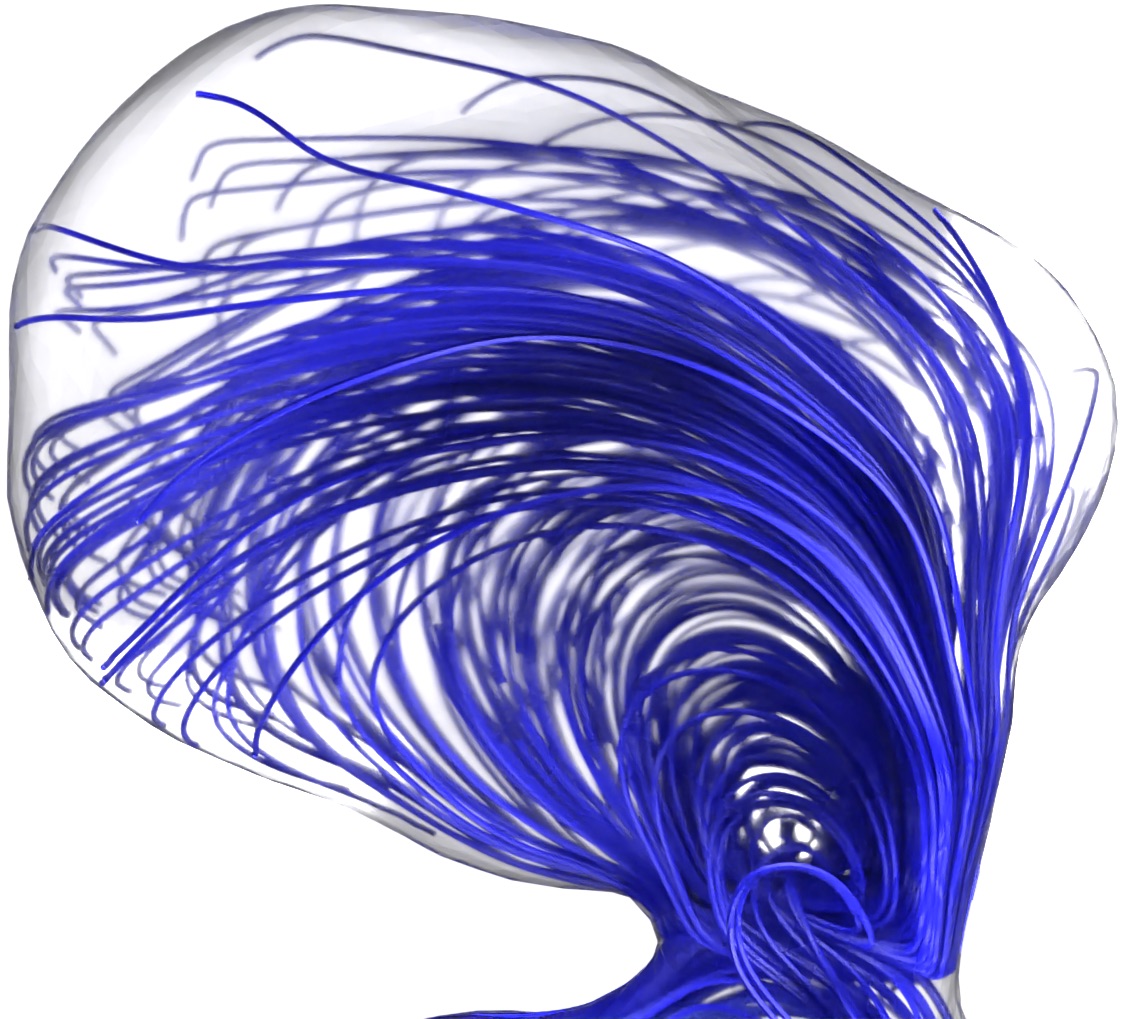}}\hspace*{1cm}\subfigure[Streamlines within the \textit{coiled} aneurysm, i.e. two times the same situation where on the right the coil is just made invisible.]{\includegraphics[width=4cm]{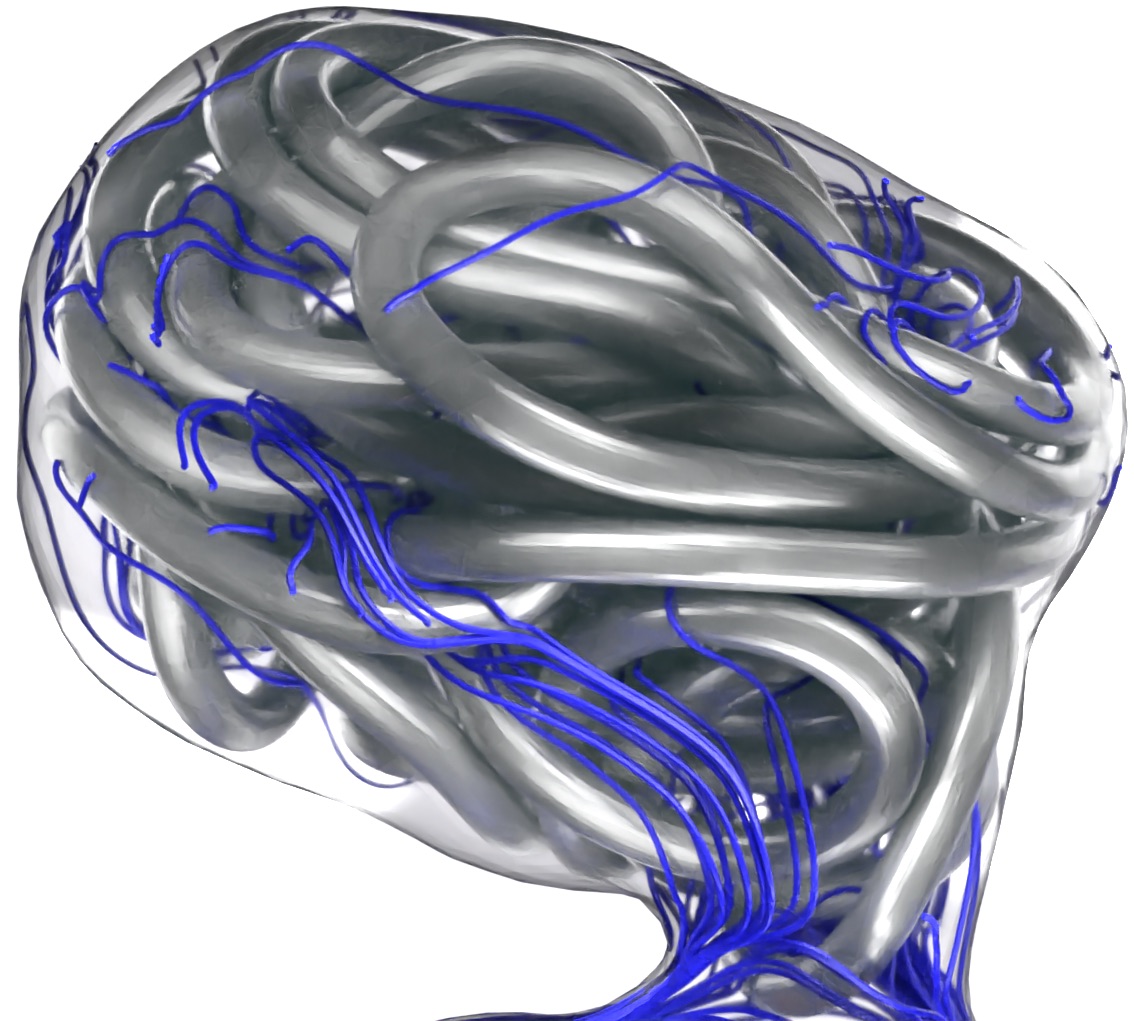}\includegraphics[width=4cm]{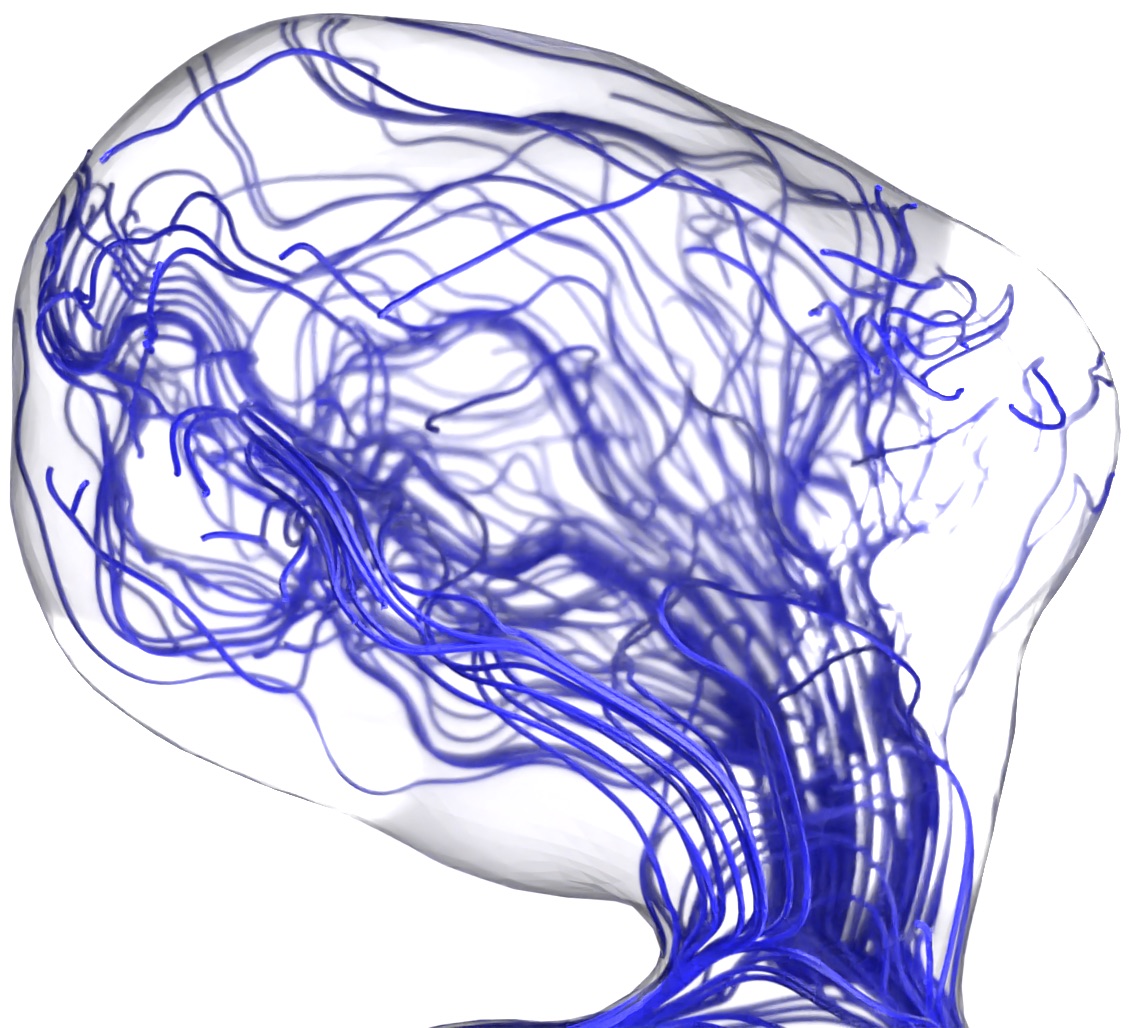}}
    \caption{Comparison of flow patterns between untreated and coiled aneurysm. The flow behavior within the aneurysm is changed by the presence of the coil from a developed vortex to only creeping flow in between the coil's windings.\label{fig:Streamlines}}
\end{center}
\end{figure}
The presence of the coil reduces the flow field to a creeping flow through the remaining tight spaces within the coil.
Similarly, Figure~\ref{fig:Crosssection} compares cross sections through the velocity magnitude field~$\|\vel(\vec{x})\|_2$
to obtain an impression of the change in the distribution of high flow velocities within the aneurysm due to the presence of the coil.
\begin{figure}
\begin{center}
    \subfigure[Cross section through the velocity field in the \textit{un}treated situation following the main vessel with linearly scaled colorbar highlighting the pipe flow character within the vessel while within the aneurysm the velocity magnitude is comparably low.]{\includegraphics[width=10cm]{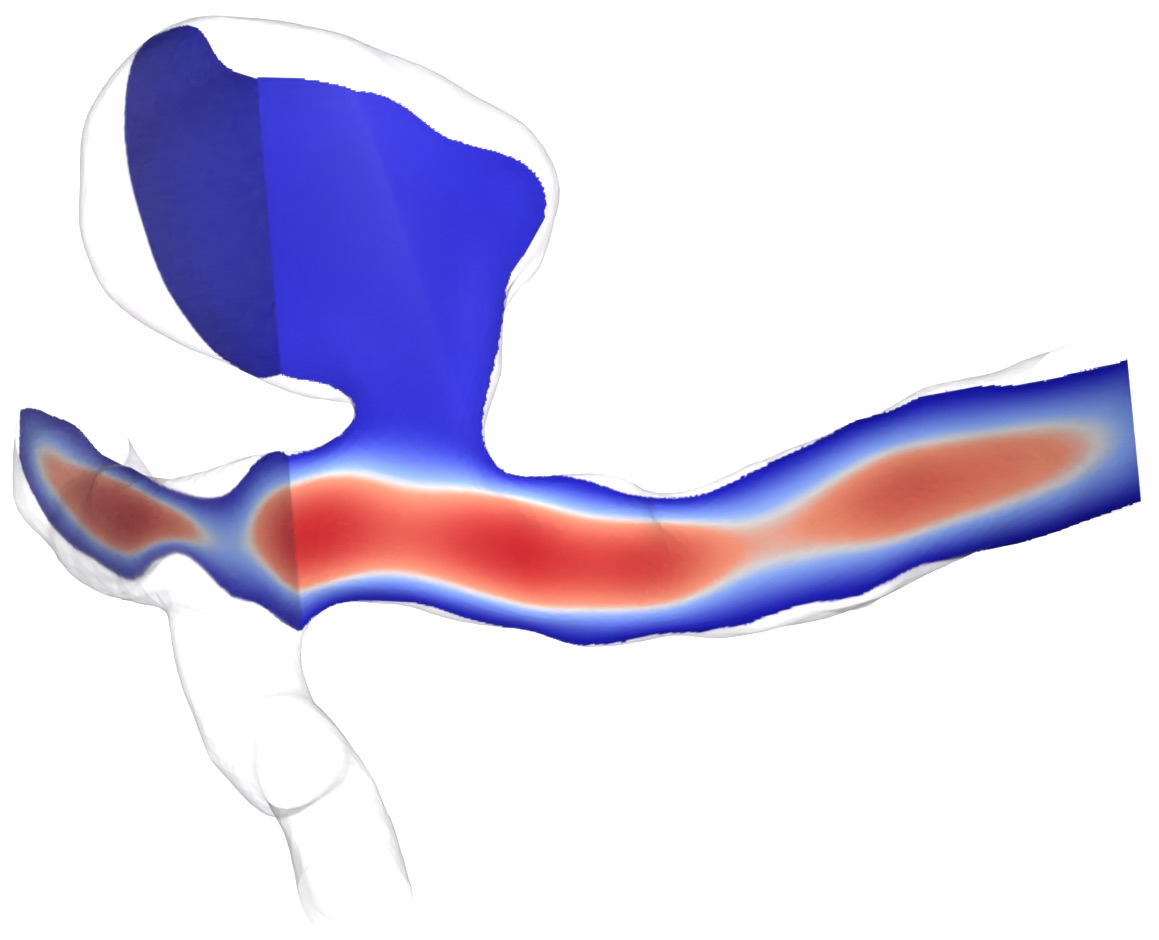}\hspace*{1cm}\raisebox{1cm}{\includegraphics[width=2.75cm]{images/CBAR_Velocity.png}}}\\
    \hspace*{1cm}\subfigure[Zoom into the untreated aneurysm with logarithmic color scaling.]{\includegraphics[width=4cm]{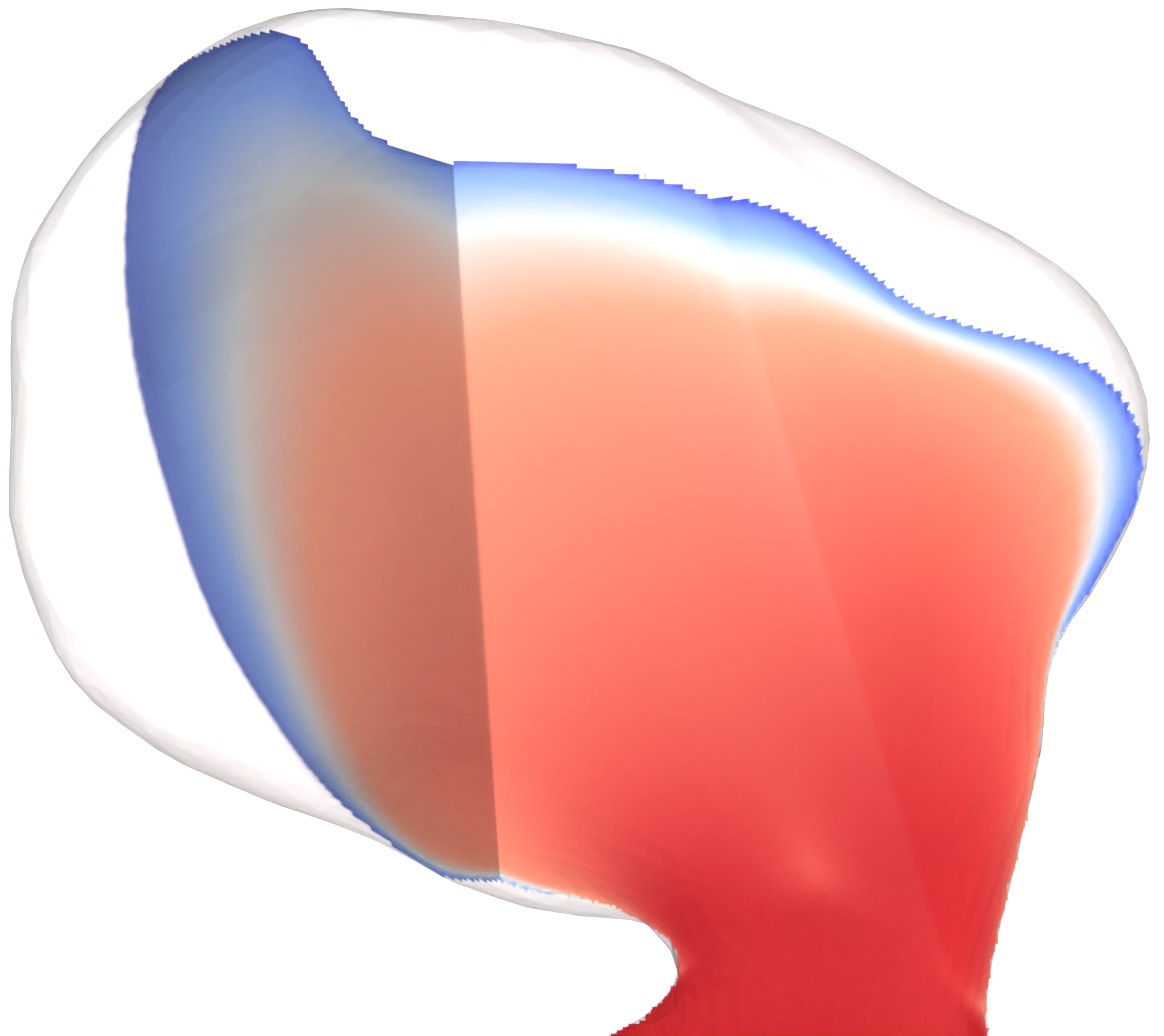}}\hspace*{1cm}\subfigure[Zoom into the untreated aneurysm with logarithmic color scaling. Both images show the same situation where the coil has been made invisible in the right picture.]{\includegraphics[width=4cm]{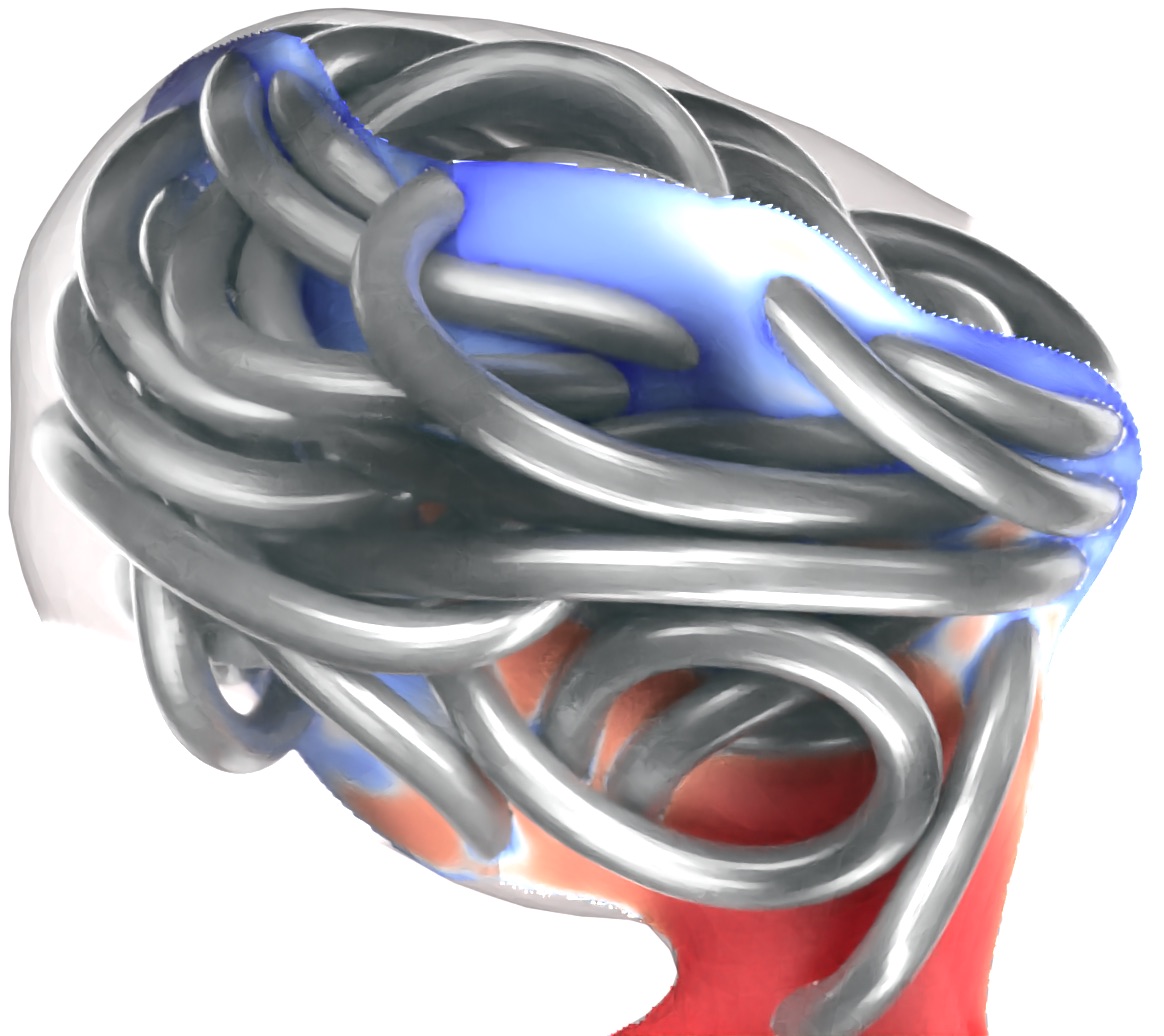}\includegraphics[width=4cm]{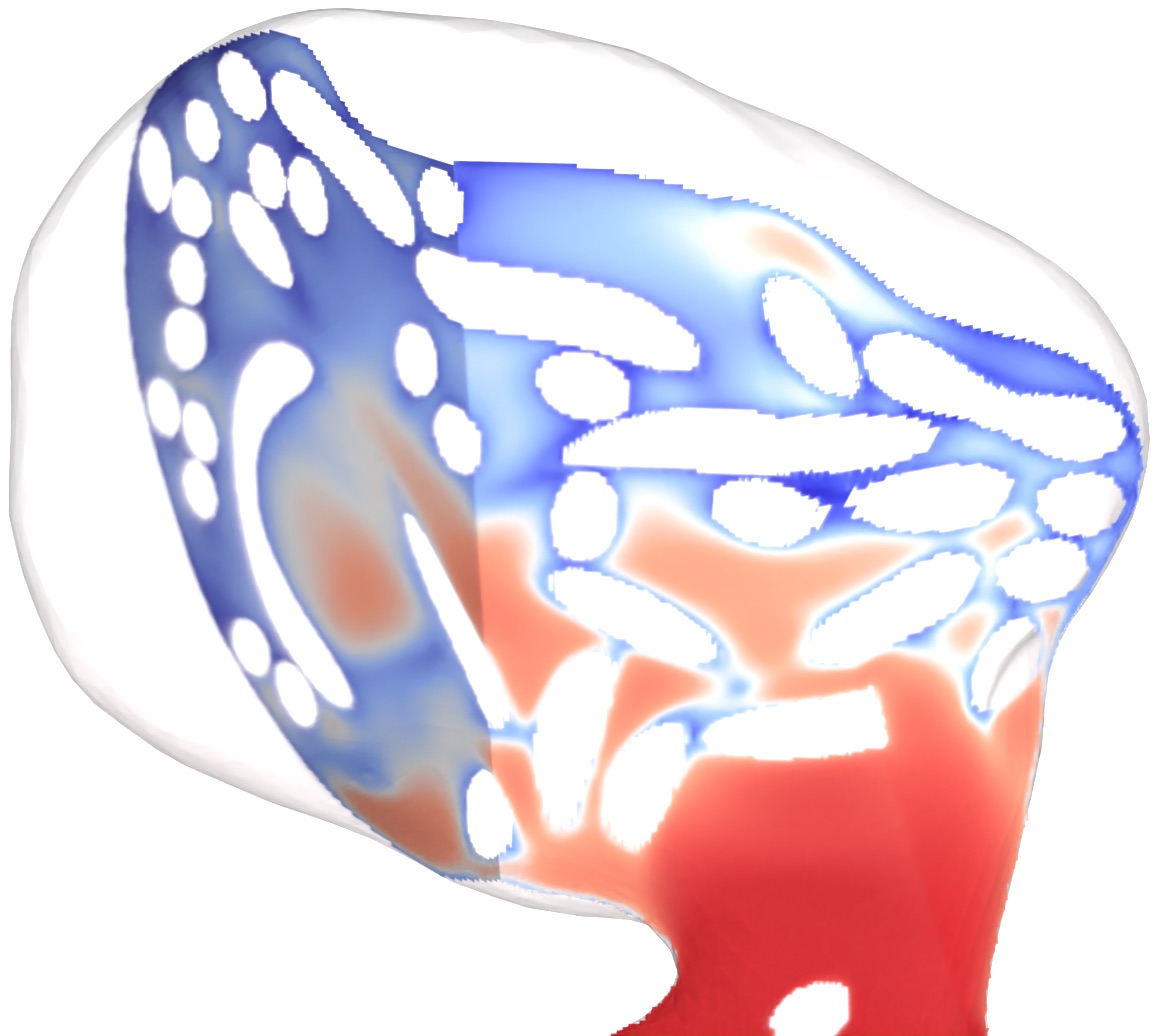}}\hspace*{-10mm}\includegraphics[width=4.5cm]{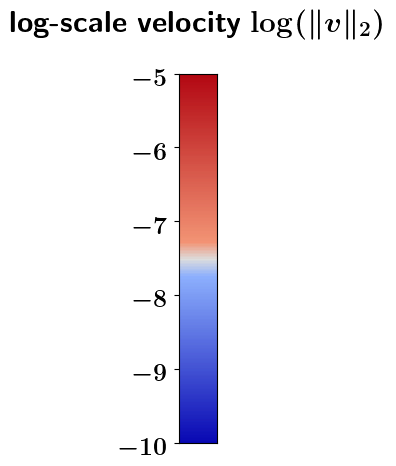}
    \caption{Comparison of velocity field magnitudes on domain cross sections. In order to highlight the structure of the comparably low velocity magnitudes within the aneurysm in subfigures (b) and (c), we switch to a logarithmic scale. Due to the presence of the coil, the regions of larger velocity magnitude are driven towards the neck, respectively out of the aneurysm cavity reducing the overall flow and hence perfusion of the aneurysm.\label{fig:Crosssection}}
\end{center}
\end{figure}
Effectively, the presence of the coil slows down the flow field by at least one or more orders of magnitude.
Finally in Figure~\ref{fig:WSS}, we evaluate the wall shear stress, which is only well defined on the surface of the aneurysm, 
to assess the success of the coiling procedure to reduce said quantity in order to alleviate the aneurysm's rupture risk.
\begin{figure}
\begin{center}
    \subfigure[Wall shear stress distribution on the aneurysm's surface \textit{before} coiling treatment.]{\includegraphics[width=4cm]{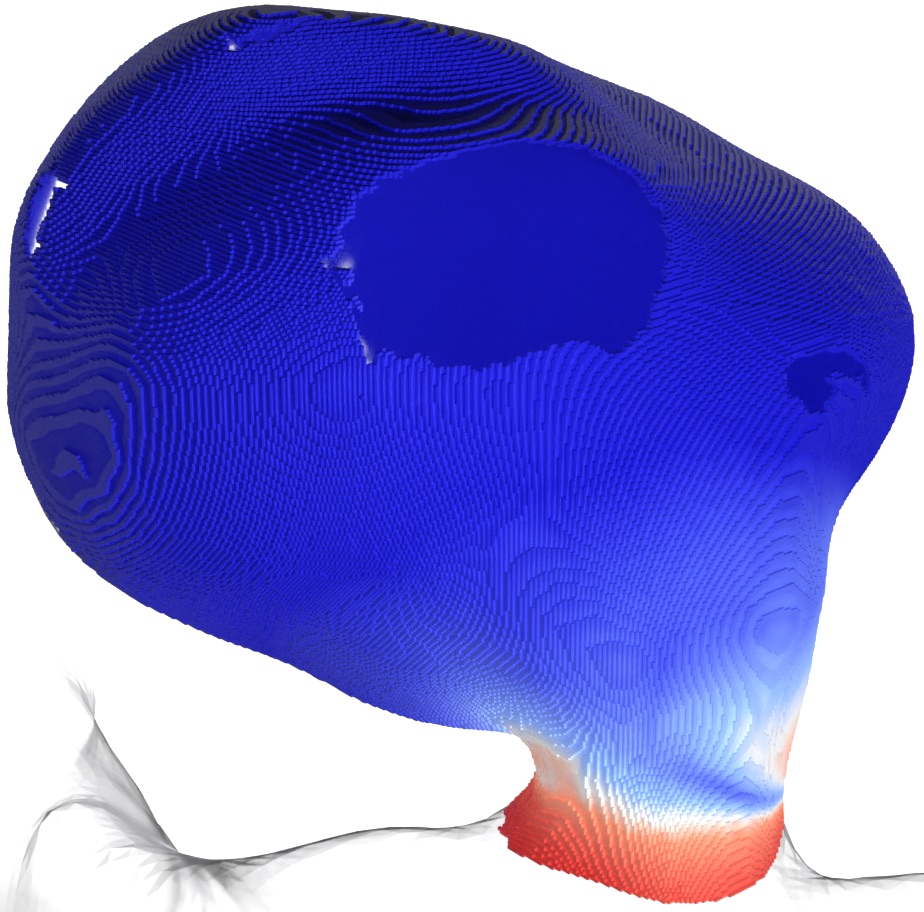}}\hspace*{1cm}\subfigure[{\WSS} distribution in the \textit{coiled} aneurysm.]{\includegraphics[width=4cm]{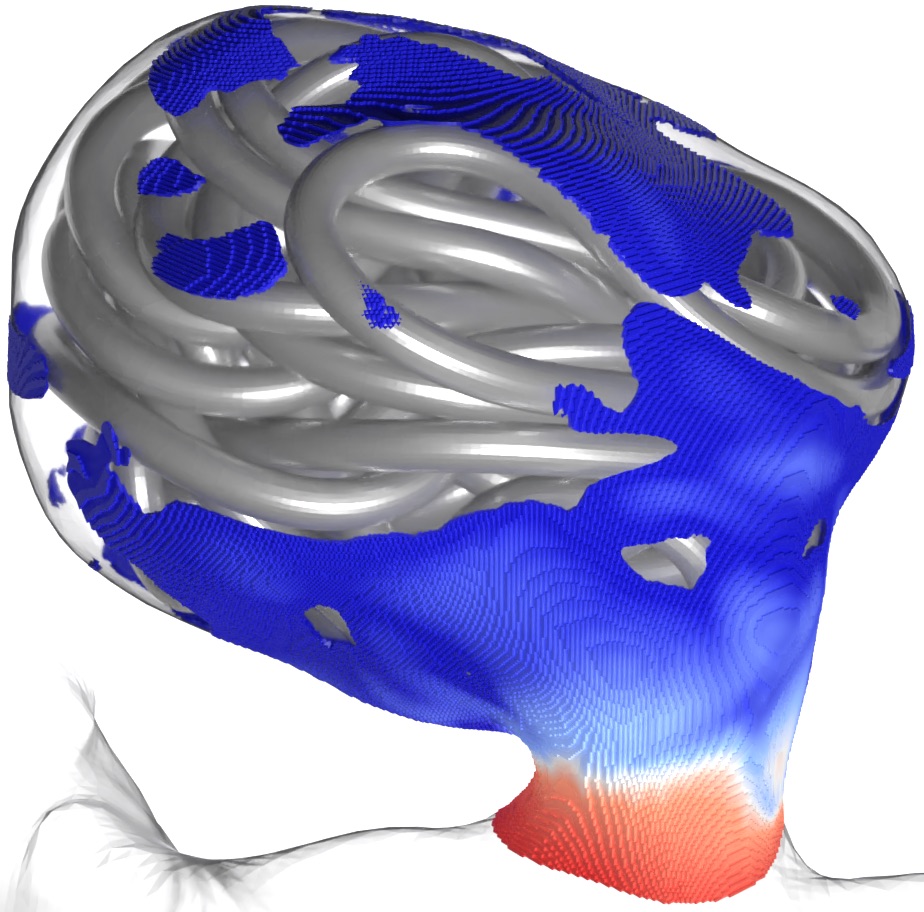}}\hspace*{1cm}\includegraphics[width=1.75cm]{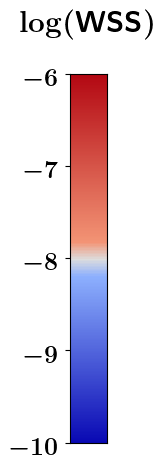}
    \caption{Comparison of wall shear stress distributions on the aneurysm's surface due to the coiling treatment. To increase the visibility of the wall shear stress dynamics, we have switched to a logarithmic scale with values lying below the lower end of the colorbar being transparent.}
    \label{fig:WSS}
\end{center}
\end{figure}
In line with the reduction of the flow velocities,
shear rates and, thus, also wall shear stresses are greatly reduced by the presence of the coil. Since the wall shear stress computation is done in the LBM setting, some artifacts are present because of the corresponding discretization. To get rid of them, the Gaussian smoothing operator is applied to the wall shear stress field after the simulation is finished.


Overall, the numerical simulations are capable to reproduce the intended effect of coil insertion,
namely the reduction of wall shear stresses to reduce the rupture risk as well as reductions in flow velocities to aid thrombus formation inside the {\CA}.
Another account on this topic can again be found in~\cite{horvat2024lattice}.

\section{Concluding remarks and outlook}
\label{sec:Conclusion}

In this paper,
we have outlined various components to build mathematical models and numerical simulation workflows
with the goal to assess different endovascular treatment options for cerebral aneurysms
and predict their long-term success for enhanced patient well-being.

Since we target patient-specific predictions,
we first have described the necessary steps from medical imaging data to ready-to-use representations of computational grids.
On these grids, we can either use a {\LBM} method to examine the flow field within a pure CFD simulation of cerebral aneurysms,
which we have demonstrated for the case of endovascular coiling,
or employ a finite element discretization of a poro-elastic continuum model to study the impact of endovascular devices on the flow field within a deformable aneurysm domain in a homogenized manner,
showcased for a simplified aneurysm treated by a flow diverter.
Then, we have proposed models for three types of endovascular devices,
in particular a mechanical model for endovascular coils
and two purely geometric models for {\web} devices and stents/flow diverters, respectively.

The presented one domain modeling approach based on a poro-elastic medium
appears as a simple, but promising approach to model the fluid and 
structural behavior of a {\FD} in a homogenized fashion. 
However, before clinically relevant questions can be answered 
with the suggested modeling approach, a few considerations must still be addressed.
The estimation of a proper permeability is still a topic of research, 
since it may also include anisotropic effects of the fluid flow resulting from the {\FD}'s geometry. Additionally, the effective stiffness of a {\FD} should be taken into 
consideration within the modeling approach. 
For the considered material parameters, the poro-elastic medium behaved rather stiff, 
resulting in a relatively low deformation of the {\FD}. 
This raises the question of how to assess the bending of the thin wires 
of a {\FD} due to the interaction of blood and artery wall and include the 
exact structural behavior into the modeling approach.

The {\LBM} solver can easily be used to study flow fields in patient-specific geometries,
yet limited to rigid vessel geometries so far.
In future work, we plan to include moving boundaries into the {\LBM} framework in order to account for the vessel's deformation due to the pulsatile nature of the heart beat.

With the presented models, we are able to compute flow fields in cerebral aneurysms just after an endovascular intervention.
For now, we had to limit ourselves to two-dimensional problems for the {\PM} model.
Extension to three-dimensional geometries is straightforward from a modeling point of view,
however requires to tackle performance questions,
for example the design of effective preconditioners for the system of linear equations arising in each nonlinear iteration of the {\PM} solver.
Having successfully applied the concept of physics-based block preconditioning to surface-coupled multi-physics couplings
such as contact mechanics~\cite{Wiesner2018a,Wiesner2021a} or fluid-solid interaction~\cite{Mayr2020a},
we intend to apply this principle also to the present case of poro-elastic media.
Furthermore, we are working on the fully-coupled embedding of patient-specific {\CA} models using the proposed {\PM} approach
into reduced-dimensional representations of the cerebral blood flow (based on ideas from our prior work~\cite{fritz20221d}).
The {\LBM} models for coiled {\CAs} are currently limited to fixed domains and require an extension to domains with moving boundaries
to also account for vessel dilation due to the blood pulsation.
In addition, increased voxel resolutions are required to facilitate a more accurate modeling of the thin coiling wires,
which in turns requires to ramp up the use of computational resources.
While the developers of the {\LBM} solver \textit{waLBerla} have demonstrated its application to very fine meshes,
this computational power needs to be transferred to the application at hand.

The present methods and results will serve as a starting point to simulate the subsequent thrombus formation, which is driven by biochemical reaction processes,
but dominantly governed by the blood flow velocity inside the aneurysm cavity.
Ultimately, we hope to cover the entire spectrum of flow field prediction, thrombus formation, and assessment of the quality and stability of the aneurysm occlusion, such that we ultimately arrive at our goal to assist the attending neuroradiologist in management and treatment of cerebral aneurysms.

\section*{Acknowledgement}
The work in this manuscript has been carried out within the priority programme ``SPP 2311:  Robust coupling of continuum-biomechanical in silico models to establish active biological system models for later use in clinical applications -- Co-design of modeling, numerics and usability''
and has been funded by the Deutsche Forschungsgemeinschaft (DFG, German Research Foundation) -- 465242983 (KI 2101/8-1, PO 1883/6-1, WO 671/20-1).
Fabian Holzberger, Medeea Horvat, Markus Muhr, Natalia Nebulishvili and Barbara Wohlmuth also gratefully acknowledge the financial support
partially provided by the Deutsche Forschungsgemeinschaft (DFG, German Research Foundation) under the grant number WO 671/11-1.
Martin Frank, Matthias Mayr and Alexander Popp gratefully acknowledge the computing resources provided by the Data Science~\& Computing Lab at the University of the Bundeswehr Munich.

\section*{Conflict of interest}
On behalf of all authors, the corresponding author states that there is no conflict of interest.



\bibliographystyle{model1-num-names} 
\bibliography{references}




\end{document}